 \documentclass[letterpaper,10pt]{article}

\usepackage{amsfonts,amsmath,verbatim,graphicx, subfigure,longtable}
\usepackage{amssymb,amsthm,amscd}
\usepackage{pstricks}
\usepackage{epstopdf}
\usepackage{color}
\usepackage{transparent}
\usepackage{tikz}
\usetikzlibrary{fit, positioning}

\def\wh{\widehat}

\newtheorem{thm}{Theorem}[section]
\newtheorem{cor}[thm]{Corollary}
\newtheorem{prop}[thm]{Proposition}
\newtheorem{lem}[thm]{Lemma}
\newtheorem{definition}[thm]{Definition}
\newtheorem{assumption}[thm]{Assumption}

\newtheorem{remark}[thm]{Remark}

\newenvironment{pf}{\noindent\emph{Proof \,}}{\mbox{}\qed}

\newcommand{\toL}{\,{\buildrel \mathcal{L} \over \longrightarrow}\,}

\numberwithin{equation}{section}
\def\R{{\mathbb R}}     
\def\Z{{\mathbb Z}}     
\def\P{{\mathbb P}}     
\def\E{{\mathbb E}}     
\def\D{{\mathcal E}}    
\def\F{{\mathcal F}}    
\def\X{{\mathfrak X}}   

\def\A{{\mathcal A}}   
\def\L{{\mathfrak L}}   

\def\eps{\varepsilon}
\def\<{{\langle}}
\def\>{{\rangle}}
\def\1{{\bf 1}}         

\renewcommand{\bar}{\overline}

\begin{document}
\allowdisplaybreaks

\title{\Large \bf
Hydrodynamic Limits and Propagation of Chaos for Interacting  Random Walks in
Domains
 }

\author{{\bf Zhen-Qing Chen} \quad and \quad {\bf Wai-Tong (Louis) Fan}}
\date{\today}
\maketitle

\begin{abstract}
A  new non-conservative stochastic reaction-diffusion system in which two families of random walks in two adjacent domains interact near the interface is introduced and studied in this paper.
Such a system  can be used to model the transport of positive and negative charges in a solar cell or the population dynamics of two segregated species under competition. We show that in the macroscopic limit, the particle densities converge to the solution of a coupled nonlinear heat equations.
For this, we first prove that propagation of chaos holds by establishing the uniqueness of a new BBGKY hierarchy. A local central limit theorem for reflected diffusions in bounded Lipschitz domains is also established as a crucial tool.
\end{abstract}

\tableofcontents

\section{Introduction}

\begin{figure}[h]
	\begin{center}
    \vspace{-0.5em}
	\includegraphics[scale=.22]{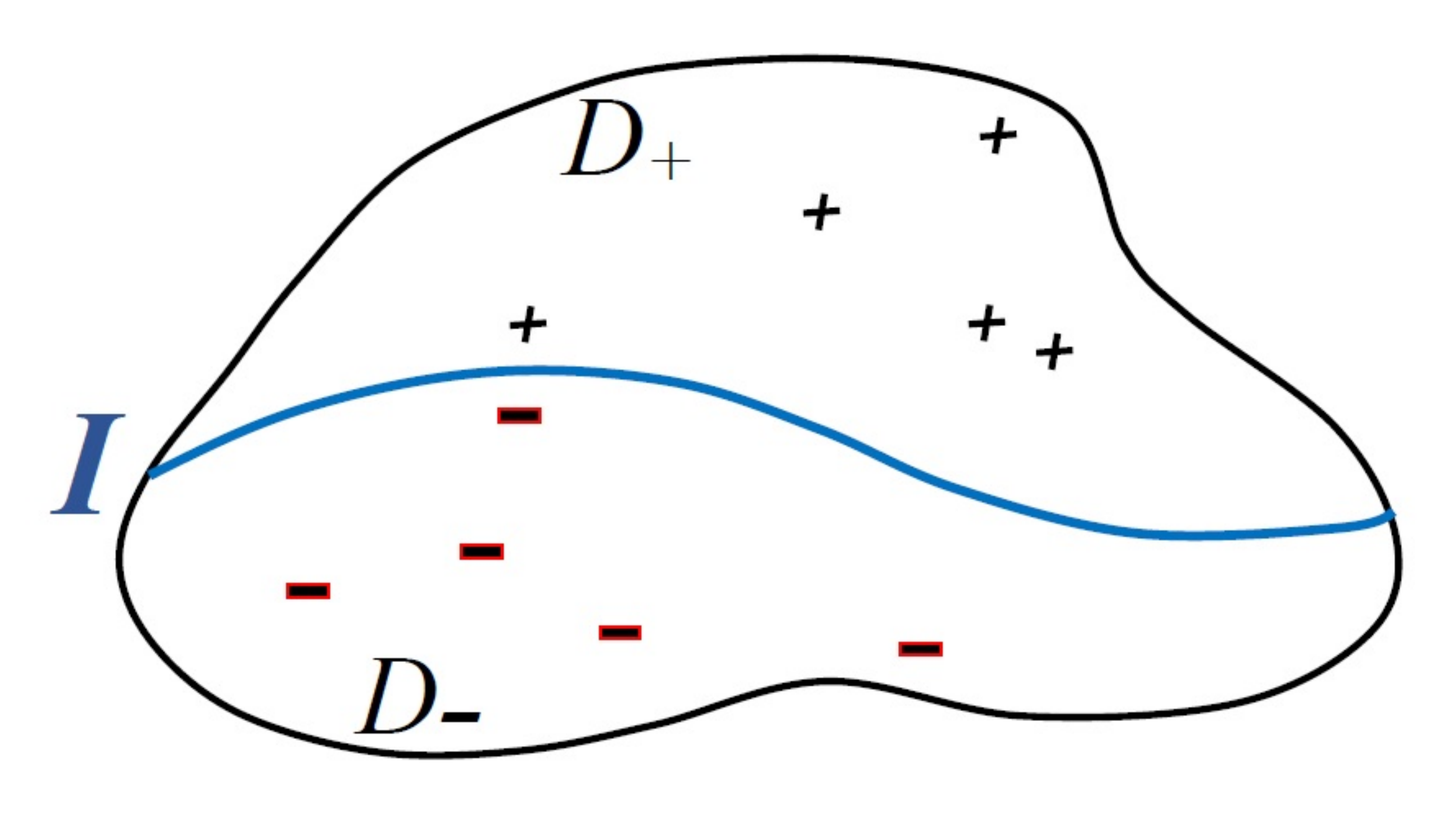}
    \vspace{-1em}
    \caption{$I$ is the interface of $D_{+}$ and $D_-$}\label{fig:Interface}
    \vspace{-1em}
	\end{center}
	\end{figure}

The original motivation of this project is to study the transports of positive and negative charges in solar cells.
We model a solar cell by a domain in $\R^d$ that is divided into two disjoint sub-domains $D_+$ and $D_-$ by an interface $I$,
 a $(d-1)$-dimensional hypersurface, which can be possibly disconnected.
$D_+$ and $D_-$  represent the hybrid medium that confine the positive and the negative charges, respectively.
At microscopic level, positive and negative charges  are initially modeled by $N$ independent reflected Brownian motion (RBM) with drift on $D_+$ and on $D_-$ respectively.
(In this paper, they are actually modeled by $N$ independent random walks on lattices inside $D_+$ and $D_-$ that serve as discrete approximation of RBM with drifts.)
These random motions model the transport of positive (respectively negative) charges under an electric potential (See Figure \ref{fig:Interface}).

These two types of particles annihilate each other at a certain rate when they come close to each other near the interface $I$. This interaction models the annihilation, trapping, recombination and separation phenomena of the charges. The interaction distance is of microscopic order $\eps$ where $N\eps^d$ is comparable to 1, and the intensity of annihilation per pair is of order $\lambda/\eps$ where $\lambda\geq 0$ is a given parameter. This means that, intuitively and roughly speaking, according to a random time clock which runs with a speed proportional to the number of pairs (one particle of each type) of distance $\eps$, we annihilate a pair (picked uniformly among those pairs of distance less than $\eps$) with an exponential rate of parameter $\lambda/\eps$. The above scaling guarantees that in the limit, a nontrivial proportion of particles are annihilated in any open time interval. We investigate the scaling limit of the empirical distribution of positive and negative charges;
that is, the hydrodynamic limit of this interacting diffusion system. We show that in the macroscopic level, the empirical distribution converges to a deterministic
measure whose density satisfies a  system of partial differential equations that has non-linear interaction
at the interface.

The study of \emph{hydrodynamic limits} of particle systems with stochastic dynamics is of fundamental importance in many areas. This study dates back to the sixth Hilbert problem in year 1900, which concerns the mathematical treatment of the axioms of physics, and to Boltzmann's work on principles of mechanics. Proving hydrodynamic limits corresponds to establishing \emph{the law of large number} for the empirical measure of some attributes (such as position, genetic type, spin type, etc.) of the individuals in the systems. It contributes to our better understanding of the asymptotic behavior of many phenomena, such as chemical reactions \cite{tK71}, population dynamics \cite{rDsL94, tK81, rMmN92}, super-conductivity \cite{htY91}, quantum dynamics \cite{lEbShtY07}, fluid dynamics \cite{fG05}, etc. It reveals fascinating connections between the microscopic stochastic systems and deterministic partial differential equations that describe the macroscopic pictures. It also provides approximations via stochastic models to some partial differential equations that are hard or impossible to solve directly.

Since the work of Boltzmann and Hilbert, there have been many different lines of research on stochastic particle systems. Various models were constructed and different techniques were developed to establish hydrodynamic limits. Among those techniques, the entropy method \cite{GPV88} and the relative entropy method \cite{htY91} are considered to be general methods. Unfortunately these methods do not seem to work for our model due to the singular interaction near the interface.

Many models studied in literature are conservative, for example exclusion processes \cite{cKcL98, KOV89} and Fleming-Viot type systems \cite{BHM00, BQ06}.  Reaction-diffusion systems (R-D system) constitute a class of models that  are typically non-conservative. These are systems which have hydrodynamic limits of the form $\frac{\partial u}{\partial t}= \frac{1}{2}\Delta u + R(u)$ (a reaction-diffusion equation, or R-D equation in short), where $R(u)$ is a function in $u$ which is thought of as the reaction term. R-D systems arise from many different contexts and have been studied by many authors. For instance, for the case $R(u)$ is a polynomial in $u$, these systems contain the Schl\"ogl's model and were studied in \cite{pD88a, pD88b, pK86, pK88} on a cube with Neumann boundary conditions, and in \cite{BDP92, BDPP87} on a periodic lattice. Recently, perturbations of the voter models which contain the Lotka-Volterra systems are considered in \cite{CDP11}. In addition to results on hydrodynamic limit, \cite{CDP11} also established general conditions for the existence of non-trivial stationary measures and for extinction of the particles.

Our model is a non-conservative stochastic particle system which consists of two types of particles. In \cite{BQ06}, Burdzy and Quastel studied an annihilating-branching system of two types of particles, for which the total number of particles of each type remains constant over the time. Its hydrodynamic limit is described by a linear heat equation with zero average temperature.
In contrast, besides being non-conservative, our model gives rise to a  system of nonlinear differential equations that seems to be new.
Moreover, the interaction between two types of particles is singular near the interface of the two media, which gives rise to a boundary integral term in the hydrodynamic limit. The approach of this paper provides some new tools that are potentially useful for the study of other non-equilibrium systems.

We now give some more details on the discrete approximation of the spatial motions in our modeling.
We approximate $D_{\pm}$ by square lattices $D^{\eps}_{\pm}$ of side length $\eps$, and then approximate reflected diffusions on $D_{\pm}$ by continuous time random walks (CTRWs) on $D^{\eps}_{\pm}$. The rigorous formulation of the particle system is captured by the operator $\mathfrak{L}^{\eps}$ in (\ref{E:generator0}).

Let $X^{\pm}_{i}(t)$ be the position of the particle with index $i$ in $\bar{D}_{\pm}$ at time $t$. We prescribe each particle a mass $1/N$ and consider the normalized empirical measures
\begin{eqnarray*}
\X^{N,+}_t(dx) := \frac{1}{N}\sum_{\alpha: \alpha\sim t}\1_{X^{+}_{\alpha}(t)}(dx) \quad \text{and} \quad
\X^{N,-}_t(dy) := \frac{1}{N}\sum_{\beta: \beta\sim t}\1_{X^{-}_{\beta}(t)}(dy).
\end{eqnarray*}
Here $\1_y(dx)$ stands for the Dirac measure concentrated at the point $y$,  while
$\alpha\sim t$ if and only if the particle $X^{+}_{\alpha}$ is alive at time $t$, and $\beta\sim t$ if and only if the particle $X^{-}_{\beta}$ is alive at time $t$. For fixed positive integer $N$ and $t>0$, $\X^{N,\pm}_t$ is a random measure on $\bar{D}_{\pm}$. We want to study the asymptotic behavior, when $N\to\infty$ (or equivalently $\eps \to 0$), of the evolution in time $t$ of the pair $(\X^{N,+}_t,\,\X^{N,-}_t)$.

\subsection{Main results}

Our first main result (Theorem \ref{T:conjecture}) implies the following.  Suppose each particle in $D_{\pm}$ is approximating a RBM with gradient drift $\frac{1}{2}\,\nabla (\log \rho_{\pm})$, where $\rho_{\pm}\in C^1(\bar{D}_{\pm})$ is strictly positive. Then under appropriate assumptions on the initial configuration $(\X^{N,+}_0,\,\X^{N,-}_0)$, the normalized empirical  measure $(\X^{N,+}_t(dx),\,\X^{N,-}_t(dy))$ converges in distribution to a deterministic measure
$$
(u_+(t,x)\rho_+(x)dx,\,u_-(t,y)\rho_-(y)dy)
$$
for all   $t>0$, where $(u_+,u_-)$ is the solution of the following coupled heat equations:
    \begin{equation}\label{E:coupledpde:+}
        \left\{\begin{aligned}
        \dfrac{\partial u_+}{\partial t}           &= \frac{1}{2}\Delta u_+ +\frac{1}{2}\nabla (\log \rho_{+})\cdot\nabla u_+    & &\qquad\text{on } (0,\infty)\times D_+  \\
        \dfrac{\partial u_+}{\partial \vec{n}_+} &=\frac{\lambda}{\rho_+}\,u_+u_-\,\1_{I}  & &\qquad\text{on }  (0,\infty)\times \partial D_+
        \end{aligned}\right.
    \end{equation}
and
    \begin{equation}\label{E:coupledpde:-}
        \left\{\begin{aligned}
        \dfrac{\partial u_-}{\partial t}           &= \frac{1}{2}\Delta u_- +\frac{1}{2}\nabla (\log \rho_{-})\cdot\nabla u_-  & &\qquad\text{on } (0,\infty)\times D_-  \\
        \dfrac{\partial u_-}{\partial \vec{n}_-} &=\frac{\lambda}{\rho_-}\,u_+u_-\,\1_I  & &\qquad\text{on }  (0,\infty)\times \partial D_- ,
        \end{aligned}\right.
    \end{equation}
where $\vec{n}_\pm$ is the inward unit normal vector field on $\partial D_{\pm}$ of $D_{\pm}$ and $\1_I$ is the indicator function on $I$. Note that $\rho_{\pm}=1$ corresponds to the particular case when there is no drift.

The above result tells us that for any fixed time $t>0$, the probability distribution of a randomly picked  particle in $D^{\eps}_{\pm}$ at time $t$ is close to $c_{\pm}(t)u_{\pm}(t,x)$ when $N$ is large, where $c_{\pm}(t)=(\int_{D_{\pm}} u_{\pm}(t))^{-1}$ is a normalizing constant. In fact, the above convergence holds at the level of the path space. That is, the full trajectory (and hence the joint law at different times) of the particle profile converges to the deterministic scaling limit described by (\ref{E:coupledpde:+}) and (\ref{E:coupledpde:-}), not only its distribution at a given time.

\noindent {\bf Question}:  How about the limiting joint distribution of more than one particles ?

Our second main result (Theorem \ref{T:correlation}) answers this question. It asserts that \textbf{propagation of chaos} holds true for our system; that is, when the number of particles tends to infinity, their positions appear to be independent of each other.  More precisely, suppose $n$ and $m$ unlabeled particles in $D^{\eps}_+$ and $D^{\eps}_-$, respectively,  are chosen uniformly among the living particles at time $t$. Then, as $N\to\infty$, the probability joint density function for their positions converges to
$$c_{(n,m)}(t)\,\prod_{i=1}^{n}u_{+}(t,r_i) \prod_{j=1}^{m}u_{-}(t,s_j)$$
uniformly for $(\vec{r},\vec{s})\in \bar{D}_+^n\times \bar{D}_-^m$ and for $t$ in any compact time interval, where $c_{(n,m)}(t)$ is a normalizing constant.

\subsection{Key ideas}

A key step in our proof of propagation of chaos (Theorem \ref{T:correlation}) is Theorem \ref{T:Uniqueness_hierarchy}. The latter   establishes uniqueness of solution for the infinite system of equations satisfied by the correlation functions of the particles in the limit $N\to\infty$. Such infinite system of equations is sometimes called \emph{BBGKY-hierarchy} in statistical physics. Our  BBGKY hierarchy involves boundary terms on the interface, which  is new to the literature.
Our proof of uniqueness involves a representation and manipulations of the hierarchy in terms of trees. This technique is related to but different from that in \cite{lEbShtY07} which used Feynman diagrams. It is potentially useful in the study of other stochastic models involving  coupled differential equations.

To establish hydrodynamic limit result (Theorem \ref{T:conjecture}), we employ the classical tightness plus finite dimensional distribution approach. Tightness of $(\X^{N,+},\,\X^{N,-})$ in the Skorokhod space is proved in Theorem \ref{T:tight}. This together  with the propagation of chaos result (Theorem \ref{T:correlation}) establishes the  hydrodynamic limit of the interacting
random walks.

Two new tools for discrete approximation of random walks in domains are developed in this article. Namely, the local central limit theorem (local CLT) for reflected random walk on bounded Lipschitz domains (Theorem \ref{T:LCLT_CTRW}) and the `discrete surface measure' (Lemma \ref{L:DiscreteApprox_SurfaceMea}). We believe these tools are potentially useful in many discrete schemes which involve reflected Brownian motions.

Weak convergence of simple random walk on $D^{\eps}_{\pm}$ to RBM has been established for general bounded domains in \cite{BC08} and \cite{BC11}.
However, we need more for our model; namely a local convergence result which guarantees that the convergence rate is uniform up to the boundary. For this, we establish the local CLT. We further generalize the weak convergence result and the local limit theorem to deal with RBMs with gradient drift. There are two reasons for us to consider gradient drift. First, it is physically natural to assume the particles are subject to an electric potential. Second, the maximal extension theorem, \cite[Theorem 6.6.9]{CF12}, which is a crucial technical tool used in \cite{BC08} and \cite{BC11}, has established only in  symmetric setting. The proof of the local CLT is based on a `discrete relative isoperimetric inequality' (Theorem \ref{T:Isoperimetric_Discrete}) which leads to the Poincar\'e inequality and the Nash inequality. The crucial point is that these two inequalities are uniform in $\eps$ (scaling of lattice size) and is invariant under the dilation of the domain $D \mapsto aD$.

The paper is organized as follows. In section 2, we introduce the stochastic model and some preliminary facts that will be used later. We then prove the existence and uniqueness of solution for the coupled PDE. The main results, Theorem \ref{T:conjecture} and Theorem \ref{T:correlation}, will be rigorously formulated. We also mention various extensions of our main results in Remark \ref{Rk:GeneralizationResults}. Section 3 and section 4 contains the proof of Theorem \ref{T:correlation} and Theorem \ref{T:conjecture} respectively. Section 5 is devoted to the proofs of the discrete relative isoperimetric inequality and the local CLT.

\section{Notations and Preliminaries}

For the reader's convenience, we list our notations here:

\begin{longtable}{ll}
$\mathbb{Z}$ &  set of all integers\\
$\mathbb{Z}_+$ & $\{1,2,3,\cdots\}$ positive integers\\
$\mathbb{N}$ & $ \{0,1,2,\cdots\}$ non-negative integers\\
$\mathbb{R}$ & set of all real numbers\\
&\\
$\mathcal{B}(E)$ & Borel measurable functions on $E$ \\
$\mathcal{B}_b(E)$& bounded Borel measurable functions on $E$ \\
$\mathcal{B}^+(E)$& non-negative Borel measurable functions on $E$\\
$C(E)$ & continuous functions on $E$\\
$C_b(E)$ & bounded continuous functions on $E$\\
$C^+(E)$ & non-negative continuous functions on $E$\\
$C_c(E)$ & continuous functions on $E$ with compact support\\
$D([0,\infty),\,E)$ & space of $c\grave{a}dl\grave{a}g$ paths from $[0,\infty)$ to $E$ \\
&\quad equipped with the Skorokhod metric\\
&\\
$\mathcal{H}^m$ &  $m$-dimensional Hausdorff measure\\
$M_+(E)$ & space of finite non-negative Borel measures on $E$ \\
& \quad equipped with the weak topology\\
$M_{\leq 1}(E)$ & $\{\mu\in M_{+}(E):\,\mu(E)\leq 1\}$  \\
$M_{1}(E)$ (or $\mathcal{P}(E)$) & $\{\mu\in M_{+}(E):\,\mu(E)=1\}$ \\
$\eta^{\eps,\pm}_t(x)$ & number of living particles at $x\in D^{\eps}_{\pm}$ at time $t$\\
$(\eta^{\eps}_t)_{t\geq 0}$&  process with generator $\mathfrak{L}^{\eps}=\mathfrak{L}^{\eps}_0+\mathfrak{K}^{\eps}$ in Definition \ref{Def:ConfigurationProcess_eta}\\
$(\xi^{0}_t)_{t\geq 0}$ and $(\eta^{0}_t)_{t\geq 0}$& independent processes with generator $\mathfrak{L}^{\eps}_0$\\
$E^{\eps}$ & $\mathbb{N}^{D^ {\eps}_{+}} \times \mathbb{N}^{D^ {\eps}_{-}}$, state space of $(\eta^{\eps}_t)_{t\geq 0}$ \\
$\X^{N,\pm}_t(dz)$ & $\frac{1}{N}\sum_{x \in D^{\eps}_{\pm}}\eta^{\pm}_t(x)\1_x(dz)$,
the normalized empirical measure in $\bar{D}_{\pm}$ \\
$\mathfrak{E}$& $M_{\leq 1}(\bar{D}_+) \times M_{\leq 1}(\bar{D}_-)$, the state space of $(\X^{N,+}_t,\,\X^{N,-}_t)_{t\geq 0}$\\
$\{\F^X_t:\,t\geq 0\}$ & filtration induced by the process $(X_t)$, i.e. $\F^X_t=\sigma(X_s,\,s\leq t)$\\
$\1_x$ &  indicator function at $x$ or the Dirac measure at $x$\\
&\quad (depending on the context)\\
$\toL$ & weak convergence of random variables (or processes)\\
$\<f,\,\mu\>$ & $\int f(x)\,\mu(dx)$ \\
$x\vee y$ & $\max\{x,\,y\}$\\
$x\wedge y$ & $\min\{x,\,y\}$\\
$C,\,C_1,\,C_2,\,\cdots$ & positive constants\\
$I^{\eps}$ & `$\eps$-point approximation' of $I$ constructed in Lemma \ref{L:DiscreteApprox_SurfaceMea}\\
$\sigma_{\eps}$ & `discrete surface measure' constructed in Lemma \ref{L:DiscreteApprox_SurfaceMea}\\
\end{longtable}

\begin{tabular}{llllll}
Process & Semigroup & Heat kernel & Measure & Generator & State space\\
$X^{\pm}(t)$ & $P^{\pm}_t$ & $p^{\pm}(t,x,y)$ &  $\rho_{\pm}$  &  $\A^{\pm}$  & $\bar{D}_{\pm}$\\
$X^{\eps,\pm}(t)$ & $P^{\eps,\pm}_t$ & $p^{\eps,\pm}(t,x,y)$ &  $m^{\pm}_{\eps}$  &  $\A^{\pm}_{\eps}$ & $D^{\eps}_{\pm}$\\
$X_{(n,m)}(t)$ & $P^{(n,m)}_t$ & $p=p^{(n,m)}$ &  $\rho=\rho_{(n,m)}$  &  $\A^{(n,m)}$  & $\bar{D}_{+}^n\times \bar{D}_{-}^m$\\
$X_{(n,m)}^{\eps}(t)$ & $P^{(n,m),\eps}_t$ & $p^{\eps}=p^{(n,m),\eps}$ &  $m_{\eps}$  &  $\A^{(n,m)}_{\eps}$ & $(D^{\eps}_{+})^n\times (D^{\eps}_{-})^m$
\end{tabular}

where in the above,
\begin{align*}
p^{(n,m)}(t,(\vec{r},\vec{s}),(\vec{r'},\vec{s'}))&:= \prod_{i=1}^n p^{+}(t,r_i,r_i')\prod_{j=1}^m p^{-}(t,s_j,s_j')\\
\rho_{(n,m)}(\vec{r},\vec{s})&:= \prod_{i=1}^n\rho_+(r_i)\prod_{j=1}^m\rho_-(s_j).
\end{align*}

We also use the following abbreviations:

\begin{tabular}{ll}
a.s. & almost surely\\
LDCT & Lebesque dominated convergence theorem\\
CTRW & continuous time random walk\\
RBM & reflected Brownian motion\\
local CLT & local central limit theorem\\
LHS & left hand side\\
RHS & right hand side\\
WLOG & without loss of generality
\end{tabular}

\bigskip

\begin{definition}
A Borel subset $E$ of $\R^d$ is called $\mathcal{H}^m$\textbf{-rectifiable} if $E$ is a countable union of Lipschitz images of bounded subsets of $\R^m$ with $\mathcal{H}^m(E)<\infty$ (As usual, we ignore sets of $\mathcal{H}^m$ measure 0). Here $\mathcal{H}^m$ denotes the $m$-dimensional
Hausdorff measure.
\end{definition}

\begin{definition}
A \textbf{bounded Lipschitz domain} $D\subset\R^d$ is a bounded connected open set such that for any $\xi\in \partial{D}$, there exits $r_{\xi}>0$ such that $B(\xi,r_{\xi})\cap D$ is represented by $B(\xi,r_{\xi}) \cap  \{(y',y^{d})\in \R^d: \phi_{\xi}(y')<y^d\}$ for some coordinate system centered at $\xi$ and a Lipschitz function $\phi_{\xi}$ with Lipschitz constant $M$, where $M=M_D>0$ does not depend on $\xi$ and is called the Lipschitz constant of $D$.
\end{definition}

\begin{assumption}\label{A:Setting}
$D_{\pm}$ are given adjacent bounded Lipschitz domains in $\R^d$ such that $I:= \bar{D}_{+}\cap \bar{D}_{-}=\partial D_+\cap \partial D_-$ is a finite union of disjoint connected $\mathcal{H}^{d-1}$-rectifiable sets, $\rho_{\pm}\in W^{(1,2)}(D_{\pm})\cap C^1(\bar{D}_{\pm})$ are given functions which are strictly positive, $\lambda>0$ is a fixed parameter.
\end{assumption}

\subsection{Interacting random walks in domains}

In this subsection, we describe the interacting random walk model. We start with
  some key ingredients needed in  discrete approximation.

\subsubsection{Discrete approximation of surface measure}

To capture the boundary behavior of our processes near the interface $I$ in the discrete scheme, we need a discrete approximation of the surface measure $\sigma$
on $I$. The construction of $I^{\eps}$ and $\sigma_{\eps}$ in the following lemma is a key to our approximation scheme. For us,  $\mathbb{N}:= \{0,1,2,\cdots\}$ denotes the set of non-negative integers.

\begin{lem}\label{L:DiscreteApprox_SurfaceMea}
Suppose $D$ is a bounded Lipschitz domain of $\R^d$. Let $I\subset \partial D$ be closed, connected and $\mathcal{H}^{d-1}$-rectifiable. Let $\eps_j = 2^{-j}$ for $j\in \mathbb{N}$. Then there exist finite subsets $I^{(j)}=I^{\eps_j}$ of $I$ and functions $\sigma_{(j)}=\sigma_{\eps_j}: \,I^{(j)}\rightarrow [\eps^{d-1}/C,\, C\eps^{d-1}]$ such that (a) and (b) below hold simultaneously:
\begin{enumerate}
\item[(a)]  \begin{equation}\label{E:NumberOfPoints_SurfaceMea}
                \sup_{x\in \bar{D}}\, \#\,\left( I^{(j)}\cap B(x,\,s)\right) \leq C\,\left(\frac{s}{\eps_j}\vee 1 \right)^{d-1}\, \quad \forall\,s\in(0,\infty),\,j \in\mathbb{N},
            \end{equation}
            where $\# A $ denotes the number of elements in the finite set $A$, $B(x,\,s)=\{y\in \R^d:\,|y-x|<s\}$ is the ball with radius $s$ centered at $x$,
 and $C$ is a constant that depends only on $D$.
\item[(b)] For any equi-continuous and uniformly bounded family $\F\subset C(I)$,
            \begin{equation}\label{E:WeakConverge_SurfaceMea}
                \lim_{j\to \infty}\, \sup_{f\in\F} \Big|\,\sum_{I^{(j)}}f\,\sigma_{(j)}- \int_{I}f\,d\sigma\,\Big|=0.
            \end{equation}
\end{enumerate}
\end{lem}

\begin{pf}
We can always split $I$ into small pieces. The point is to guarantee that each piece is not too small, so that $\sigma_{(j)}/\eps^{d-1}\geq C$ and that (\ref{E:NumberOfPoints_SurfaceMea}) holds. Since $I$ is $\mathcal{H}^{d-1}$-rectifiable, we have
$$
C^{-1}\,R^{d-1} \leq \sup_{x\in I} \mathcal{H}^{d-1}(I \cap B(x,R)) \leq C\,R^{d-1}
$$
for  $R\in (0, 1]$, where $C$ does not depend on $R$.
Since $I$ is closed, it is regular with dimension $d-1$ in the terminology of section 1 of \cite{DS91}. Hence by \cite{gD88} or section 2 of \cite{DS91}, we can build ``dyadic cubes'' for $I$. More precisely, there exists a family of partitions $\{\Delta_j\}_{j\in \Z}$ of $I$ into ``cubes'' $Q$ such that
\begin{enumerate}
\item[(i)]   if $j\leq k$, $Q\in \Delta_j$ and $Q'\in\Delta_k\,$, then either $Q\cap Q' =\emptyset$ or $Q\subset Q'$;
\item[(ii)]   if $Q\in \Delta_j\,$, then $$C^{-1}\,2^j \leq {\rm diam} ( Q)  \leq C\,2^j\quad \text{ and}$$  $$C^{-1}\,2^{j(d-1)}\leq \mathcal{H}^{d-1}(Q) \leq C\,2^{j(d-1)};$$
\item[(iii)]   $$\mathcal{H}^{d-1}(\{x\in Q:\, {\rm dist} (x,\,I\setminus Q)\leq r2^j\}) \leq C\,r^{1/C}\,2^{j(d-1)}$$
        for all $Q\in \Delta_j$ and $r>0$.
\end{enumerate}
Here the constant $C$ is independent of $j,\,Q,$ or $r$.
Note that $\mathcal{H}^{d-1}$ is the surface measure $\sigma$ of $\partial D$ and that property (iii) tells us that the cubes have relatively small boundary. In particular, (iii) implies $\sigma (\partial Q \cap I)=0$ for all cube $Q$.

Suppose $\Delta_j= \left\{U_i^{(j)}\right\}_{i=1}^{k_{j}}$.
We pick one point $z^{(j)}_i$ from each $U_i^{(j)}$ to form the set $I^{(j)}$. Finally, we define $\sigma_{(j)}(z^{(j)}_i):= \sigma(U_i^{(j)})$. It follows from (ii) that $\sigma_{(j)}\in[\eps^{d-1}/C,\, C\eps^{d-1}]$ for some $C$ which depends only on $D$. The inequality (\ref{E:NumberOfPoints_SurfaceMea}) follows from $C^{-1}\,\eps_j^{d-1}\leq \sigma\,(U_i^{(j)})$ and the Lipschitz property of $\partial D$. It remains to check (\ref{E:WeakConverge_SurfaceMea}).

Fix any $\eta>0$. There exists $\lambda=\lambda(\eta)>0$ such that $|f(x)-f(y)|<\eta$ whenever $|x-y|<\lambda$. Hence for $j$ large enough (depending only on $\lambda$),
\begin{equation*}
\Big|\int_{I}g\,d\sigma - \sum_{I^{(j)}}g\,\sigma_{(j)}\Big|
 = \left|\sum_{i}\,\Big(\int_{U^{(j)}_i}g\,d\sigma - g(z^{(j)}_i)\sigma(U^{(j)}_i)\Big)\right| \leq \eta\,\sum_{i}\sigma(U^{(j)}_i)=\eta\,\sigma(I).
\end{equation*}
The desired convergence (\ref{E:WeakConverge_SurfaceMea}) now follows.
\end{pf}

\begin{remark} \rm
(\ref{E:WeakConverge_SurfaceMea}) implies that we have the weak convergence $\sum_{z\in I^{(j)}}\,\sigma_{(j)}\,\delta_{z} \to \sigma\big|_{I}$
on the space $M_{+}(I)$ of positive finite measure Borel measures on $I$. Here $\delta_{z}$ is the dirac delta measure at $z$, and $\sigma\big|_{I}$ is the surface measure restricted to $I$. (\ref{E:NumberOfPoints_SurfaceMea}) is a control on the number of points locally in $I^{\eps_j}$. We call $I^{\eps}$ the `$\eps$-point approximation' of $I$ and $\sigma_{\eps}$ the `discrete surface measure' associated to $I^{\eps}$. \qed
\end{remark}

\begin{remark} \rm
The above lemma remains  true if $I$ is the finite union of disjoint closed connected and $\mathcal{H}^{d-1}$-rectifiable subsets of $\partial D$. This enables us to deal with disconnected interface $I$. \qed
\end{remark}

\subsubsection{Reflected diffusion and random walk approximation }\label{UnderlyingMotion}

We now describe the motion of each underlying particle. First we fix a  bounded Lipschitz domain $D\subset \R^d$ and any $\eps >0$. Without loss of generality, we assume that the origin $0\in D$. Let  $\bar{\eps \Z^d}$ be the union of all closed line segments joining nearest neighbors in $\eps \Z^d$, and $(D^{\eps})^{*}$ the connected component of $D\cap \bar{\eps \Z^d}$ that contains the point $0$.  Set  $D^{\eps} = (D^{\eps})^{*} \cap \eps\Z^d$.
We can view $D^\eps$ as the vertices of a graph whose edges coming from $(D^{\eps})^{*}$.
 We also denote the graph-boundary $\partial D^{\eps}:= \{x\in D^{\eps}:\,v_{\eps}(x)<2d\}$,  where $v_{\eps}(x)$ is the degree of $x$ in $D^{\eps}$.

Suppose $\rho\in W^{1,2}(D)\cap C^1(\bar{D})$ is strictly positive. Define
$$\D (f,g) := \dfrac{1}{2} \int_{D}\nabla f(x)\cdot  \nabla g(x)\,\rho(x) \,dx.$$
Since $D$ is  Lipschitz, $(W^{1,2}(D), \D)$ is a regular Dirichlet form on $L^2(D;\rho)$ and so there is a $\rho$-symmetric diffusion $X$ associated with it (cf. for example \cite{zqChen93}).

\begin{definition}
We call $X$ the $(I_{d\times d},\rho)$\textbf{-reflected diffusion}, where $I_{d\times d}$ is the $d\times d$ identity matrix. When $\rho=1$, $X$ is called the \textbf{reflected Brownian motion} (RBM) in $D$. Hence a $(I_{d\times d},\rho)$-reflected diffusion is a RBM in$D$ with drift $\frac{1}{2}\nabla\log\rho$.
\end{definition}
The $L^2$-infinitesimal generator of $X$ is
$$\A=\frac{1}{2\rho}\,\nabla\cdot\,(\rho\,\nabla) = \frac{1}{2}\Delta + \frac{1}{2}\,\nabla (\log \rho)\cdot\nabla.
$$
Moreover, $X$ has the Skorokhod representation:
\begin{eqnarray}
dX_t = dB_t + \frac{1}{2}\nabla\log\rho(X_t)dt + \vec{n}(X_t)dL_t
 \quad \hbox{for } t\geq 0,\,\P_x \hbox{-a.s. } x\in \bar{D},
\end{eqnarray}
where $\vec{n}$ is the inward unit normal of $\partial D$ and $L$ is the positive continuous additive function (PCAF) of $X$ whose Revuz measure is $\frac{1}{2}\sigma$ (c.f.\cite{zqChen93}).
We call $L$ the boundary local time of $X$.

Next, we define $X^{\eps}$ to be a continuous time random walk (CTRW) on $D^{\eps}$ with exponential waiting time of parameter $\frac{d}{\eps^{2}}$ and one step transition probabilities
$$
p_{xy}:= \frac{\mu_{xy}}{\sum_{y}\mu_{xy}},
$$
 where $\{\mu_{xy}:\,x,y\in D^{\eps}\}$ are symmetric weights (conductances) to be constructed in two steps as follows:
First, for every $x\in D^{\eps}\setminus \partial D^{\eps}$ and $i=1,2,\cdots,d$, define
\begin{eqnarray}
\mu_{x,x+\eps \vec{e_i}} &:=& \left(1+\frac{1}{2}\ln\frac{\rho(x+\eps \vec{e_i})}{\rho(x)}\right)\,\left(\frac{\rho(x)+\rho(x+\eps \vec{e_i})}{2}\right)\,\frac{\eps^{d-2}}{2} \label{E:Conductance_BaisedRW+}\\
\mu_{x,x-\eps \vec{e_i}} &:=& \left(1+\frac{1}{2}\ln\frac{\rho(x)}{\rho(x-\eps \vec{e_i})}\right)\,\left(\frac{\rho(x)+\rho(x-\eps \vec{e_i})}{2}\right)\,\frac{\eps^{d-2}}{2} \label{E:Conductance_BaisedRW-}
\end{eqnarray}
Clearly, $\mu_{xy}=\mu_{yx}$ for all $x,y\in D^{\eps}\setminus \partial D^{\eps}$.
Note that since $\rho$ is $C^1$ and strictly positive on $\bar D$, when $\eps$ is sufficiently small,
$\mu_{x,x+\eps \vec{e_i}}$ and $\mu_{x,x-\eps \vec{e_i}}$ are strictly positive for every
 $x\in D^{\eps}\setminus \partial D^{\eps}$ and $i=1,2,\cdots,d$.
Second, we define
\begin{equation*}
\mu_{xy} :=
\begin{cases}
\mu_{yx}, &\text{ if  } x\in \partial D^{\eps},\,y\in D^{\eps}\setminus \partial D^{\eps} \\
\frac{\eps^{d-2}}{2}, &\text{ if  } x,y\in \partial D^{\eps} \text{  are adjacent in } D^{\eps}
\end{cases}
\end{equation*}
Now $\mu_{xy}=\mu_{yx}$ for all $x,y\in D^{\eps}$. A heuristic reason of the above construction can be found in \cite{wtF14}.

We call $X^{\eps}$ the $\eps$\textbf{-approximation} of $X$. Clearly, $X^{\eps}$ is symmetric with respect to the measure $m_{\eps}$ defined by
\begin{equation*}
m_{\eps}(x):= \frac{\eps^2}{d}\sum_{y}\mu_{xy}.
\end{equation*}

Since $\rho\in C^{1}(\bar{D})$, there exists a constant $C>0$ such that
\begin{eqnarray}\label{E:Approximate_m(x)} C^{-1}\leq \inf_{x}\frac{m_{\eps}(x)}{\eps^d}\leq \sup_{x}\frac{m_{\eps}(x)}{\eps^d}\leq C.
\end{eqnarray}
Moreover, $\lim_{\eps \to 0}\frac{m_{\eps}(x^{\eps})}{\eps^d}= \rho(x)$ whenever $x^{\eps}\in D^{\eps}$ converges to $x\in D$.

A special but important case is when $\rho\equiv 1$. In this case,  $X$ is simply the reflected Brownian motion on $D$, and
$X^{\eps}$ is a simple random walk on the graph $D^\eps$. It is proved in \cite{BC08} that $X^\eps$ converges weakly
to the reflected Brownian motion $X$ as $\eps \to 0$.

Recall that by Assumption \ref{A:Setting}, we are given $\rho_{\pm}\in W^{(1,2)}(D_{\pm})\cap C^1(\bar{D}_{\pm})$. We denote by $X^{\pm}$  a $(I_{d\times d},\rho_{\pm})$-reflected diffusion in $D_{\pm}$, and by $X^{\eps,\pm}$ the $\eps$-approximation of $X^{\pm}$.

\subsubsection{Random walks with interaction}

\begin{figure}[h]
	\begin{center}
    \vspace{-0.5em}
	\includegraphics[scale=.22]{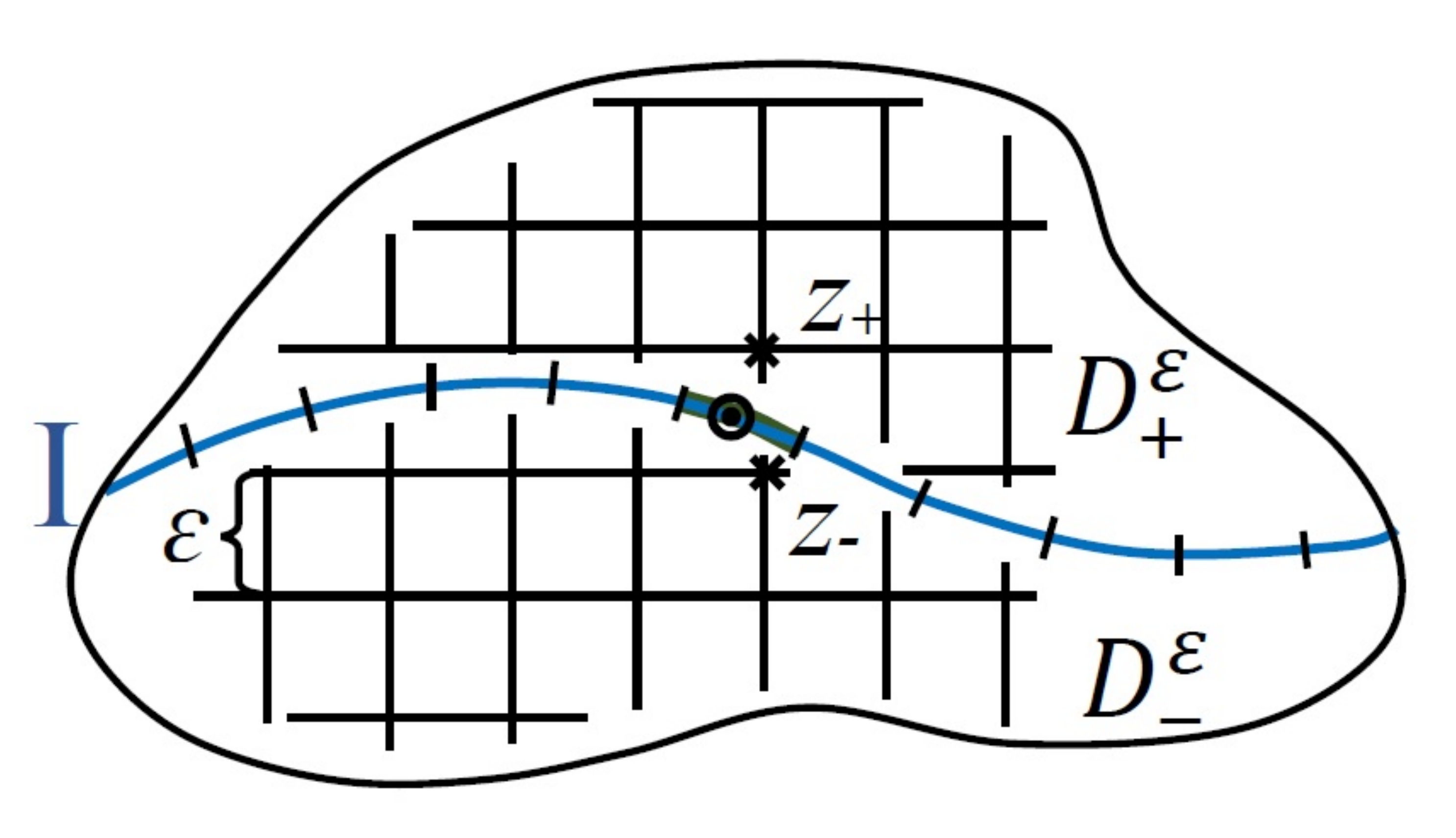}
    \vspace{-1em}
    \caption{$z\in I^{\eps}\subset I,\;z_{\pm}\in D^{\eps}_{\pm}$}\label{fig:KillingPairsDiscrete}
    \vspace{-1em}
	\end{center}
	\end{figure}

Fix $\eps=\eps_j=2^{-j}$ ($j\in \mathbb{N}$) and $N=2^{jd}$ such that $N \eps^d=1$. Assume there are $N$ ``+'' particles in $D^ {\eps}_{+}$ and $N$ ``$-$'' particles in $D^ {\eps}_{-}$ at $t=0$. Each particle moves as an independent CTRW $X^{\eps,\pm}$ (see the previous subsection) in its respective domain $D^{\eps}_{\pm}$. Let $I^{\eps}$ be the finite subset of $I$ defined in Lemma \ref{L:DiscreteApprox_SurfaceMea}. For each $z\in I^{\eps}$, pick an $z_{+}\in D^{\eps}_{+}$ and an $z_{-} \in D^{\eps}_{-}$ which are closest to $z$ (See Figure \ref{fig:KillingPairsDiscrete}). A pair of particles of opposite charges at $(z_+,z_-)$ is being killed with a certain rate to be explained. Note that for $\eps$ small enough, we have $\sup_{z\in I^{\eps}}|z_{\pm}-z|\leq 2M\eps$, where $M$ is the Lipschitz constant of $I$.

The state space of the particle system is the collection of configurations
\begin{equation}\label{e:2.7}
E^{\eps} :=\left\{ \eta^\eps=(\eta^{\eps, +}, \eta^{\eps, -}): \eta^{\eps, \pm } :  D^\eps_\pm\to   \mathbb{N} \right\}  .
\end{equation}
The state of the particle system at time $t$ is a random element $\eta^{\eps}_{t}=(\eta^{\eps,+}_{t},\eta^{\eps,-}_{t}) \in E^{\eps}$. Here $\eta^{\eps,\pm}_t(x)$ stands for the number of $``\pm''$ particles at $x \in D^ {\eps}_{\pm}$ at time $t$. We omit $\eps$ and $N$ for convenience when there is no ambiguity. For example, we write $\eta_t$ and $m(x)$ in place of $\eta^{\eps}_t$ and $m_{\eps}(x)$ respectively. The function
 $\xi$ such that $\xi (x)=1$ and $\xi (y)=0$ for $y\neq x$ is denoted as $\1_x$.

\begin{definition}\label{Def:ConfigurationProcess_eta}
$\eta_{t}$ is defined to be the unique strong Markov process which has the generator $\mathfrak{L}=\mathfrak{L}^{\eps}$ given by
\begin{equation}\label{E:generator0}
\mathfrak{L}^{\eps}:=\mathfrak{L}^{\eps}_0+\mathfrak{K}^{\eps}\,,
\end{equation}
where $\mathfrak{L}^{\eps}_0$ is the generator of two families of independent random walks in $D^\eps_+$ and $D^\eps_-$, respectively,
 with no annihilation between them, namely
\begin{eqnarray}
\mathfrak{L}^{\eps}_0 f(\eta)
& := & \frac{d}{\eps^2}\,\sum_{x,y \in D^{\eps}_{+}}\eta^+(x)p^+_{xy}\{f(\eta^{+}-\1_x+\1_y,\eta^-)-f(\eta)\} \nonumber  \\
&+&  \frac{d}{\eps^2}\,\sum_{x,y \in D^{\eps}_{-}} \eta^-(x)p^-_{xy}\{f(\eta^{+},\eta^{-}-\1_x+\1_y)-f(\eta)\} \label{e:2.9}
\end{eqnarray}
and $\mathfrak{K}^{\eps}$ is the operator corresponding to annihilation between particles of opposite signs at the interface $I^\eps$, namely
\begin{equation}\label{e:2.10}
\mathfrak{K}^{\eps} f(\eta)
 :=  \dfrac{\lambda}{\eps} \, \sum_{z\in I^{\eps}} \,\Psi_{\eps}(z)\,\eta^+(z_+)\eta^-(z_-)\,\{f(\eta^{+}-\1_{z_+},\eta^{-}-\1_{z_-})-f(\eta)\},
\end{equation}
where $p^{\pm}_{xy}$ is the one-step transition probabilities for the CTRW $X^{\eps,\pm}$
 on $D^\eps_\pm$ (without any interaction) and
\begin{equation}\label{e:2.8}
\Psi_{\eps}(z):= \frac{\sigma_{\eps}(z)}{\eps^{d-1}}\,\dfrac{\eps^{2d}}{m(z_+)m(z_-)}
\end{equation}
with $\sigma_{\eps}$ and $I^{\eps}$ being constructed by Lemma \ref{L:DiscreteApprox_SurfaceMea}.
\end{definition}

The expression for $\mathfrak{K}^{\eps}$ comes from the underlying assumptions of the model: First, the term $\eta^+(z_+)\eta^-(z_-)$ is combinatorial in nature. Since there are $\eta^+(z_+)\eta^-(z_-)$ pairs of particles at position $(z_+,z_-)$, the chance of killing is proportional to the number of ways of selecting a pair of particles near the interface. Second, each pair of particles near $I$ disappears at rate $(\lambda/\eps)\,\Psi_{\eps}(z)$ where
$\lambda$ is a parameter. Intuitively, in the limit, the amount of annihilation in a neighborhood of a point is proportional to the surface area of the interface $I$ in that neighborhood. The scaling $1/\eps$ is suggested by the observation that there are about $1/\eps$ "layers" starting from the interface $I$, so that the chance for a particle to arrive near $I$ is of order $\eps$. $\Psi_{\eps}(z)$ is comparable to 1 and  can be viewed as a normalizing constant with respect to the lattice. This choice (\ref{e:2.8}) is justified in the proof of Theorem \ref{T:hierarchy}.

\subsection{Discrete heat kernel}

Throughout this subsection, $D$ is a bounded Lipschitz domain in $\R^d$, $\rho\in W^{(1,2)}(D)\cap C^1(\bar{D})$ is strictly positive, $X$ is a $(I_{d\times d},\rho)$-reflected diffusion. It is well known (cf. \cite{BH91, GSC11} and the references therein) that $X$ has a transition density $p(t,x,y)$ with respect to  the symmetrizing measure $\rho(x)dx$ (i.e., $\P_x(X_t\in dy)=p(t,x,y)\,\rho(y)dy$ and $p(t,x,y)=p(t,y,x)$), that $p$ is locally H\"older continuous and hence $p\in C((0,\infty)\times \bar{D}\times \bar{D})$, and that we have the followings: there are constants $c_1\geq 1$ and $c_2\geq 1$ such that
   \begin{equation}\label{E:Gaussian2SidedHKE}
        \dfrac{1}{c_1 t^{d/2}}\,\exp\left(\frac{-c_2 |y-x|^2}{t}\right)
\leq p(t,x,y) \leq \dfrac{c_1}{t^{d/2}}\,\exp\left(\frac{-|y-x|^2}{c_2\,t}\right)
        \end{equation}
for every $(t,x,y)\in(0,\infty)\times \bar{D}\times \bar{D}$. Using (\ref{E:Gaussian2SidedHKE}) and the Lipschitz assumption for the boundary, we can check that
   \begin{equation}\label{E:boundary_strip_boundedness}
            \sup_{x\in\bar{D}}\,\sup_{0<\delta\leq \delta_0}\,\frac{1}{\delta}\int_{D^{\delta}}p(t,x,y)\,dy \leq \frac{C_1}{t^{1/2}}+C_2
\quad \hbox{for } t>0\quad \text{and hence}
        \end{equation}
   \begin{equation}\label{E:Surface_integral_boundedness}
            \sup_{x\in\bar{D}}\,\int_{\partial D}p(t,x,y)\,\sigma(dy) \leq \frac{C_1}{t^{1/2}}+C_2 \quad \hbox{for } t>0,
        \end{equation}
where $C_1,\,C_2,\,\delta_0>0$ are constants which depends only on $d$, $D$ and $\rho$.

On other hand, suppose $g\in \mathcal{B}_b( [0,T]\times\partial D)$. Then for $t\in[0,T]$ and  $x\in \bar{D}$,
        \begin{equation}\label{E:localtime_def}
            \E^x\left[\int_0^tg(s,X_s)dL_s\right] = \dfrac{1}{2}\int_0^t\int_{\partial D} g(s,y)p(s,x,y) \sigma(dy)ds.
        \end{equation}

Now let $X^{\eps}$ be the $\eps$-approximation of $X$ with symmetrizing measure $m_{\eps}$.
The transition density $p^{\eps}$ of $X^{\eps}$ with respect to  the measure $m_\eps$ is
\begin{equation}
p^{\eps}(t,x,y) := \dfrac{\P^{x} (X^{\eps}_t=y)}{ m_\eps(y)}, \quad  t>0,\, x,\,y\in D^{\eps}.
\end{equation}
Clearly, $p^{\eps}$ is strictly positive and is symmetric in $x$ and $y$.

We will prove in  Section \ref{S:5} that $p^{\eps}$ enjoys two-sided Gaussian bound and is jointly H\"older continuous
uniform in $\eps \in (0, \eps_0)$ for some $\eps_0>0$,
and that $p^{\eps}$ converges to $p$ uniformly on compact subsets of $(0,\infty)\times \bar{D}\times \bar{D}$. In rigorous terms, we have the following three results. The important point is that the constants involved in these results do not depend on $\epsilon$.

\begin{thm}\label{T:UpperHKE}(Gaussian upper bound)
There exist $C_k=C_k(d,D,\rho,T)>0$, $k=1, 2$, and $\eps_0=\eps_0(d,D,\rho,T) \in (0, 1]$ such that for every $\eps\in(0,\eps_0)$ and
$x,y \in D^{\eps}$,
\begin{equation}\label{e:2.14}
p^{\eps}(t,x,y) \leq \dfrac{C_1}{(\eps\vee t^{1/2})^d}\,\exp\left(- C_2 \frac{|x-y|^2}{t}\right)
\quad \hbox{for }  t\in [\eps,T],
\end{equation}
 and
\begin{equation}\label{e:2.15}
p^{\eps}(t,x,y) \leq \dfrac{C_1}{(\eps\vee t^{1/2})^d}\,\exp\left(-C_2 \frac{|x-y|}{ t^{1/2}}\right)
\quad \hbox{for } t\in (0, \eps).
\end{equation}
\end{thm}

Observe that  \eqref{e:2.14} implies that \eqref{e:2.15} also holds for $t\in [\eps, T]$.

\begin{thm}\label{T:LowerHKE}(Gaussian lower bound)
There exist $C_k=C_k(d,D,\rho,T)>0$, $k=1, 2$, and $\eps_0=\eps_0(d,D,\rho,T) \in (0, 1]$ such that for every $\eps\in(0,\eps_0)$ and
$x,y \in D^{\eps}$,
\begin{equation}\label{e:2.19}
p^{\eps}(t,x,y) \geq \dfrac{C_1}{(\eps\vee t^{1/2})^d}\,\exp\left(- C_2 \frac{|x-y|^2}{t}\right)
\quad \hbox{for }  t\in (0,T].
\end{equation}
\end{thm}

\begin{thm}\label{T:HolderCts}(H\"older continuity)
There exist positive constants $\gamma=\gamma (d,D,\rho)$, $\eps_0(d,D,\rho)$ and $C(d,D,\rho)$ such that for all $\eps\in(0,\eps_0)$, we have
\begin{equation}\label{E:HolderCts2}
|p^{\eps}(t,x,y)-p^{\eps}(t',x',y')| \leq C\, \dfrac{ (|t-t'|^{1/2}+ |x-x'|+ |y-y'|)^\gamma }
{ (t\wedge t')^{\sigma/2}\,(1 \wedge t\wedge t')^{d/2}}.
\end{equation}
\end{thm}

\begin{thm}\label{T:LCLT_CTRW}(Local CLT)
$$\lim_{n \to\infty} \sup_{t\in [a,b]} \sup_{x,y \in \bar{D}}\Big|p^{(2^{-n})}(t,x,y)\,-\,p(t,x,y)\Big| =0 $$
for any compact interval $[a,b]\subset (0,\infty)$.
\end{thm}

To establish the  tightness of $\{(\X^{N,+},\,\X^{N,-})\}$, we need the following uniform estimate for the heat kernel  $p^{\eps} (t, x, y)$ of CTRW on $D^\eps$ near the boundary of $D^\eps$.  It is the discrete analog of \eqref{E:Surface_integral_boundedness}.

\begin{lem}\label{L:q_near_I}
    There exist $C=C(d,D,\rho,T)>0$ and $\eps_0=\eps_0(d,D,\rho)>0$ such that
    \begin{equation}
    \sup_{x\in D^{\eps}}\,\eps^{d-1}\,\sum_{y\in \partial D^{\eps}}p^{\eps}(t,x,y) \leq \frac{C}{\eps \vee t^{1/2}}
    \end{equation}
    for all $t\in[0,T]$ and $\eps\in(0,\eps_0)$. Here $\partial D^{\eps}$ is the graph-boundary of $D^{\eps}$, which is all the vertices in $ D^{\eps}$ with degree less than $2d$.
\end{lem}

\begin{pf}
Fix $\theta\in[0,T]$. By the Gaussian upper bound in Theorem \ref{T:UpperHKE}, we have
\begin{eqnarray*}
&& \sum_{y\in \partial D^{\eps}}p^{\eps}(\theta,x,y) \\
&\leq& \frac{C_1}{(\eps \vee \theta^{1/2})^d}\,\sum_{y\in \partial D^{\eps}}\exp{\left(\frac{-|y-x|}{\eps \vee \theta^{1/2}}\right)} \\
&=& \frac{C_1}{(\eps \vee \theta^{1/2})^d}\,\int_0^{\infty}\,|\{y\in D^{\eps}:\,|f(y)|>r\}|\,dr \quad\text{by setting } f(y)=\1_{\partial D^{\eps}}(y)\,\exp{\left(\frac{-|y-x|}{\eps \vee \theta^{1/2}}\right)}\\
&=& \frac{C_1}{(\eps \vee \theta^{1/2})^d}\,\int_0^{1}\,|\{\partial D^{\eps}\cap B(x,\,(\eps \vee \theta^{1/2})(-\ln r))\}|\,dr \quad (\text{since }f\leq 1)\\
&=& \frac{C_1}{(\eps \vee \theta^{1/2})^{d+1}}\,\int_0^{\infty}\,|\{\partial D^{\eps}\cap B(x,\,s)\}|\,\exp{\left(\frac{-s}{\eps \vee \theta^{1/2}}\right)}\,ds
\quad (\text{where }s=(\eps\vee \theta^{1/2})(-\ln r) ) \\
&\leq& \frac{C_1}{(\eps \vee \theta^{1/2})^d}\,\vee\,\frac{C_2}{\eps^{d-1}(\eps \vee \theta^{1/2})^{d+1}}\,\int_0^{\infty} s^{d-1}\,\exp{\left(\frac{-s}{\eps \vee \theta^{1/2}}\right)}\,ds \\
&\leq& \frac{1}{\eps^{d-1}}\,\left(\frac{C_1}{\eps \vee \theta^{1/2}}\,\vee \,\frac{C_2}{\eps \vee \theta^{1/2}}\,\int_0^{\infty} w^{d-1}e^{-w}dw \right) \quad
( \text{ where } w= \frac{s}{\eps \vee \theta^{1/2}}).
\end{eqnarray*}
Here $C_i$ are all constants which depend only on $d$, $D$ and $T$. Note that in the second last line, we used the fact that
$ |\{\partial D^{\eps}\cap B(x,\,s)\}| \leq C((s/\eps)^{d-1}\vee 1) $, which follows from Lemma \ref{L:DiscreteApprox_SurfaceMea}. The proof is now complete.
\end{pf}

In general, we use $"\pm"$ for quantities related to $X^{\pm}$. For example,
$\A^{\pm}$, $(P^{\pm}_t)_{t\geq 0}$ and $p^{\pm}(t,x,y)$ denote the generator, semigroup and transition density of
the reflected diffusion $X^{\pm}$ in $D_\pm$ with respect to  $\rho_{\pm}(x) dx$.
In addition, we use $"\eps"$ for quantities related to the CTRWs in the discrete domains $D^{\eps}_{\pm}$. For example, $p^{\eps,\pm}(t,x,y)$ denotes the transition density of the CTRW $X^{\eps,\pm}$  on $D^\eps_\pm$ with respect to
the measure $m_{\eps}^{\pm}$. We also denote $p^{\eps}(t,(\vec{r},\vec{s}),(\vec{r'},\vec{s'})):= \prod_{i=1}^{n}p^{\eps}(t,r_i,r_i') \prod_{j=1}^{m}p^{\eps}(t,s_j,s_j')$ for $(\vec{r},\vec{s})\in (D^{\eps}_+)^n \times (D^{\eps}_-)^m$.

By applying Lemma \ref{L:q_near_I} to $p^{\eps, \pm}(t, x, y)$, then by the boundedness of $\Phi_\eps$ in \eqref{e:2.8}, Theorem \ref{T:LCLT_CTRW}, Lemma \ref{L:DiscreteApprox_SurfaceMea} and LDCT, we have the following approximation for the local time of $X^{\pm}$ on $I$.
\begin{prop}(Discrete local time)
\begin{equation}\label{E:q_near_I}
\lim_{\eps\to 0}\int_0^t \eps^{d-1}\,\sum_{z\in I^{\eps}}p^{\eps,\pm}(\theta,x,z_{\pm})\,\sigma_\eps (z)\,d\theta = \int_0^t\int_{I} p^{\pm}(\theta,x,z)\sigma(dz)\,d\theta.
\end{equation}
\end{prop}

\subsection{Hydrodynamic limit: system of nonlinear PDEs}
In this subsection, we provide suitable notion of solutions for the coupled PDE (\ref{E:coupledpde:+}) and (\ref{E:coupledpde:-}), and then prove the existence and uniqueness of the solution. Throughout this subsection, $D$ is a bounded Lipschitz domain, $\rho\in W^{(1,2)}(D)\cap C^1(\bar{D})$ is strictly positive, $X$ is a $(I_{d\times d},\rho)$-reflected diffusion, $\{P_t\}$ and $p(t,x,y)$ are the semigroup and the transition density of $X$, with respect to the measure $\rho(x)dx$.

Observe that (\ref{E:coupledpde:+}) is a second order parabolic equation for $u_+$ with Robin boundary condition, and (\ref{E:coupledpde:-}) is a similar equation for $u_-$. This leads us to consider the following Robin boundary problem, where $g\in \mathcal{B}_b([0,\infty)\times \partial D)$.
\begin{align} \label{E:mixedBVP_drift}
    \begin{cases}
        \dfrac{\partial u}{\partial t} = \frac{1}{2}\Delta u +\frac{1}{2}\nabla (\log \rho)\cdot\nabla u  & \hbox{for } x\in D ,\, t>0,  \\
        \dfrac{\partial u}{\partial \vec{n}} = \frac{1}{\rho}\,gu & \hbox{for }  x\in \partial D,\, t>0 ,  \\
        u(0,x)=\varphi(x) &\hbox{for }  x\in D.
    \end{cases}
\end{align}

By It\'o's formula and the Skorokhod representation for $X$, we see that a classical solution of (\ref{E:mixedBVP_drift}), should it exists, has the probabilistic representation
\begin{equation}\label{E:ProbabilisticRep_Robin}
u(t,x):= \E^{x} \left[\varphi(X_t) e^{-\int_0^tg(t-s,X_s)dL_s} \right].
\end{equation}
where $L$ is the boundary local time of $X$.

\begin{definition}\label{Def:Mild_Sol_ProbabilisticRep_Robin}
$u$ defined by (\ref{E:ProbabilisticRep_Robin}) is called a probabilistic solution of (\ref{E:mixedBVP_drift}).
\end{definition}

First we show that the function $u$  defined by \eqref{E:ProbabilisticRep_Robin}   is continuous.

\begin{lem}\label{L:properties_u}
Suppose $\varphi\in \mathcal{B}_b(\bar{D})$, $g\in \mathcal{B}^+_b([0,T]\times \partial D)$ and $u$ is defined by (\ref{E:ProbabilisticRep_Robin}). Then $u\in C((0,T]\times \bar{D})$. Moreover, if $\varphi\in C(\bar{D})$, then $u\in C([0,T]\times \bar{D})$.
\end{lem}
\begin{pf}
Observe that for any $r\in[0,t]$,
\begin{eqnarray}\label{E:Robin_continuity}
u(t,x)
&=& \E^{x}\left[\varphi(X_t) e^{-\int_r^t g(t-s,X_s)dL_s}\, e^{-\int_0^r g(t-s,X_s)dL_s}\right] \notag\\
&=& \E^{x}\left[\varphi(X_t) e^{-\int_r^tg(t-s,X_s)dL_s}\right] + \E^{x}\left[\varphi(X_t) e^{-\int_r^tg(t-s,X_s)dL_s} \left( e^{-\int_0^r g(t-s,X_s)dL_s}-1 \right) \right] . \notag\\
\end{eqnarray}
By Markov property, the first term is
$$\E^{x} \left[ \E^{X_r}[\varphi(X_{t-r}) e^{-\int_0^{t-r} g(t-r-s,X_s)dL_s} ] \right] = \E^{x}[u(t-r,X_r)].$$
Since $X$ has the strong Feller property (see \cite{BH91}) and $u$ is bounded,  $x\mapsto \E^{x}[u(t-r,X_r)]$ is continuous on $\bar{D}$ for any fixed $t>0$ and $r\in(0,t)$.

The second term of (\ref{E:Robin_continuity}) converges to 0 uniformly on $(0,T]\times \bar{D}$, as $r\to 0$. This is because its absolute value is bounded by
\begin{eqnarray*}
&& \|\varphi\|_\infty \E^{x}\left[ 1- e^{-\int_0^r g(t-s,X_s)dL_s} \right] \\
&\leq& \|\varphi\|_\infty
 \E^{x}\left[\int_0^r g(t-s,X_s) dL_s\right] \quad \text{by mean-value theorem} \\
&\leq& \|\varphi\|_\infty \,\|g\|_\infty \,\frac{1}{2}\,\int_0^r\int_{\partial D}p(s,x,y)\,\sigma(dy)\,ds \\
&\leq& \|\varphi\|_\infty \,\|g\|_\infty \,(2C_1\sqrt{r}+C_2r).
\end{eqnarray*}
Hence, $u$ is continuous in $x\in \bar{D}$.

By a similar calculation as in (\ref{E:Robin_continuity}), we have
\begin{eqnarray*}
u(t+a,x)-u(t,x)
&=&\E^{x}[u(t,X_a)-u(t,x)] \\
&& + \E^{x}\left[\varphi(X_{t+a}) e^{-\int_a^{t+a}g(t+a-s,X_s)dL_s}
\left(e^{-\int_0^a g(t+a-s,X_s)dL_s}-1\right) \right].
\end{eqnarray*}
Hence,
\begin{eqnarray*}
&& |u(t+a,x)-u(t,x)|  \\
&\leq&  \E^{x}[|u(t,X_a)-u(t,x)|] + \|\varphi\|_{\infty}\,\E^{x}\left[\int_0^a g(t+a-s,X_s)dL_s\right] \quad \text{by mean-value theorem} \\
&\leq& \int_{D}|u(t,z)-u(t,x)|p(a,x,z)\,dz + \|\varphi\|_\infty \,\|g\|_\infty \,\frac{1}{2}\,\int_0^a\,\int_{\partial D}p(s,x,z)\,\sigma(dz)ds.
\end{eqnarray*}
Both terms go to 0 uniformly in $x\in \bar{D}$ as $a$ goes to 0. (In fact, the first term goes to 0 uniformly since the semigroup $P_t$ is strongly continuous on $C(\bar{D})$. For the second term, $\int_0^a\,\int_{\partial D}p(s,x,z)\,\sigma(dz)ds \leq 2C_1\sqrt{a}+C_2a$ also goes to 0 uniformly in $x$.) Hence $u$ is continuous in $t\in(0,T]$ uniformly in $x\in \bar{D}$. Therefore, $u\in C((0,T]\times \bar{D})$. If $\varphi\in C(\bar{D})$, we can extend the above argument to show that $u\in C([0,T]\times \bar{D})$.
\end{pf}

\begin{remark}\rm  In fact, one can allow $g$ to be unbounded and show that  the conclusion of Lemma \ref{L:properties_u} remains true
if $g \sigma  $ satisfies a Kato class condition:
$$
\lim_{a\to 0} \sup_{x\in \bar{D}}\int_0^a\,\int_{\partial D}p(s,x,z)\, |g(t+a-s,z)|
\sigma(dz)ds  = 0 .
$$
\qed
\end{remark}

\begin{prop}\label{prop:MildSol_RobinPDE}
Suppose $\varphi\in C(\bar{D})$ and $g\in \mathcal{B}^+_b([0,T]\times \partial D)$. Then
$$
u(t,x) := \E^{x} \left[\varphi(X_t) e^{-\int_0^tg(t-s,X_s)dL_s} \right]
$$
 is the unique element in $C([0,T]\times \bar{D})$ that satisfies the following integral equation:
\begin{equation}
u(t,x)=P_t\varphi(x)-\dfrac{1}{2}\int_0^t \int_{\partial D}p(t-r,x,y)g(r,y)u(r,y)\,\sigma(dy)\,dr.
\end{equation}
\end{prop}

\begin{pf}
By Lemma \ref{L:properties_u}, $u(t,x) := \E^{x}[\varphi(X_t) e^{-\int_0^tg(t-s,X_s)\,dL_s}]$ lies in $C([0,T]\times \bar{D})$. Moreover, by Markov property and (\ref{E:localtime_def}) we have
\begin{eqnarray*}
u(t,x) &=& \E^{x}[\varphi(X_t)] - \E^{x}[\varphi(X_t)\,(1-e^{-\int_0^t g(t-s,X_s)\,dL_s})] \\
&=& P_t\varphi(x)- \E^x\left[\varphi(X_t)\, e^{-\int_r^t g(t-s,X_s)\,dL_s}\Big|^{r=t}_{r=0} \right] \\
&=& P_t\varphi(x)- \E^x\left[\varphi(X_t)\, \int_0^t g(t-r,X_r)\,e^{-\int_r^t g(t-s,X_s)\,dL_s}\,dL_r \right] \\
&=& P_t\varphi(x)- \E^x\left[\int_0^t g(t-r,X_r)\,\E^{X_r}\left[\varphi(X_{t-r})\,e^{-\int_0^{t-r} g(t-r-s,X_s)dL_s}\right]\,dL_r \right] \\
&=& P_t\varphi(x)- \E^x\left[\int_0^t g(t-r,X_r)\,u(t-r,X_r)\,dL_r \right] \\
&=& P_t\varphi(x)- \dfrac{1}{2}\int_0^t\int_{\partial D}p(r,x,y)g(t-r,y)u(t-r,y)\,\sigma(dy)dr.
\end{eqnarray*}
Hence $u$ satisfies the integral equation. It remains to prove uniqueness. Suppose $\tilde{u}\in C([0,T]\times \bar{D})$ also satisfies the integral equation. Then $w=u-\tilde{u}\in C([0,T]\times \bar{D})$ solves
\begin{equation}
w(t,x)=-\dfrac{1}{2}\int_0^t \int_{\partial D}p(t-r,x,y)g(r,y)w(r,y)\,\sigma(dy)dr.
\end{equation}
By a Gronwall type argument and (\ref{E:Surface_integral_boundedness}), we can show that $w=0$. More precisely, let $\psi(s)=\sup_{x\in\bar{D}}|w(s,x)|$. Then
$$0\leq \psi(t)\leq \int_0^t \psi(r)\left(\frac{A}{\sqrt{t-r}}+B\right)\,dr \quad\forall\,t\geq 0.$$
Note that
$$\int_0^t \psi(r)\left(\frac{A}{\sqrt{t-r}}+B\right)\,dr=
\frac{\partial}{\partial t}\int_0^t \psi(r)\left(2A\sqrt{t-r}+Bt\right)\,dr - \psi(t)Bt.$$
Combining the above two inequalities, we have
$$(1+Bt)\psi(t)\leq \frac{\partial}{\partial t}\int_0^t \psi(r)\left(2A\sqrt{t-r}+Bt\right)\,dr. $$
Integrating both sides with respect to  $t$ on the interval $[0,t_0]$, we have
$$0\leq \int_0^{t_0} (1+Bt)\psi(t)\,dt \leq \int_0^{t_0} \psi(r)\left(2A\sqrt{t_0-r}+Bt_0\right)\,dr.$$
From this we have $\psi=0$ on $[0,t_0]$, where $t_0>0$ is small enough so that $2A\sqrt{t_0}+Bt_0<1$. Let $\tilde{\psi}(t)=\psi(t+t_0)$, we can show that $$0\leq \tilde{\psi}(t)\leq \int_0^t \tilde{\psi}(r)\left(\frac{A}{\sqrt{t-r}}+B\right)\,dr \;\text{ for all }t\geq 0$$
and repeat the argument to obtain $\tilde{\psi}=0$ on $[0,t_0]$ (i.e. $\psi=0$ on $[0,2t_0]$). Inductively, we obtain $\psi=0$ on $[0,T]$.
\end{pf}

Now we come to our coupled equation.

\begin{prop}\label{prop:MildSol_CoupledPDE}
For $T>0$, consider the Banach space $\Lambda_T=C([0,T]\times\bar{D}_+)\times C([0,T]\times\bar{D}_-)$ with norm $\|(u,v)\|:= \|u\|_{\infty}+\|v\|_{\infty}$. Suppose $u_+(0)=f\in C(\bar{D}_+)$ and $u_-(0)=g\in C(\bar{D}_-)$. Then there is a unique element $(u_+,u_-)\in\Lambda_T$ which satisfies the coupled integral equation
 \begin{align}\label{E:IntegralRep_CoupledPDE}
    \begin{cases}
    u_+(t,x)=P^+_tf(x)-\frac{\lambda}{2}\int_0^t \int_{I}p^+(t-r,x,z)[u_+(r,z)u_-(r,z)]d\sigma(z)\,dr\\
    u_-(t,y)=P^-_tg(y)-\frac{\lambda}{2}\int_0^t \int_{I}p^-(t-r,y,z)[u_+(r,z)u_-(r,z)]d\sigma(z)\,dr.
    \end{cases}
 \end{align}
Moreover, $(u_+,u_-)$ satisfies
 \begin{align}\label{E:ProbabilisticRep_CoupledPDE}
    \begin{cases}
        u_+(t,x)= \E^{x} \big[f(X^+_t)e^{-\lambda\int^t_0u_-(t-s,X^+_s)dL^+_s}\big]\\
        u_-(t,y)= \E^{y}\big[g(X^-_t)e^{-\lambda\int^t_0u_+(t-s,X^-_s)dL^-_s}\big],
    \end{cases}
 \end{align}
where $L^{\pm}$ is the boundary local time of $X^{\pm}$ on the interface $I$.
\end{prop}

Functions $(u_+, u_-)$ satisfying equation \eqref{E:IntegralRep_CoupledPDE}
will be called a weak solution of  \eqref{E:coupledpde:+} and \eqref{E:coupledpde:-},
as it can be shown that they are weakly differentiable and solve the equations in the distributional sense.

\medskip

\noindent {\it Proof of Proposition \ref{prop:MildSol_CoupledPDE}}.
Define the operator $S$ on $\Lambda_T$ by
$S(u,v)=(S^+v,S^-u)$, where
$$
S^+v(t,x)=\E^{x}\left[f(X^+_t)e^{-\lambda\int^t_0v(t-s,X^+_s)dL^+_s}\right]
 \quad \hbox{for }  (t,x)\in [0,T]\times\bar{D}_+ ,
$$
$$S^-u(t,y)=\E^{y} \left[g(X^-_t)e^{-\lambda\int^t_0u(t-s,X^-_s)dL^-_s} \right]
 \quad \hbox{for }  (t,y)\in [0,T]\times\bar{D}_- .
$$
Lemma \ref{L:properties_u} implies that $S$ maps into $\Lambda_T$. Moreover, for $ (t,x)\in [0,T]\times\bar{D}_+$,
\begin{eqnarray*}
| (S^+v_1-S^+v_2)(t,x) |  &=& \left| \E^{x} \left[f(X^+_t)\left(e^{-\lambda\int^t_0v_1(t-s,X_s)dL^+_s}-e^{-\lambda\int^t_0v_2(t-s,X^+_s)dL^+_s}\right) \right] \right| \\
&\leq& \|f\|_{\infty} \E^{x}\left[ \left|\lambda\int^t_0v_1(t-s,X_s)dL^+_s-\lambda\int^t_0v_2(t-s,X^+_s)\,dL^+_s \right| \right] \\
&=&     \|f\|_{\infty}\lambda\, \E^{x} \left[\int^t_0|v_1(t-s,X^+_s)-v_2(t-s,X^+_s)|\,\,dL^+_s
\right]  \\
&\leq&  \lambda|\,\|f\|_{\infty}\,\|v_1-v_2\|_{\infty}\,\E^{x} [L^+_t]\\
&=&     \lambda\,\|f\|_{\infty}\, \|v_1-v_2\|_{\infty}\,\frac{1}{2}\int^t_0\int_{I}p^+(s,x,y)\sigma(dy)ds\\
&\leq&  C_1\lambda\, \sqrt{T}\|f\|_{\infty} \, \|v_1-v_2\|_{\infty} .
\end{eqnarray*}

A similar result holds for $S^-u_1-S^-u_2$. Hence,
\begin{eqnarray*}
\|S(u_1,v_1)-S(u_2,v_2)\|_\infty
&=& \|S^+v_1-S^+v_2\|_{\infty} + \|S^-u_1-S^-u_2\|_{\infty} \\
&\leq& C_1\lambda\, \sqrt{T}\|u_0\|_{\infty} \, \|v_1-v_2\|_{\infty} + C_2\lambda\,\sqrt{T}\|v_0\|_{\infty} \, \|u_1-u_2\|_{\infty}\\
&\leq& \gamma \|(u_1,v_1)-(u_2,v_2)\|
\end{eqnarray*}
for some $\gamma<1$ when $T$  is small enough.

Hence there is a $T_0>0$ such that $S:\, \Lambda_{T_0}\rightarrow \Lambda_{T_0}$ is a contraction map. By Banach fixed point theorem, there is a unique element $(u^{\star},v^{\star}) \in \Lambda_{T_0}$ such that $(u^{\star},v^{\star})=S(u^{\star},v^{\star})$. By Proposition \ref{prop:MildSol_RobinPDE}, $(u^{\star},v^{\star})$ is the unique weak solution to the coupled PDE on $[0,T_0]$.

Repeat the above argument, with $u_0(\cdot)$ replaced by $u^{\star}(T_0,\cdot)$, and $v_0(\cdot)$ replaced by $v^{\star}(T_0,\cdot)$. We see that, since $\|u^{\star}(T_0,\cdot)\|_{\infty} \leq \|u_0\|_{\infty}$, $\|v^{\star}(T_0,\cdot)\|_{\infty} \leq \|v_0\|_{\infty}$ and $C_i \, (i=1,2)$ are the same, we can extend the solution of the coupled PDE uniquely to $[T_0,2T_0]$. Iterating the argument, we have for any $T>0$, the coupled PDE has a unique weak solution in $\Lambda_T$. Invoke Proposition \ref{prop:MildSol_RobinPDE} once more, we obtain the desired implicit probabilistic representation (\ref{E:ProbabilisticRep_CoupledPDE}).

Finally, by using Markov property as in the proof of Proposition \ref{prop:MildSol_RobinPDE}, we see that (\ref{E:ProbabilisticRep_CoupledPDE}) and (\ref{E:IntegralRep_CoupledPDE}) are equivalent.   \qed

\subsection{Main results (rigorous statements)}

In this paper, we always assume the scaling $N\eps^d=1$ holds, so that the interacting random walk model is parameterized by a single parameter $N$ which is the initial number of particles in each of $D^{\eps}_+$ and $D^{\eps}_-$.
More precisely, for each fixed $N$, we set $\eps= N^{-1/d}$ and let $(\eta^{\eps}_t)_{t\geq 0}$ be a Markov process having generator $\mathfrak{L}^{\eps}$ defined in \eqref{E:generator0} and having initial distribution satisfies $\sum_{x\in D^\eps_+} \eta^{\eps,+}_0(x)= \sum_{y\in D^\eps_-} \eta^{\eps,-}_0(y)=N$. We define the empirical measures
$$
\X^{N,\pm}_t(dz) := \dfrac{1}{N}\sum_{x \in D^{\eps}_{\pm}}\eta^{\eps,\pm}_t(x)\1_x(dz).
$$
It is clear that $(\X^{N,+}_t, \X^{N,-}_t)_{t\geq 0}$ is a continuous time Markov process (inheriting from that of $\eta_t$) with state space
$$\mathfrak{E}:= M_{\leq 1}(\bar{D}_{+}) \times M_{\leq 1}(\bar{D}_{-}),$$ where
$M_{\leq 1}(E)$ denotes the space of non-negative Borel measures on $E$ with mass at most 1. $M_{\leq 1}(E)$ is a closed subset of $M_+(E)$, where the latter denotes the space of finite non-negative Borel measures on $E$ equipped with the following metric:
\begin{equation}\label{E:Metric_MD}
\|\mu-\nu\| := \sum_{k=1}^{\infty} \dfrac{1}{2^k}\,\dfrac{|\<\mu,\phi_k\>-\<\nu,\phi_k\>|}{1+|\<\mu,\phi_k\>-\<\nu,\phi_k\>|},
\end{equation}
where $\{\phi_k : k\geq1\}$ is any countable dense subset of $C(E)$. The topology induced by this metric is equivalent to the weak topology (i.e. $\|\mu_n-\mu\|\rightarrow 0$ if and only if $\<\mu_n,\,f\> \rightarrow \<\mu,\,f\>$ for all $f\in C(E)$). Under this metric, $M_+(\bar{D})$ is a complete separable metric space, hence so are $\mathfrak{E}$ and the Skorokhod space $D([0,T], \,\mathfrak{E})$ (see e.g. Theorem 3.5.6 of \cite{EK86}). Here is our first main result.

\begin{thm}\label{T:conjecture} (Hydrodynamic Limit)
    Suppose Assumption \ref{A:Setting} holds and the sequence of initial configurations $\eta^\eps_0$ satisfies the following conditions:
    \begin{enumerate}
    \item[\rm (i)]  $\X^{N,\pm}_0 \toL u^{\pm}_0(z)dz$ in $M_{\leq 1}(\bar{D}_{\pm})$, where $u^{\pm}_0\in C(\bar{D}_{\pm})$. (Note that  $\X^{N,\pm}_0$ has unit mass for all $N$.)
    \item[\rm (ii)] $\varlimsup_{N\to\infty}\sup_{z\in D^{\eps}_{\pm}}\E \left[\left(\eta^{\eps, \pm}_0(z)\right)^2\right]<\infty.$
        \end{enumerate}
    Then for any $T>0$, as $\eps\to 0$ along the sequence $\eps_j=2^{-j}$, we have
    $$(\X^{N,+}, \,\X^{N,-})\toL (\nu^+, \,\nu^-) \in D([0,T], \,\mathfrak{E}),$$
    where $(\nu^+,\,\nu^-)$ is the deterministic element in $C([0,T],\,\mathfrak{E})$ such that
$$
(\nu^+_t(dx), \,\nu^-_t(dy))=(u_+(t,x)\,\rho_+(x) dx,\, u_-(t,y) \,\rho_-(y)dy)
$$
 for all $t\in [0,T]$, and $(u_{+},u_{-})$ is the unique weak  solution of the coupled PDEs (\ref{E:coupledpde:+}) and (\ref{E:coupledpde:-}) with initial value $(u^+_0,\,u^-_0)$.
\end{thm}

Theorem \ref{T:conjecture} gives the limiting probability distribution of one particle randomly
picked in $D^{\eps}_{\pm}$ at time $t$. This is the $1$-particle distribution in the terminology of statistical physics.

{\bf Question}: What is the limiting joint distribution of more than one particles?

Before stating the answer,  we need to introduce a standard tool in the study of stochastic particle systems: the notion of correlation functions\footnote{More precisely, we will be using correlation functions for \emph{unlabeled} particles. We refer the readers to \cite{LX80} for the relation between labeled and unlabeled correlation functions.}. Recall that the state space of $\eta^{\eps}=(\eta^{\eps}_t)_{t\geq 0}$ is $E^{\eps}$ defined in \eqref{e:2.7}. We denote by
$$
\Omega^{\eps}_{n,m} := \Big\{\xi=(\xi^{+}, \xi^{-}) \in E^{\eps}: |\xi^{+}|:= \sum_{x}\xi^{+}(x)=n, \ |\xi^{-}|:= \sum_{y}\xi^{-}(y)=m \Big\}
$$
the set of configurations with $n$ and $m$ particles in $D^{\eps}_+$ and  $D^{\eps}_-$ respectively. We then define $A: E^{\eps}\times E^{\eps} \rightarrow \R$ in such a way that whenever $\xi\in \Omega^{\eps}_{n,m}$,
\begin{equation}\label{Def:A_xi_eta}
A(\xi,\eta) := A^+(\xi^+,\eta^+)A^-(\xi^-,\eta^-) := \prod_{x\in D_+}A_{\xi^+(x)}^{\eta^+(x)} \,\prod_{x\in D_-}A_{\xi^-(x)}^{\eta^-(x)},
\end{equation}
where $n\mapsto A_k^n$ is the Poisson polynomial\footnote{The notation $A_k^n$ is suggested by the fact that $\E[A_k^{\varrho}]=\theta^k$ when $\varrho$ is a Poisson random variable with mean $\theta$.} of order $k$, namely $A_0^n:= 1$ and $A_k^n :=n(n-1)\cdots (n-k+1)$ for $k\geq 1$ (in particular, $A_k^n=0$ for $k>n$).
Note that $A_k^n$ is the number of permutations of $k$ objects chosen from  $n$ distinct objects.
So $A(\xi , \eta)$ is the total number of possible site to site pairings between labeled particles
having configuration $\xi$ with a subset of labeled particles having  configuration $\eta$.
An alternative representation of (\ref{Def:A_xi_eta}) will be given  in (\ref{E:CombinatoricCorrelation}).

\textbf{Convention: }
For $(\vec{r},\vec{s})\in (D^\eps_+)^n \times (D^\eps_-)^m$ and $\eta \in \Omega^\eps_{N, M}$,  we define
$A((\vec{r},\vec{s}),\,\eta)$ to be $A(\xi, \eta)$ with
$\xi=(\sum_{i}\delta_{r_i}, \sum_{j}\delta_{s_j})$.

\begin{definition}\label{D:2.19}
Let $\P^{\eta}$ is the law of a process with generator $\mathfrak{L}^{\eps}$ and initial distribution $\eta$ satisfying
\begin{equation}\label{e:2.31}
\sum_{x\in D^\eps_+} \eta^{ +} (x) = \sum_{y\in D^\eps_-} \eta^{-} (y)=\eps^{-d}.
 \end{equation}
 For all $t\geq 0$, we define
\begin{equation}\label{e:2.27}
\gamma^{\eps}(\xi,t):= \gamma^{\eps,\,(n,m)}(\xi,t) := \dfrac{\eps^{d(n+m)}}{\alpha_{\eps}(\xi)}\,\E^{\eta}[A(\xi,\eta_t)]
\end{equation}
for all $\xi\in \Omega^{\eps}_{n,m}$, where
\begin{equation}\label{e:2.29}
\alpha_{\eps}(\xi):=m_{\eps}(\vec{r},\vec{s}) := \prod_{i=1}^n m^+_{\eps}(r_i)\,\prod_{j=1}^m m^{-}_{\eps}(s_j).
\end{equation}
when $\xi=(\sum_{i}\delta_{r_i}, \sum_{j}\delta_{s_j})$. By convention, we also have $\gamma^{\eps}((\vec{r},\vec{s}),t):=\gamma^{\eps}(\xi,t)$. Note that $\gamma^{\eps}$ depends on the initial configuration of $\eta$.
\end{definition}

Intuitively, suppose we randomly pick $n$ and $m$ living particles in $D_+$ and $D_-$ respectively at time $t$, then $(\vec{r},\vec{s})\mapsto \gamma^{\eps,(n,m)}((\vec{r},\vec{s}),\,t)$ is the joint probability density function for their positions, up to a normalizing constant.
Therefore, it is natural that $\gamma^{\eps,(n,m)}$ defined by \eqref{e:2.27} is called the $(n,m)$-{\bf particle correlation function}.

The next is our second main result, \textbf{propagation of chaos}, for our system.
It says that when the number of particles tends to infinity, they appears to be independent of each other.
Mathematically, the correlation function factors out in the limit $N\to\infty$.

\begin{thm}\label{T:correlation} (Propagation of Chaos)
    Under the same condition as in Theorem \ref{T:conjecture}, for all $n,\,m\in \mathbb{N}$ and any compact interval $[a,b]\subset (0,\infty)$,
    $$\lim_{\eps\to 0} \,\sup_{(\vec{r},\vec{s})\in \bar{D}_+^n\times \bar{D}_-^m }\, \sup_{t\in[a,b]} \, \Big|\gamma^{\eps}((\vec{r},\vec{s}),t)-\prod_{i=1}^{n}u_{+}(t,r_i) \prod_{j=1}^{m}u_{-}(t,s_j)\Big| =0\,,$$
    where $(u_{+},u_{-})$ is the weak solution of the coupled PDE.
\end{thm}

To investigate the intensity of killing near the interface, we  define $J^{N,\pm} \in D([0,\infty),M_+(\bar{D}_{\pm}))$ by
\begin{eqnarray}
J^{N,+}_t(A) &:=& \eps^{d-1}\sum_{z\in I^{\eps}} \Psi(z)\, \eta^+_t(z_+)\eta^-_t(z_-)\,\1_A (z_+) \quad \hbox{for } A\subset \overline D_+,  \label{e:2.34} \\
J^{N,-}_t(B)  &:=&  \eps^{d-1}\sum_{z\in I^{\eps}} \Psi(z)\, \eta^+_t(z_+)\eta^-_t(z_-)\,\1_B (z_-)
\quad \hbox{for } B\subset \overline D_-.  \label{e:2.35}
\end{eqnarray}
Clearly, $\<J^{N, +}_t, 1\>=\<J^{N,-}_t, 1\>$, which measures the
number of encounters of the two types of particles near  $I$.
  An immediately corollary of Theorem \ref{T:correlation} is the following, which is what we need to identify the limit of $(\X^{N,+},\,\X^{N,-})$.
\begin{cor} \label{cor:correlation}
    For any fixed $t\in(0,\infty)$ and $\phi\in C(\bar{D}_{\pm})$,
    \begin{eqnarray*}
    \lim_{N\to\infty}\,\E[\<J^{N,\pm}_t,\phi\>] &=& \frac{1}{2}\,\int_{I}\,u_{+}(t, y)u_{-}(t, y)\,\phi(y)\,\sigma(d y),  \\
    \lim_{N\to\infty}\,\E[(\<J^{N,\pm}_t,\phi\>)^2]&=& \left( \frac{1}{2}\,\int_{I}\,u_{+}(t,y)u_{-}(t, y)\,\phi(y)\,\sigma(dy)
    \right)^2 , \\
    \lim_{N\to\infty}\,\E[\<\mathfrak{X}^{\pm,N}_t,\phi\>] &=& \int_{D_{\pm}}u_{\pm}(t, y)\,\phi(y)\,\rho_{\pm}(y)\,dy,  \\
    \lim_{N\to\infty}\,\E[(\<\mathfrak{X}^{\pm,N}_t,\phi\>)^2] &=& \left(\int_{D_{\pm}}u_{\pm}(t, y)\,\phi(y)\,\rho_{\pm}(y)\,dy\right)^2 .
    \end{eqnarray*}
\end{cor}

\begin{pf}
We only need to apply Theorem \ref{T:correlation} for the cases $(n,m)=(1,1)$ and $(n,m)=(1,0)$. By definition,
$$\gamma^{\eps}(\1_{r},t)=\frac{\eps^d}{m^+(r)}\,\E^{\eta}[\eta^+_t(r)] \quad \text{and} \quad \gamma^{\eps}(\1_{r}+\1_{s},t)=\frac{\eps^{2d}}{m^+(r)m^-(s)}\,\E^{\eta}[\eta^+_t(r)\eta^-_t(s)].$$
Using (\ref{E:Approximate_m(x)}) and Lemma \ref{L:DiscreteApprox_SurfaceMea}, we get the first two equations via Theorem \ref{T:correlation}. Using (\ref{E:Approximate_m(x)}) and the assumption that $\rho_{\pm}\in C(\bar{D}_{\pm})$, we have the last two equations again by
 Theorem \ref{T:correlation}.
\end{pf}

\begin{remark} \rm (Conditions on $\eta_0$)
The  two conditions for the initial configuration $\eta_0$ in Theorem \ref{T:conjecture} are mild and natural. They are satisfied, for example, when each particle has the same random initial distribution $\frac{u^{\pm}_0(z)}{\sum_{D_{\pm}}u^{\pm}_0}$. Condition (ii) guarantees that, asymptotically, there is no "blow up" of number of particles at any site. More precisely, this technical condition is imposed so that we have
\begin{equation}\label{E:Jump_SecondMoment}
\sup_{t\geq 0} \E \left[\<1,\,J^{N, +}_t\>^2\right] \leq C<\infty \quad \text{for sufficiently large }N.
\end{equation}
The above can be easily checked by comparing with the process $\bar{\eta}$ that has no annihilation (i.e. $\bar{\eta}$ has generator $\mathfrak{L}^{\eps}_0$). Alternatively, we can use the comparison result (\ref{E:comparison_correlation}) to prove (\ref{E:Jump_SecondMoment}). \qed
\end{remark}

\begin{remark}\label{Rk:GeneralizationResults} \rm
(Generalization) We can generalize our results in a number of ways by the same argument. For example, the initial number of particles in $D_+$ and $D_-$ can be different, the condition $N\eps^d=1$ can be relaxed to $\lim_{N\to \infty}N\eps^d\to 1$ where $\eps$ depends on $N$. The annihilation constant $\lambda$ can be replaced by a space and time dependent function $\lambda(t,x)\in C([0,\infty)\times I)$. The diffusion coefficients in $D_{+}$ and $D_{-}$ can be different. The condition ``$\X^{N,\pm}_0$ has mass one for all $N$'' can be replaced by ``the mass of $\X^{N,\pm}_0$ is uniformly bounded in $N$''.  More generally, the same method can be extended to deal with similar models with more than two types of particles. \qed
\end{remark}

The remaining part of this paper is devoted to the proof of Theorem \ref{T:conjecture} and Theorem \ref{T:correlation}. We first prove Theorem \ref{T:correlation} because the proof of Theorem \ref{T:conjecture} relies on Theorem \ref{T:correlation}.

\section{Propagation of Chaos}

\subsection{Duality}

The starting point of our analysis is the discrete integral equation for $\gamma^{\eps}$ in Lemma \ref{L:correlation_NotSelfDual}. At the heart of its proof is the dual relation in Lemma \ref{L:Duality}, which says that the two independent processes $\xi^{0}=(\xi^{0}_t)_{t\geq 0}$ and $\eta^{0}=(\eta^{0}_t)_{t\geq 0}$ of independent ransom walks \emph{with no interaction} are dual to each other with respect to the function $\frac{A(\xi,\eta)}{\alpha_{\eps}(\xi)}$, where $\xi,\,\eta \in E^{\eps}$. Such kind of dual formula for the whole grid $\Z^d$ appeared in \cite{BDPP87} and in Chapter 15 of \cite{mfChen03}.

\begin{lem}\label{L:Duality}(Duality for independent processes)
Let $\xi^{0}=(\xi^{0}_t)_{t\geq 0}$ and $\eta^{0}=(\eta^{0}_t)_{t\geq 0}$ be independent continuous time Markov processes on $E^{\eps}$ with generator $\mathfrak{L}^{\eps}_0$ defined in Definition \ref{Def:ConfigurationProcess_eta}. Then we have
\begin{equation}\label{E:Duality}
\E\left[\dfrac{A(\xi^{0}_{t},\eta^{0}_0)}{\alpha_{\eps}(\xi^{0}_{t})}\right]=\E\left[\dfrac{A(\xi^{0}_{0},\eta^{0}_t)}{\alpha_{\eps}(\xi^{0}_{0})}\right]\quad \text{for every }t\geq 0.
\end{equation}
\end{lem}

\begin{pf}
Assume $\xi^0_0\in \Omega^{\eps}_{n,m}$ and $\eta^0_0\in \Omega^{\eps}_{N,M}$. Then we have $\xi^0_t\in \Omega^{\eps}_{n,m}$ and $\eta^0_t\in \Omega^{\eps}_{N,M}$ for all $t\geq 0$. Without loss of generality, we may assume $N\geq n\geq 1$ and $M\geq m\geq 1$ as otherwise both sides inside expectations of (\ref{E:Duality}) are zero by the definition of $A(\xi,\eta)$.

Denote $U$  the map that sends $(\vec{r},\vec{s})\in (D^{\eps}_{+})^n \times (D^{\eps}_{-})^m$ to $(\sum_{i}\delta_{r_i},\,\sum_{j}\delta_{s_j})\in\Omega^{\eps}_{n,m}$ for any $(n,m)$. We first focus on $D_+$ in Step 1 and Step 2 below.

\textbf{Step 1. } For any $\vec{r}\in (D^{\eps}_+)^n$ and $\eta^+ \in \Omega^{\eps}_{N,0}$, fix some  $\vec{x}^+=(x^+_1, \dots, x^+_N)\in U^{-1} (\eta^+)$.
Then by the definition \eqref{Def:A_xi_eta} of   $A$,
\begin{eqnarray}\label{E:Correlation_equivalentDef}
A^{+}(\vec{r},\eta^+)=  \sharp\, \left\{    \vec{i}:
 \vec{x}^{+}_{\vec{i}}=\vec{r} \right\},
\end{eqnarray}
where   $n$-tuples $\vec{i}:=(i_1,\cdots,\,i_n)$ consist of distinct positive integers in the set $\{1,\,2,\,\cdots,\,N\}$, $\vec{x}^{+}_{\vec{i}}:=(x^+_{i_i},\cdots,\,x^+_{i_n})$ and $\sharp S$ denotes the number of elements in the finite set $S$.

\textbf{Step 2. }
Denote by $\P^{\eta^+}_0$ the law of the \emph{unlabeled} process $(\eta^{0}_t)_{t\geq 0}$ starting from $\eta^+ \in \Omega^{\eps}_{N,0}$ and has generator $\mathfrak{L}^{\eps}_0$.
 Let $\vec{x}^+=(x^+_{1},\cdots,\,x^+_{N})\in U^{-1}(\eta^+)$,
and $\vec{X}^{+,\eps}_t:=(X^{+,\eps}_1(t),\,\cdots,\,X^{+,\eps}_N (t))$ be
independent CTRWs in $D^\eps_+$  starting from $\vec{x}$,
whose law will be denoted as  $\P^{\vec{x}^+}$.
Then by (\ref{E:Correlation_equivalentDef}), we have
\begin{eqnarray}\label{E:Duality_2}
    \E^{\eta^+}_0[A(\vec{r},\,\eta^{0}_t)]
=  \E[\,\sharp\, \{ \hbox{$n$-tuples } \vec{i}:
\,  \vec{X}^{+, \eps}_{\vec{i}}(t)=\vec{r}\}\,]
    = \sum_{\vec{i}:  \hbox{ $n$-tuples}}
\P^{\vec{x}^+_{\vec{i}}} (\vec{X}^{+,\eps}_{\vec{i}} (t)= \vec{r}).
\end{eqnarray}
where $\P^{\vec{x}^+_{\vec{i}}}$ is the law of $\{ \vec{X}^{+,\eps}_{\vec{i}} (t); t\geq 0\}$.
Denote by $p^\eps (\theta, \vec{z}. \vec{w})$ the transition density of $n$ independent
CTRWs in $D^\eps_+$.   By Chapman-Kolmogorov equation, we have for any $\theta\in[0,t]$,
$$\P^{\vec{x}^+_{\vec{i}}}\,( \vec{X}^{+,\eps}(t)= \vec{r}\,)
=\sum_{\vec{z}\in (D^{\eps}_{+})^n}p^{\eps}(\theta,\vec{x}^+_{\vec{i}},\vec{z})p^{\eps}(t-\theta,\vec{z},\vec{r})\,m(\vec{z})\,m(\vec{r}).$$
Putting this into (\ref{E:Duality_2}), we have
\begin{eqnarray}\label{E:Duality_3}
    \E^{\eta^+}_0[A(\vec{r},\,\eta^{0}_t)] &=&
    m_\eps (\vec{r})\,\sum_{\vec{z}}\sum_{\vec{i}}\P^{\vec{x}^{+,\eps}_{\vec{i}}}\,
    (\vec{X}^{+}_{\vec{i}}(\theta)= \vec{z}\,) \,p^{\eps}(t-\theta,\vec{z},\vec{r}) \nonumber \\
    &=&  m_\eps (\vec{r})\,\sum_{\vec{z}}\E^{\eta^+}_0[A(\vec{z},\,\eta^{0}_{\theta})]\,
    p^{\eps}(t-\theta,\vec{z},\vec{r})\quad\text{by }(\ref{E:Duality_2})\text{  again}
     \nonumber \\
    &=&  m_\eps (\vec{r})\,\sum_{\vec{z}}\E^{\eta^+}_0[A(\vec{z},\,\eta^{0}_{\theta})]\,
    p^{\eps}(t-\theta,\vec{r},\vec{z})\quad\text{by symmetry of }p^{\eps} \nonumber \\
    &=&  m_\eps (\vec{r})\,\E\,\left[\dfrac{A(\vec{Y}^{+,\eps}_{t-\theta} ,\eta^{0}_{\theta})}{m_\eps (\vec{Y}^{+,\eps}_{t-\theta})}\right],
\end{eqnarray}
where $\E$ is the expectation corresponding the probability measure under which each coordinate processes of  $\{\vec{Y}^{+,\eps}_t;  t\geq 0\}$ are independent CTRWs with $\vec{Y}^{+,\eps}_0=\vec{r}$ and
are independent of  $(\eta^{0}_t)_{t\geq 0}$.

\textbf{Step 3. }Now we work on $D_+\times D_-$. For any $(\vec{r},\vec{s})\in (D^{\eps}_{+})^n \times (D^{\eps}_{-})^m$ and $\eta=(\eta^+,\eta^-)\in \Omega^{\eps}_{N,M}$, take $\vec{x}= (x^+_{1},\cdots,\,x^+_{N},\,x^-_{1},\cdots,\,x^-_{M})\in U^{-1}(\eta)$.
As in step 1, we have
\begin{eqnarray}\label{E:CombinatoricCorrelation}
   A^{+}((\vec{r},\vec{s}),\,\eta ) =
 \sharp\, \left\{  (\vec{i}, \vec{j}):  ( \vec{x}^{+}_{\vec{j}},   \vec{x}^{-}_{\vec{j}})=(\vec{r},\vec{s}) \right\},
\end{eqnarray}
where $\vec{i}$ runs over all $n$-tuples $\vec{i}:=(i_1,\cdots,\,i_n)$ consisting of distinct positive integers in the set $\{1,\,2,\,\cdots,\,N\}$, and $\vec{j}$  over all $m$-tuples $\vec{j}:=(j_1,\cdots,\,j_m)$ consisting of distinct positive integers in the set $\{1,\,2,\,\cdots,\,M\}$.

Denote by $\P^{\eta}_0$ the law of the \emph{unlabeled} process $(\eta^{0}_t)_{t\geq 0}$ starting from $\eta \in \Omega^{\eps}_{N,M}$ and has generator $\mathfrak{L}^{\eps}_0$. Since all processes on $D^{\eps}_+$ are independent of those on $D^{\eps}_-$, we can proceed as in step 2 (via (\ref{E:CombinatoricCorrelation})) to obtain
\begin{equation}\label{e:3.7}
 \E^{\eta}_0[A((\vec{r},\vec{s}),\,\eta^{0}_t)]= m_\eps (\vec{r},\vec{s})\,\E\,\left[\dfrac{A(\vec{Y}^{\eps}_{t-\theta},\eta^{0}_{\theta})}{m_\eps(\vec{Y}^{\eps}_{t-\theta})}\right]
\quad \hbox{for } \theta\in[0,t],
\end{equation}
where $\vec{Y}^{\eps}:=(Y^{+,\eps}_1,\,\cdots,\,Y^{+,\eps}_n,\,Y^{-,\eps}_1,\,\cdots,\,Y^{-,\eps}_m)$ is independent of $\eta^0$ with $\vec{Y}^{\eps}_0=(\vec{r}, \vec{s})$,
and its  components are mutually independent CTRWs on $D^\eps_\pm$, respectively.
The proof is now complete by taking $\theta=0$.
\end{pf}

We now formulate the discrete integral equations that we need. Recall the definition of $\mathfrak{K}^{\eps}$ from \eqref{e:2.10}  and the definition of
$\gamma^{\eps}((\vec{r},\vec{s}),t)$
from \eqref{e:2.27}.

\begin{lem}\label{L:correlation_NotSelfDual}(Discrete integral equation for $\gamma^{\eps}$)
    For any $\eps>0$, $t>0$, $(\vec{r},\vec{s})\in (D^{\eps}_+)^n\times (D^{\eps}_-)^m$, non-negative integers $n,\,m$ and initial distribution $\eta_0$, we have
    \begin{eqnarray}
    \gamma^{\eps}((\vec{r},\vec{s}),t)
    &=& \sum_{(\vec{r'},\vec{s'})}\, \gamma^{\eps}((\vec{r'},\vec{s'}),0)\, p^{\eps}(t,(\vec{r},\vec{s}),(\vec{r'},\vec{s'}))\,m(\vec{r'},\vec{s'}) \nonumber \\
    && \quad +  \int_0^t \sum_{(\vec{r'},\vec{s'})}\,p^{\eps}(t-s,(\vec{r},\vec{s}), (\vec{r'},\vec{s'}))\,\E[\mathfrak{K}^{\eps}A((\vec{r'},\vec{s'}),\eta_s)]\,\eps^{d(n+m)}\,ds, \label{e:3.7}
    \end{eqnarray}
    where $\mathfrak{K}^{\eps}$ acts on the $\eta$-variable of $A((\vec{r},\vec{s}),\eta)$.
\end{lem}

\begin{pf}
Starting from (\ref{E:Duality}), we can obtain Lemma \ref{L:correlation_NotSelfDual} by `integration by parts' as follows.

Let $\P_{(\xi^0)}$ and $\P_{(\eta^0)}$ be the laws of $\xi^0$ and $\eta^0$ respectively. (\ref{E:Duality}) is equivalent to saying that for any $\xi$ and $\eta$, we have
\begin{equation}\label{E:correlation_NotSelfDual_1}
\E_{(\xi^0)}\left[\dfrac{A(\xi^{0}_{w},\eta)}{\alpha_{\eps}(\xi^{0}_{w})}\Big|\,\xi^0_0=\xi \right]=
\E_{(\eta^0)}\left[\dfrac{A(\xi,\eta^{0}_w)}{\alpha_{\eps}(\xi)}\Big|\,\eta^0_0=\eta\right]\quad \text{for every }w\geq 0.
\end{equation}
Taking $w=t-s$, we see that (\ref{E:correlation_NotSelfDual_1}) is in turn equivalent to
\begin{equation}\label{E:correlation_NotSelfDual_2}
F^{(\xi)}_s(\eta):= P^{(\xi^0)}_{t-s}\left( \dfrac{A(\cdot,\eta)}{\alpha_{\eps}(\cdot)}\right)(\xi)
=P^{(\eta^0)}_{t-s}\left( \dfrac{A(\xi,\cdot)}{\alpha_{\eps}(\xi)}\right)(\eta)=: G^{(\eta)}_s(\xi)
\quad \text{for every }s\in[0,t] \text{ and }t\geq 0,
\end{equation}
where $P^{(\xi^0)}_{t}$ and $P^{(\eta^0)}_{t}$ are the transition semigroup  of $\xi^0$ and $\eta^0$, respectively, and they act on the $\xi$ and $\eta$ variables in $\frac{A(\xi,\eta)}{\alpha_{\eps}(\xi)}$, respectively. Therefore, with $\mathfrak{L}^{\eps}_0$ acting on the $\eta$ variable, we have
\begin{equation}\label{e:3.9}
    \frac{\partial}{\partial s}\,F^{(\xi)}_s(\eta) = \frac{\partial}{\partial s}\,G^{(\eta)}_s(\xi)
= -\mathfrak{L}^{\eps}_0\,P^{(\eta^0)}_{t-s}\left( \dfrac{A(\xi,\cdot)}{\alpha_{\eps}(\xi)}\right)(\eta)
= -\mathfrak{L}^{\eps}_0\,F^{(\xi)}_s(\eta)  .
\end{equation}

Recall that $\eta_t$ is the configuration process of our interacting system with generator $\mathfrak{L}^{\eps}_0+ \mathfrak{K}^{\eps}$ (see Definition \ref{Def:ConfigurationProcess_eta}). Fix $\xi$ and consider the function $(s,\eta)\mapsto F_s(\eta):=F^{(\xi)}_s(\eta)$. We have
$$M_s:= F_s(\eta_s)-F_0(\eta_0)-\int_0^s \left(\frac{\partial F_r}{\partial r}+\mathfrak{L}^{\eps}_0F_r+ \mathfrak{K}^{\eps}F_r\right)(\eta_r)\,dr$$
is a $\F^{\eta}_s$-martingale for $s\in[0,t]$. By
\eqref{e:3.9}
 and the fact that $E^{\eta}[M_t]=E^{\eta}[M_0]=0$, where $\P^{\eta}$ is the law of $(\eta_t)_{t\geq 0}$ starting from $\eta$, we have
\begin{equation*}
0= \E^{\eta}\left[\dfrac{A(\xi,\eta_t)}{\alpha_{\eps}(\xi)}\right]
- P^{(\xi^0)}_{t}\left( \dfrac{A(\cdot,\eta)}{\alpha_{\eps}(\cdot)}\right)(\xi)
- \int_0^t\E^{\eta}\left[\mathfrak{K}^{\eps}\,
P^{(\xi^0)}_{t-r}\left( \dfrac{A(\cdot,\eta_r)}{\alpha_{\eps}(\cdot)}\right)(\xi)
\right] \, dr
\end{equation*}
for all $\xi$ and $\eta$. This is equivalent to the stated equation in the lemma.

\end{pf}

It is clear that $\mathfrak{K}^{\eps}A(\xi,\eta) \leq 0$. Hence, as an immediate consequence of Lemma \ref{L:correlation_NotSelfDual}, we have the following comparison result:
\begin{equation}\label{E:comparison_correlation}
\gamma^{\eps}((\vec{r},\vec{s}),t) \leq \sum_{(\vec{r'},\vec{s'})}\, \gamma^{\eps}((\vec{r'},\vec{s'}),0)\, p^{\eps}(t,(\vec{r},\vec{s}),(\vec{r'},\vec{s'}))\,m(\vec{r'},\vec{s'})
\end{equation}
for all $t>0$ and $(\vec{r},\vec{s})\in (D^{\eps}_{+})^n\times (D^{\eps}_{-})^m$.

\subsection{Annihilation near the interface}

For any $\xi=(\xi^+,\xi^-)\in E^{\eps}$, we let
$\xi^{+}_{(x)}=\xi^{+}(x)\1_{x}$, the element that has only $\xi^{+}(x)$ number of particles at $x$, and none elsewhere. Similarly,  we denote  $\xi^{-}(y)\1_{y}$ by $\xi^{-}_{(y)}$. Set  $\xi_{(x,y)}=(\,\xi^{+}(x)\1_{x} , \xi^{-}(y)\1_{y}\,)$, the element that has only $\xi^{+}(x)$ number of particles at $x$, $\xi^{-}(y)$ number of particles at $y$, and none elsewhere.

\begin{lem} \label{L:correlation_NotSelfDualk}
Let $\mathfrak{K}^{\eps}$ be the operator defined in \eqref{e:2.10} and acts on the $\eta$-variable of $A(\xi, \eta)$. Then
\begin{equation}\label{E:correlationk}
\mathfrak{K}^{\eps}A(\xi,\eta)=\sum_{z\in I^{\eps}}A(\xi-\xi_{(z_+,z_-)},\,\eta)\cdot\mathfrak{K}^{\eps}A(\xi_{(z_+,z_-)},\eta) .
\end{equation}
Moreover, if $\xi\in \Xi:= \{\xi:\;\xi^{\pm}(z_{\pm})\leq 1
\hbox{ for every }  z\in I^{\eps}\}$, then
\begin{eqnarray}
\mathfrak{K}^{\eps}A(\xi,\eta)
&=& -\dfrac{\lambda}{\eps} \, \sum_{\substack {z\in I^{\eps}: \, \xi^{+}(z_+)=1}} \,\Psi_\eps (z)\,A(\xi+\1_{(0, z_-)},\eta)     \label{E:KA(xi,eta)_1}\\
&& -\dfrac{\lambda}{\eps} \, \sum_{\substack {z\in I^{\eps}: \, \xi^{-}(z_-)=1}} \,\Psi_\eps (z)\,A(\xi+\1_{(z_+, 0)},\eta)    \label{E:KA(xi,eta)_2}\\
&& -\dfrac{\lambda}{\eps} \, \sum_{\substack {z\in I^{\eps}: \, \xi(z_+,z_-)=(1,1)}} \,\Psi_\eps (z)\,A(\xi,\eta) .      \label{E:KA(xi,eta)_3}
\end{eqnarray}
\end{lem}

\begin{pf}
Observe that $A(\xi-\xi_{(x,y)},\eta)\,A(\xi_{(x,y)},\eta)=A(\xi,\eta)$. Consequently
\begin{eqnarray*}
&& \frac{\lambda}{\eps} \Psi_\eps (z) \eta^+(z_+) \eta^-(z_- )
\left( A(\xi,\,\eta-\1_{(z_+,z_-)})-A(\xi,\eta)  \right) \\
&=& \frac{\lambda}{\eps} \Psi_\eps (z) \eta^+(z_+) \eta^-(z_- )
A(\xi-\xi_{(z_+,z_-)},\eta) \left( A(\xi_{(z_+,z_-)},\eta-\1_{(z_+,z_-)})-A(\xi_{(z_+,z_-)},\eta)\right) \\
&=& A(\xi-\xi_{(z_+,z_-)},\eta)  \mathfrak{K}^{\eps}A(\xi_{(z_+,z_-)},\eta).
\end{eqnarray*}
Thus \eqref{E:correlationk} holds.   On other hand,
$$
   \mathfrak{K}^{\eps}A(\xi_{(z_+,z_-)},\eta)
=\Psi_\eps (z)\,\dfrac{\lambda}{\eps}\,\eta^{+}(z_+)\eta^{-}(z_-) \times
            \begin{cases}
              -1,   &\text{if  }  \xi(z_+,z_-)=(1,0) \text{ or } (0,1)  \\
              1-\eta^{+}(z_+)-\eta^{-}(z_-),    &\text{if  }  \xi(z_+,z_-)=(1,1)
            \end{cases}.
$$
Observe also that for $x\in D^\eps_+$ and $y\in D^\eps_-$,
$$
A \left(\xi-\xi_{(x,y)}, \,\eta\right)\,\eta^{+}(x)\eta^{-}(y)=A \left(\xi-\xi_{(x,y)}+ \1_{(x, y)}, \,\eta \right)
$$
and
\begin{eqnarray*}
&& A(\xi-\xi_{(x,y)})\,\eta^{+}(x)^2 \, \eta^{-}(y)\\
&=& A \left(\xi-\xi_{(x,y)}, \eta \right)\, \left( \eta^{+}(x)^2-\eta^{+}(x)+\eta^{+}(x) \right)\,\eta^{-}(y)\\
&=& A(\xi-\xi_{(x,y)}, \eta )\,  A (2\1_x, \eta^+ ) \,\eta^{-}(y) +A(\xi-\xi_{(x,y)}, \eta )\, \eta^{+}(x) \,\eta^{-}(y)\\
&=& A( \xi -\xi_{(x, y)}+2\1_{(x, 0)}+\1_{(0,y)}, \eta) + A(\xi, \eta ).
\end{eqnarray*}
Similarly,
$$
A(\xi-\xi_{(x,y)}, \eta)\,\eta^{+}(x)\eta^{-}(y)^2 =
A( \xi -\xi_{(x, y)}+\1_{(x, 0)}+2 \1_{(0,y)}, \eta) + A(\xi, \eta ).
$$
From the above calculations and (\ref{E:correlationk}), we see that for $\xi \in \Xi$,
\begin{eqnarray*}
\mathfrak{K}^{\eps}A(\xi,\eta)
&=&   -\dfrac{\lambda}{\eps} \, \sum_{\substack {z\in I^{\eps}: \, \xi(z_+,z_-)=(1,0)}} \,\Psi_\eps (z)\,A(\xi +\1_{(0, z_-)}, \, \eta)  \\
&& -\dfrac{\lambda}{\eps} \, \sum_{\substack {z\in I^{\eps}: \, \xi(z_+,z_-)=(0,1)}} \,\Psi_\eps (z)\,A(\xi +\1_{(z_+, 0)}, \, \eta) \\
&& -\dfrac{\lambda}{\eps} \, \sum_{\substack {z\in I^{\eps}:  \, \xi(z_+,z_-)=(1,1)}} \,\Psi_\eps (z)
\left( A(\xi,\eta)+A(\xi+\1_{(0, z_-)},\eta)
+A(\xi + \1_{(z_+, 0)}, \eta)\right),
\end{eqnarray*}
which gives the desired formula.
\end{pf}

\subsection{Uniform bound and equi-continuity}

We extend to define $\gamma^{\eps,(n,m)}(\cdot,t)$ continuously on $\bar{D}_{+}^n \times \bar{D}_{-}^m$ while preserving the supremum and the infinmum in each small $\eps$-cube. We can accomplish this by the interpolation described in \cite{BK08} or \cite{SZ97}, or by a sequence of harmonic extensions along simplexes with increasing dimensions (described in  \cite{wtF14}).
Recall that the definition of $\gamma^{\eps,(n,m)}(\cdot,t)$ depends on the initial configuration
$\eta_0$ of the interacting random walks (see Definition \ref{D:2.19}), which  has the normalization \eqref{e:2.31}.

\begin{thm}\label{T:equicts}
    There exists $\eps_0>0$ such that for any $(n,m)\in \mathbb{N}\times \mathbb{N}$, the family of functions $\{\gamma^{\eps}((\vec{r},\vec{s}),t)\}_{\eps \in(0,\eps_0)}$ is uniformly bounded and equi-continuous on $\bar{D}_{+}^n \times \bar{D}_{-}^m \times (0,\infty)$, which is uniform in the initial configuration $\eta_0$ that satisfies \eqref{e:2.31}.
\end{thm}

\begin{pf}
We first prove uniform boundedness. By (\ref{E:comparison_correlation}) and the Gaussian upper bound in Theorem \ref{T:UpperHKE}, we have
\begin{eqnarray*}
\gamma^{\eps}((\vec{r},\vec{s}),t) &\leq& \sum_{(\vec{r'},\vec{s'})}\, \gamma^{\eps}((\vec{r'},\vec{s'}),0)\, p^{\eps}(t,(\vec{r},\vec{s}),(\vec{r'},\vec{s'}))\,m(\vec{r'},\vec{s'})\\
&\leq& \left(\frac{C}{t^{d/2}}\right)^{n+m}\,   \sum_{(\vec{r'},\vec{s'}) \in \bar{D}^n_+\times \bar{D}^m_-}\,
 A \left( (\vec{r'},\vec{s'}), \eta_0\right)\,\eps^{d(n+m)} \quad\text{whenever } \eps\in(0,\eps_0) .
\end{eqnarray*}
Since the initial distribution $\eta_0=(\eta^+_0,\eta^-_0)$ has the property
 that $\sum_{x\in D^\eps_+} \eta^+_0 (x)= \sum_{y\in D^\eps_-} \eta^-_0 (y)=
\eps^{-d}$,  we have
\begin{eqnarray}\label{E:UniformBound_eta0}
&& \sum_{(\vec{r},\vec{s}) \in \bar{D}^n_+\times \bar{D}^m_-}A((\vec{r},\vec{s}),\eta_0)\eps^{d(n+m)} \nonumber \\
&=& \left(\sum_{\vec{r}\in \bar{D}^n_+}A^+(\vec{r},\eta^+_0)\right)\left(\sum_{\vec{s}\in \bar{D}^m_-}A^-(\vec{s},\eta^-_0)\right)\,\eps^{d(n+m)} \notag\\
&\leq& \left(\sum_{\vec{r}\in \bar{D}^n_+}\prod_{i=1}^n\eta^+_0 (r_i)\right)\left(\sum_{\vec{s}\in
\bar{D}^m_-}\prod_{j=1}^m\eta_0^-(s_j)\right)\eps^{d(n+m)} \quad\text{ since } A_k^n\leq n^k \notag\\
&\leq & \prod_{i=1}^n\left(\sum_{r_i \in \bar{D}_+}\eta_0^+ (r_i)\eps^d\right)\,\prod_{j=1}^m\left(\sum_{s_j \in \bar{D}_-}\eta_0^-(s_j)\eps^d\right)=1 .
\end{eqnarray}
Thus there exist $\eps_0=\eps_0(d,D,\rho)$ and $C=C(d,D,\rho)>0$ such that for all $t\in(0,\infty)$,  $(n,m)\in \mathbb{N}\times \mathbb{N}$ and $\eps \in (0, \eps_0)$,
\begin{equation}\label{E:uniformbound}
\sup_{\eps\in (0,\eps_0)}\sup_{\xi\in \Omega^{\eps}_{n,m}} \gamma^{\eps}(\xi,t) \leq \left(\frac{C}{t^{d/2}}\right)^{n+m} .
\end{equation}

We next show that both terms on the right hand side of \eqref{e:3.7} are equi-continuous. Recall that we can rewrite the equation
\eqref{e:3.7} as
\begin{equation*}
\gamma^{\eps}((\vec{r},\vec{s}),t)= F^{\eps}((\vec{r},\vec{s}),t)+G^{\eps}((\vec{r},\vec{s}),t) ,
\end{equation*}
where
\begin{eqnarray*}
F^{\eps}((\vec{r},\vec{s}),t) &:=& \sum_{(\vec{r'},\vec{s'})}\, \gamma^{\eps}((\vec{r'},\vec{s'}),0)\, p^{\eps}(t,(\vec{r},\vec{s}),(\vec{r'},\vec{s'}))\,m(\vec{r'},\vec{s'}), \\
G^{\eps}((\vec{r},\vec{s}),t) &:=& \int_0^t \sum_{(\vec{r'},\vec{s'})}\,p^{\eps}(t-s,(\vec{r},\vec{s}), (\vec{r'},\vec{s'}))\,\E[\mathfrak{K}^{\eps}A((\vec{r'},\vec{s'}),\eta_s)]\,\eps^{d(n+m)}\,ds .
\end{eqnarray*}
Now let $(\vec{r},\vec{s})$, $(\vec{p},\vec{q})\in (D^{\eps}_{+})^n \times (D^{\eps}_{-})^m$ and $0<t<\ell\leq \infty$.
For the first term,
\begin{eqnarray*}
&& \Big|F^{\eps}((\vec{r},\vec{s}),t) - F^{\eps}((\vec{p},\vec{q}),\ell)\Big|  \\
&=& \Big|\sum_{(\vec{r'},\vec{s'})}\left(p^{\eps}(t,(\vec{r},\vec{s}),(\vec{r'},\vec{s'}))-p^{\eps}(\ell,(\vec{p},\vec{q}),(\vec{r'},\vec{s'}))\right)
\,\E[A((\vec{r'},\vec{s'}),\eta_0)]\eps^{d(n+m)}\Big|\\
&\leq& \left(\sup_{(\vec{r'},\vec{s'})}\big|p^{\eps}(t,(\vec{r},\vec{s}),(\vec{r'},\vec{s'}))-p^{\eps}(\ell,(\vec{p},\vec{q}),(\vec{r'},\vec{s'}))\big|
\right)\E_0\left[\sum_{(\vec{r'},\vec{s'})}A((\vec{r'},\vec{s'}),\eta_0)\eps^{d(n+m)}\right]\\
&\leq& \sup_{(\vec{r'},\vec{s'})}\big|p^{\eps}(t,(\vec{r},\vec{s}),(\vec{r'},\vec{s'}))-p^{\eps}(\ell,(\vec{p},\vec{q}),(\vec{r'},\vec{s'}))\big| ,
\end{eqnarray*}
where we have used (\ref{E:UniformBound_eta0}) in the last line. By the uniform  H\"older continuity of $p^{\eps}(t,(\vec{r},\vec{s}),(\vec{r'},\vec{s'}))$ (see Theorem \ref{T:HolderCts} below) and the fact that $p^{\pm}(t,x,y)\in C((0,\infty)\times \bar{D}_{\pm}\times\bar{D}_{\pm})$, we see that $\{F^{\eps}\}$ is equi-continuous at $((\vec{r},\vec{s}),t)$.
For the second term, note that
\begin{equation}\label{E:equicts_G}
G^{\eps}((\vec{p},\vec{q}),\ell)- G^{\eps}((\vec{r},\vec{s}),t) = \int_{t}^{\ell}H^{(1)}(s)\,ds + \int_0^{t}H^{(2)}(s)\,ds ,
\end{equation}
where $$H^{(1)}(s) := \sum_{(\vec{r'},\vec{s'})}\,p^{\eps}(\ell-s,(\vec{p},\vec{q}),(\vec{r'},\vec{s'}))\,\E[\mathfrak{K}^{\eps}A((\vec{r'},\vec{s'}),\eta_s)]\,\eps^{d(n+m)}\text{ and }$$
$$H^{(2)}(s) := \sum_{(\vec{r'},\vec{s'})}\,\left[p^{\eps}(\ell-s,(\vec{r},\vec{s}),(\vec{r'},\vec{s'})) -
p^{\eps}(t-s,(\vec{p},\vec{q}),(\vec{r'},\vec{s'}))\right]\,\E[\mathfrak{K}^{\eps}A((\vec{r'},\vec{s'}),\eta_s)]\,\eps^{d(n+m)} .
$$
In the remaining, we will show that $G^{\eps}$ is equi-continuous. We first deal with $H^{(1)}$ in (\ref{E:equicts_G}).

As in (\ref{E:UniformBound_eta0}), we have
\begin{eqnarray}\label{E:equicts_G_1}
    &&\sum_{(\vec{r'},\vec{s'})}\,p^{\eps}(\theta_1,(\vec{p},\vec{q}),(\vec{r'},\vec{s'}))\,A((\vec{r'},\vec{s'}),\eta_{\theta_2})\,\eps^{d(n+m)}\\
    &\leq& \prod_{i=1}^n\left(\sum_{r_i'}p^{\eps}(\theta_1,p_i,r_i')\eta^+_{\theta_2}(r_i')\eps^d\right)\,
    \prod_{j=1}^m\left(\sum_{s_j'}p^{\eps}(\theta_1,q_j,s_j')\eta^-_{\theta_2}(s_j')\eps^d\right). \notag
\end{eqnarray}
On other hand, using (\ref{E:comparison_correlation}), the Chapman Kolmogorov equation and assumption (ii) for $\eta_0$, in this order, we have
\begin{equation}\label{E:equicts_G_2}
\sup_{\theta_1,\,\theta_2>0}\;\sup_{a\in D^{\eps}_+}\,\E\left[\left(\sum_{x\in D^{\eps}_+}p^{\eps}(\theta_1,a,x)\eta^+_{\theta_2}(x)\eps^d\right)^2 \right] \leq C
\end{equation}
for large enough $N$, where $C>0$ is a constant.
\begin{eqnarray*}
&& \sum_{(\vec{r'},\vec{s'})}\,p^{\eps}(\theta_1,(\vec{p},\vec{q}),(\vec{r'},\vec{s'}))\,\Big|\E\left[\mathfrak{K}^{\eps}A((\vec{r'},\vec{s'}),\eta_{\theta_2})\right]\Big|\,\eps^{d(n+m)}\\
&\leq&
\E\left[\frac{\lambda}{\eps}\sum_{z\in I}\eta^{+}_{\theta_2}(z_+)\eta^{-}_{\theta_2}(z_-)\,
\sum_{(\vec{r'},\vec{s'})}\,p^{\eps}(\theta_1,(\vec{p},\vec{q}),(\vec{r'},\vec{s'}))\,2A((\vec{r'},\vec{s'}),\eta_{\theta_2})\,\eps^{d(n+m)}\right]\\
&\leq&
2\E\left[\<1,\,J^N_{\theta_2}\>\,
\prod_{i=1}^n\left(\sum_{r_i'}p^{\eps}(\theta_1,p_i,r_i')\eta^+_{\theta_2}(r_i')\eps^d\right)\,
    \prod_{j=1}^m\left(\sum_{s_j'}p^{\eps}(\theta_1,q_j,s_j')\eta^-_{\theta_2}(s_j')\eps^d\right)
\right]\\
&\leq& C \quad\text{uniformly for }\theta_1>0,\,\theta_2>0,\,(\vec{p},\vec{q})\in (D^{\eps}_{+})^n \times (D^{\eps}_{-})^m \text{ and }\eps>0 \text{ small enough}.
\end{eqnarray*}
We have used (\ref{E:equicts_G_1}) for the second inequality. The last inequality follows from H\"older's inequality, (\ref{E:Jump_SecondMoment}) and (\ref{E:equicts_G_2}). Therefore for any $(n,m)$,
\begin{equation}\label{E:equicts_G_3}
\sup_{\theta_1,\,\theta_2>0}\;\sup_{(\vec{p},\vec{q})\in (D^{\eps}_{+})^n \times (D^{\eps}_{-})^m}\;
\sum_{(\vec{r'},\vec{s'})}\,p^{\eps}(\theta_1,(\vec{p},\vec{q}),(\vec{r'},\vec{s'}))\,\Big|\E\left[\mathfrak{K}^{\eps}A((\vec{r'},\vec{s'}),\eta_{\theta_2})\right]\Big|\,\eps^{d(n+m)}\leq C
\end{equation}
for large enough $N$, where $C>0$ is a constant. Hence $\int_{t}^{\ell}|H^{(1)}(s)|\,ds\leq C(\ell-t)\to 0$ as $\ell\to t$, uniformly for $(\vec{p},\vec{q})$, $s\in(t,\ell)$ and $\eps$ small enough. Finally, we deal with $H^{(2)}$. For any $h\in(0,t)$, we have
\begin{equation*}
\bigg|\int_0^tH^{(2)}(s)\,ds \bigg| \leq \int_0^{t-h}|H^{(2)}(s)|ds+ \int_{t-h}^t|H^{(2)}(s)|ds .
\end{equation*}
By (\ref{E:equicts_G_3}), we have $\int_{t-h}^t|H^{(2)}(s)|ds\leq C\,h$. By the H\"older continuity of $p^{\eps}$ (cf. Theorem \ref{T:HolderCts}),
\begin{eqnarray*}
\int_0^{t-h}|H^{(2)}(s)|ds &\leq& \int_0^{t-h}
\sup_{(\vec{r'},\vec{s'})}\,\Big|p^{\eps}(\ell-s,(\vec{r},\vec{s}),(\vec{r'},\vec{s'})) -
p^{\eps}(t-s,(\vec{p},\vec{q}),(\vec{r'},\vec{s'}))\Big|\,\\
&& \qquad \qquad \E\left[\frac{\lambda}{\eps}\sum_{z\in I}\eta^{+}_{s}(z_+)\eta^{-}_{s}(z_-)\,
\sum_{(\vec{r'},\vec{s'})}\,2A((\vec{r'},\vec{s'}),\eta_s)\,\eps^{d(n+m)}\right]\,ds \\
&\leq& (t-h)\,C\frac{|\ell-t|^{\sigma_1}+ \|(\vec{r},\vec{s})-(\vec{p},\vec{q})\|^{\sigma_2}}{h^{\sigma_3}}\quad\text{for sufficiently small }\eps>0 ,
\end{eqnarray*}
where $\sigma_i\,(i=1,2,3)$ are positive constants. Since $h\in(0,t)$ is arbitrary, we see that
$\bigg|\int_0^tH^{(2)}(s)\,ds \bigg|\to 0$ as $|\ell-t|+\|(\vec{r},\vec{s})-(\vec{p},\vec{q})\|\to 0$, uniformly for small enough $\eps>0$. Hence $G^{\eps}$ is equi-continuous at an arbitrary $((\vec{r},\vec{s}),t)\in (D^{\eps}_{+})^n \times (D^{\eps}_{-})^m\times (0,\infty)$.
\end{pf}

\medskip

From Theorem \ref{T:equicts} and a diagonal selection argument, it follows that for any sequence $\eps_k \to 0$ there is a subsequence along which  $\gamma^{\eps}$ converges on $\bar{D}_{+}^n \times \bar{D}_{-}^m \times (0,T)$, uniformly on the compacts, to some $\gamma^{(n,m)}\in C(\bar{D}_{+}^n \times \bar{D}_{-}^m \times (0,T))$, for every $(n,m)\in \mathbb{N}\times \mathbb{N}$. Our goal is to show that
 $$
\gamma_t^{(n,m)}(\vec{r},\vec{s})=\prod_{i=1}^n u_+(t,r_i) \prod_{j=1}^m u_-(t,s_j) .
$$

We will achieve this by first showing that both $\Gamma=\{\gamma^{(n,m)}\}$ and $\Phi^{(n,m)}_t:= \prod_{i=1}^n u_+(t,r_i) \prod_{j=1}^m u_-(t,s_j)$ satisfy the same an infinite hierarchy of equations, and then establishing uniqueness of the hierarchy.

\subsection{Limiting hierarchy}

Note that $D_{+}^n\times D_{-}^m$ is a bounded Lipschitz domain in $\R^{(n+m)d}$, and that the boundary $\partial(D_{+}^n\times D_{-}^m)$ contains the disjoint union $\cup_{i=1}^n \partial_+^i\,\bigcup\,\cup_{j=1}^m\partial_-^j$ where
\begin{eqnarray}
\partial_+^i &:=&  \left(D_{+}\times\cdots \times (\overset{i^{th}}{\partial D_{+}}\cap I) \times \cdots \times D_{+}\right)\times D_{-}^m ,  \\
\partial_-^j &:=&  D_{+}^n \times \left(D_{-}\times\cdots \times (\overset{j^{th}}{\partial D_{-}}\cap I) \times\cdots \times D_{-}\right) .
\end{eqnarray}

We define the function $\rho=\rho_{(n,m)}:\,D_{+}^n\times D_{-}^m \rightarrow \R$ by $\rho(\vec{r},\vec{s}):= \prod_{i=1}^n\rho_+(r_i)\prod_{j=1}^m\rho_-(s_j)$. We also denote $p(t,(\vec{r},\vec{s}),(\vec{r'},\vec{s'})) := \prod_{i=1}^n p^{+}(t,r_i,r_i')\prod_{j=1}^m p^{-}(t,s_j,s_j')$, where $p^{\pm}$ is the transition density of the
reflected diffusion $ X^{\pm}$ on $\bar{D}_\pm$ with respect to  the measure $\rho_{\pm}(x) dx$. We now characterize the subsequential limits of $\{\gamma^{\eps}\}_{\eps>0}$:

\begin{thm}\label{T:hierarchy}
Let $\eta^\eps_0$ be a sequence of initial configurations that satisfy
\eqref{e:2.31} with $\eps= N^{-1/d}$, and  the conditions of Theorem \ref{T:conjecture};
that is,  their corresponding empirical measures
 $\X^{N,\pm}_0$  converges weakly to $ u^{\pm}_0(z)dz$ in $M_{\leq 1}(\bar{D}_{\pm})$ for some  $u^{\pm}_0\in C(\bar{D}_{\pm})$
 and
\begin{equation}\label{e:3.24}
\varlimsup_{\eps \to 0}  \sup_{z\in D^\eps_\pm} \E \left[ \eta_0^{\eps, \pm}  (z)^2 \right] <\infty.
\end{equation}
 Denote by
  $\Gamma^{\eps}=\{\gamma^{\eps,\,(n,m)}; \,  t\geq 0, n, m \in \mathbb{N}\}$
 the correlation functions for the interacting random walks
with initial configuration $\eta^\eps_0$.
Let $\Gamma=\{\gamma^{(n,m)}_t; \,  t\geq 0, n, m \in \mathbb{N}\}$ be any subsequential limit (as $\eps\to 0$) of $\Gamma^{\eps}=\{\gamma^{\eps,\,(n,m)}; \,
t\geq 0, n, m \in \mathbb{N}\}$.
Then the following infinite system of hierarchical equations holds:
\begin{eqnarray}\label{E:hierarchy_1}
\gamma_t^{(n,m)}(\vec{r},\vec{s})
&=& \int_{D_+^n \times D_-^m}\,\Phi^{(n,m)}(\vec{a},\vec{b})\,p(t,(\vec{r},\vec{s}),(\vec{a},\vec{b}))\,\rho(\vec{a},\vec{b})\, d(\vec{a},\vec{b}) \\
&& - \frac{\lambda}{2}\int_0^t\,\bigg(
\sum_{i=1}^n \int_{\partial_+^i} \gamma_{\theta}^{(n,m+1)}(\vec{a},(\vec{b},a_i)) \,p(t-\theta,(\vec{r},\vec{s}),(\vec{a},\vec{b}))\,\frac{\rho(\vec{a},\vec{b})}{\rho_+(a_i)}\,d\sigma_{(n,m)}(\vec{a},\vec{b}) \notag \\
&& \quad\quad +\sum_{j=1}^m \int_{\partial_-^j} \gamma_{\theta}^{(n+1,m)}((\vec{a},b_j),\vec{b})) \,p(t-\theta,(\vec{r},\vec{s}),(\vec{a},\vec{b}))\,\frac{\rho(\vec{a},\vec{b})}{\rho_-(b_j)}\,d\sigma_{(n,m)}(\vec{a},\vec{b})
\bigg)\, d\theta , \notag
\end{eqnarray}
where $d(\vec{a},\vec{b})$ is the Lebesgue measure on $\R^{n+m}$,  $\sigma_{(n,m)}$ is the surface measure of $\partial(D_{+}^n\times D_{-}^m)$ and $\Phi^{(n,m)}(\vec{a},\vec{b}):= \prod_{i=1}^{n}u^{+}_{0}(a_i) \prod_{j=1}^{m}u^{-}_0(b_j)$.
\end{thm}

\begin{remark}\label{R:3.6}
\begin{description}\rm
\item{(i)}   The equation expresses $\gamma_t^{(n,m)}$ as an integral in time involving $\gamma^{(n,m+1)}$ and $\gamma^{(n+1,m)}$, thus forming a coupled chain of equations. In statistical physics, it is sometimes called the BBGKY hierarchy\footnote{BBGKY stand for N. N. Bogoliubov, Max Born, H. S. Green, J. G. Kirkwood, and J. Yvon, who derived this type of  hierarchy of equations in the 1930s and 1940s in a series papers.}. It describes the evolution of the limiting $(n,m)$-particle correlation functions and hence the dynamics of the particles.

\item{(ii)}   By Proposition \ref{E:localtime_def}, (\ref{E:hierarchy_1}) is equivalent to
\begin{equation}\label{E:hierarchy_3}
\gamma_t^{(n,m)}(\vec{r},\vec{s}) = \E^{(\vec{r},\vec{s})}\left[\Phi^{(n,m)}(X_{(n,m)}(t)) - \lambda\int_0^t (\Upsilon\gamma_s)^{(n,m)}(X_{(n,m)}(t-s))\, dL^{(n,m)}_s\,\right] .
\end{equation}
Here $L^{(n,m)}$ is boundary local time of $X_{(n,m)}$, the symmetric reflected diffusion on $D_{+}^n\times D_{-}^m$ corresponding to $(I_{(n+m)d\times (n+m)d},\,\rho_{(n,m)})$, and $(\Upsilon v)^{(n,m)}$ is a function on $\partial(D_{+}^n\times D_{-}^m)$ defined as
\begin{equation*}
\quad (\Upsilon v)^{(n,m)}(\vec{r},\vec{s}) :=
\begin{cases}
    v^{(n,m+1)}(\vec{r},(\vec{s},r_i))\frac{\rho_{(n,m)}(\vec{r},\vec{s})}{\rho_+(r_i)} ,   &\text{ if  } (\vec{r},\vec{s})\in \partial_+^i\,  ;\\
    v^{(n+1,m)}((\vec{r},s_j),\vec{s}))\frac{\rho_{(n,m)}(\vec{r},\vec{s})}{\rho_-(s_j)},    &\text{ if  } (\vec{r},\vec{s})\in \partial_-^j\, ;\\
    0, &\text{ otherwise  }.
\end{cases}
\end{equation*}
Observe that the coordinate processes of $X_{(n,m)}$ consist of $n$ independent copies of
reflected diffusions in $\bar{D}_+$ and $m$ independent copies of reflected diffusions
in $\bar{D}_-$.

\item{(iii)}   It is easy to check by using (ii) and Proposition \ref{prop:MildSol_RobinPDE}
that
$$
{\widetilde \gamma}_t^{(n,m)}(\vec{r},\vec{s}) := \prod_{i=1}^{n}u_{+}(t,r_i) \prod_{j=1}^{m}u_{-}(t,s_j)$$ is a solution of  (\ref{E:hierarchy_1}), where $(u_+,u_-)$ is the weak solution of the coupled PDEs \eqref{E:coupledpde:+}-\eqref{E:coupledpde:-} with initial value
$(u^+_0, u^-_0)$.
\end{description}
\end{remark}

\medskip

\noindent{\it Proof of Theorem \ref{T:hierarchy}}.
Recall that $\Xi:= \{\xi:\,\xi_{\pm}(z_{\pm})\leq 1 \hbox{  for every } z\in I^{\eps}\}$.
We can rewrite \eqref{e:3.7} as
   \begin{eqnarray}\label{E:correlation}
      \gamma^{\eps}((\vec{r},\vec{s}),t)
    &=& \sum_{(\vec{a},\vec{b})}\, \gamma^{\eps}((\vec{a},\vec{b}),0)\, p^{\eps}(t,(\vec{r},\vec{s}),(\vec{a},\vec{b}))\,m(\vec{a},\vec{b}) \notag\\
    &&   +  \int_0^t \sum_{(\vec{a},\vec{b})\notin \Xi}\,p^{\eps}(t-s,(\vec{r},\vec{s}), (\vec{a},\vec{b}))\,\E[\mathfrak{K}^{\eps}A((\vec{a},\vec{b}),\eta_s)]\,\eps^{d(n+m)}\,ds \notag\\
    &&  +  \int_0^t \sum_{(\vec{a},\vec{b})\in \Xi}\,p^{\eps}(t-s,(\vec{r},\vec{s}), (\vec{a},\vec{b}))\,\E[\mathfrak{K}^{\eps}A((\vec{a},\vec{b}),\eta_s)]\,\eps^{d(n+m)}\,ds.
    \end{eqnarray}
Fix any $(n,m)\in \mathbb{N}\times \mathbb{N}$, $t>0$ and $(\vec{r},\vec{s})\in (D^{\eps}_{+})^n \times (D^{\eps}_{-})^m$. By a simple counting argument and condition
\eqref{e:3.24} for $\eta^\eps_0$, we see that the first term in (\ref{E:correlation}) equals
\begin{eqnarray*}
&& \E^{\eta^\eps_0}\left[\sum_{(\vec{a},\vec{b})\in (D^{\eps}_{+})^n \times (D^{\eps}_{-})^m}p^{\eps}(t,(\vec{r},\vec{s}),(\vec{a},\vec{b}))\,\prod_{i=1}^n\eta^+(a_i)\,\prod_{j=1}^m\eta^+(b_j)\right]+ o(N)\\
&=& \E^{\eta^\eps_0}\left[\prod_{i=1}^n\<\X^{N,+}_0,\,p^{\eps}(t,r_i,\cdot)\>\,\prod_{j=1}^m\<\X^{N,-}_0,\,p^{\eps}(t,s_j,\cdot)\>\right]+o(N),
\end{eqnarray*}
which converges to $\E^{(\vec{r},\vec{s})}[\Phi^{(n,m)}(X_{(n,m)}(t))]$ by Theorem \ref{T:LCLT_CTRW} and assumption (i) for the initial distributions.  Here $\P^{(\vec{r},\vec{s})}$ is the probability measure for $X_{(n,m)}$ starting at $(\vec{r},\vec{s})$.

We now prove that the second term in (\ref{E:correlation}) tends to $0$ as $\eps\to 0$. The integrand with respect to  $ds$ is at most
\begin{equation}\label{E:hierarchy_2ndterm}
\E\left[\frac{\lambda}{\eps}\sum_{z\in I}\eta^{+}_{\theta}(z_+)\eta^{-}_{\theta}(z_-)\,
\sum_{(\vec{a},\vec{b})\notin \Xi}\,p^{\eps}(t-\theta,(\vec{r},\vec{s}),(\vec{a},\vec{b}))\,2A((\vec{a},\vec{b}),\eta_{\theta})\,\eps^{d(n+m)}\right].
\end{equation}
Note that $\{(\vec{a},\vec{b})\notin \Xi\}$ is a subset of
\begin{equation}\label{E:hierarchy_2ndterm_1}
\bigcup_{w\in I}\,\bigg[\,\left(\bigcup_{k=2}^n\{(\vec{a},\vec{b}):\,\vec{a}(w_+)=k\}\right) \bigcup \left(\bigcup_{\ell=2}^m\{(\vec{a},\vec{b}):\,\vec{b}(w_-)=\ell\} \right)\,\bigg],
\end{equation}
and that for fixed $w\in I$ and $k\in \{2,\cdots,n\}$, we further have
\begin{equation*}
\{(\vec{a},\vec{b}):\,\vec{a}(w_+)=k\} = \bigcup_{\substack{i_1,\cdots,i_k\\\text{distinct}}}\{(\vec{a},\vec{b}):\,a_{i_1}=\cdots=a_{i_k}=w_+\}.
\end{equation*}
Now we restrict the sum over $\{(\vec{a},\vec{b})\notin \Xi\}$ in (\ref{E:hierarchy_2ndterm}) to the subset $\{(\vec{a},\vec{b}):\,a_{i_1}=\cdots=a_{i_k}=w_+\}$, where $w\in I$, $k\in \{2,\cdots,n\}$ and $(i_1,\cdots,i_k)$ are fixed. Moreover, we denote  $(a_1,\cdots,a_k)$
by  $\vec{a_k}$ and  $(a_{k+1},\cdots,a_n)$ by $\vec{a}\setminus \vec{a_k}$. Then
\begin{eqnarray*}\label{E:hierarchy_2ndterm_2}
&& \E\left[\frac{\lambda}{\eps}\sum_{z\in I}\eta^{+}_{\theta}(z_+)\eta^{-}_{\theta}(z_-)\,
\sum_{\{(\vec{a},\vec{b}):\,a_{i_1}=\cdots=a_{i_k}=w_+\}}\,p^{\eps}(t-\theta,(\vec{r},\vec{s}),(\vec{a},\vec{b}))\,2A((\vec{a},\vec{b}),\eta_{\theta})\,\eps^{d(n+m)}\right]\\
&\leq & p^{\eps}(t-\theta,(r_1,\cdots,r_k),(w_+,\cdots,w_+))\,\eps^{kd}\,\sum_{(\vec{a}\setminus \vec{a_k},\,\vec{b})}
p^{\eps}(t-\theta,\vec{r}\setminus \vec{r_k},\vec{a}\setminus \vec{a_k})p(t-\theta,\vec{s},\vec{b})\eps^{d(n+m-k)}\\
&& \quad \cdot\frac{\lambda}{\eps}\,\E\left[\sum_{z\in I}\eta^{+}_{\theta}(z_+)\eta^{-}_{\theta}(z_-)\,2A((\vec{a},\vec{b}),\eta_{\theta})\right]\\
&\leq& \frac{C\,\eps^{kd}}{(t-\theta)^{kd/2}}\,\frac{\lambda}{\eps}\,\#|I^{\eps}|\,\sup_{(\vec{a},\vec{b})}\E[\eta^{+}_{\theta}(z_+)\eta^{-}_{\theta}(z_-)\,2A((\vec{a},\vec{b}),\eta_{\theta})]\\
&\leq& \frac{\lambda\,\eps^{(k-1)d}}{(t-\theta)^{kd/2}}\,C \quad\text{ where }C=C(n,m,\theta,d,D_{\pm})\\
&\leq& \frac{\lambda\,\eps^{d}}{(t-\theta)^{kd/2}}\,C = O(\eps^d) \quad \text{ since }k\geq 2.
\end{eqnarray*}
The second to the last inequality above follows from the bound $\#|{I^{\eps}}|\leq C\,\eps^{-(d-1)}$ (see Lemma \ref{L:DiscreteApprox_SurfaceMea}) and the uniform upper bound (\ref{E:uniformbound}). Repeat the above argument for the other subsets of $\{(\vec{a},\vec{b})\notin \Xi\}$ and use the fact $\#|I^{\eps}|\leq C\,\eps^{-(d-1)}$ again (for $w\in I$ in (\ref{E:hierarchy_2ndterm_1})), we have, for any $\theta\in(0,t)$, (\ref{E:hierarchy_2ndterm}) is of order $\eps$ and hence converges to 0 uniformly for $(\vec{r},\vec{s})$, as $\eps\to 0$. The second term in (\ref{E:correlation}) then converges to 0, by (\ref{E:equicts_G_3}) and LDCT.

For the third term in (\ref{E:correlation}), we split the integrand with respect to  $d\theta$  into three terms corresponding to (\ref{E:KA(xi,eta)_1}), (\ref{E:KA(xi,eta)_2}) and (\ref{E:KA(xi,eta)_3}) respectively. The term corresponding to (\ref{E:KA(xi,eta)_1}) equals
\begin{eqnarray*}
&&  -\dfrac{\lambda}{\eps}\,\sum_{(\vec{a},\vec{b})\in \Xi}\,p^{\eps}(t-s,(\vec{r},\vec{s}), (\vec{a},\vec{b}))\,
\sum_{\substack {z\in I^{\eps} \\\vec{a}(z_+)=1}} \,\Psi(z)\,A((\vec{a},\,(\vec{b},z_-)),\,\eta_{\theta}) \,\eps^{d(n+m)}\\
&=& -\dfrac{\lambda}{\eps}\,\sum_{z\in I^{\eps}}\sum_{\substack{(\vec{a},\vec{b})\in \Xi\\\vec{a}(z_+)=1}}p^{\eps}(t-s,(\vec{r},\vec{s}),(\vec{a},\vec{b}))\,\Psi(z)\,A((\vec{a},\,(\vec{b},z_-)),\,\eta_{\theta}) \,\eps^{d(n+m)}\\
&=& -\dfrac{\lambda}{\eps}\,\sum_{z\in I^{\eps}}\Psi(z)\,\sum_{i=1}^n\sum_{\substack{(\vec{a},\vec{b})\in \Xi\\a_i=z_+}}p^{\eps}(t-s,(\vec{r},\vec{s}),(\vec{a},\vec{b}))\,A((\vec{a},\,(\vec{b},z_-)),\,\eta_{\theta}) \,\eps^{d(n+m)}\\
&=& -\dfrac{\lambda}{\eps}\,\sum_{z\in I^{\eps}}\Psi(z)\,\sum_{i=1}^n\,p^{\eps}(t-s,r_i,z_+)\,\sum_{(\vec{a}\setminus a_i,\,\vec{b})\in\Xi}p^{\eps}(t-s,(\vec{r}\setminus r_i,\vec{s}),(\vec{a}\setminus a_i,\vec{b}))\\
&& \qquad \qquad \times \frac{m((\vec{a},\,(\vec{b},z_-))}{\eps^{d(n+m+1)}}\,\gamma^{\eps}((\vec{a},\,(\vec{b},z_-)),\,\theta) \,\eps^{d(n+m)}\\
&=& -\lambda\sum_{i=1}^n\,\sum_{(\vec{a}\setminus a_i,\,\vec{b})\in\Xi}p^{\eps}(t-s,(\vec{r}\setminus r_i,\vec{s}),(\vec{a}\setminus a_i,\vec{b}))\,m(\vec{a}\setminus a_i,\vec{b})\\
&&\qquad \times \sum_{z\in I^{\eps}}\sigma_{\eps}(z)\,p^{\eps}(t-s,r_i,z_+)\,\gamma^{\eps}((\vec{a},\,(\vec{b},z_-)),\,\theta) .
\end{eqnarray*}
By Theorem \ref{T:LCLT_CTRW} and Lemma \ref{L:DiscreteApprox_SurfaceMea},
$$
\lim_{\eps \to 0}  \sum_{z\in I^{\eps}}\sigma_{\eps}(z)\,p^{\eps}(t-s,r_i,z_+)\,\gamma^{\eps}((\vec{a},\,(\vec{b},z_-)),\,\theta)
=  \int_{I}p (t-s,r_i,z)\,\gamma_{\theta}(\vec{a},(\vec{b},z))\,d\sigma(z)
$$
and the convergence is uniform for $r_i\in D^{\eps}_{+}$. Therefore, by applying Theorem \ref{T:LCLT_CTRW} again, the term corresponding to (\ref{E:KA(xi,eta)_1}) converges to
\begin{equation*}
-\lambda \sum_{i=1}^n \int_{\partial_+^i} \gamma_{\theta}(\vec{a},(\vec{b},a_i)) \,p(t-\theta,(\vec{r},\vec{s}),(\vec{a},\vec{b}))\,\frac{\rho_{(n,m)}(\vec{r},\vec{s})}{\rho_+(r_i)}\,d\vec{b}\,da_1\cdots d\sigma(a_i)\cdots da_n .
\end{equation*}

We repeat the same argument for the term corresponding to (\ref{E:KA(xi,eta)_2}). Moreover, note that the term corresponding to (\ref{E:KA(xi,eta)_3}) will not contribute to the limit as $\eps\to 0$, by the same argument we used for the second term in (\ref{E:correlation}). Therefore, the integrand of the second term in (\ref{E:correlation}) converges to
\begin{eqnarray*}
&& -\lambda \sum_{i=1}^n \int_{\partial_+^i} \gamma_{\theta}(\vec{a},(\vec{b},a_i)) \,p(t-\theta,(\vec{r},\vec{s}),(\vec{a},\vec{b}))\,\frac{\rho_{(n,m)}(\vec{a},\vec{b})}{\rho_+(a_i)}\,
d\vec{b}\,da_1\cdots d\sigma(a_i)\cdots da_n \\
&& -\lambda \sum_{j=1}^m \int_{\partial_-^j} \gamma_{\theta}((\vec{a},b_j),\vec{b})) \,p(t-\theta,(\vec{r},\vec{s}),(\vec{a},\vec{b}))\,\frac{\rho_{(n,m)}(\vec{a},\vec{b})}{\rho_-(b_j)}\,
 d\vec{a}\,db_1\cdots d\sigma(b_j)\cdots db_m .
\end{eqnarray*}
The integral for $\theta\in(0,t)$ in the third term in (\ref{E:correlation}) then converges to the desired quantity, by (\ref{E:equicts_G_3}) and LDCT. The proof is complete.
\qed

\medskip

In view of Remark \ref{R:3.6}(iii), the proof of Theorem \ref{T:correlation} (Propagation of Chaos) will be complete once we establish the uniqueness of the solution of the limiting hierarchy (\ref{E:hierarchy_1}). This will be accomplished in Theorem \ref{T:Uniqueness_hierarchy} in the next subsection.

\subsection{Uniqueness of infinite hierarchy}

Uniqueness of BBGKY hierarchy is an important issue in statistical physics. For example, it is a key step in the derivation of the cubic non-linear Schr\"odinger equation from the quantum dynamics of many body systems obtained in \cite{lEbShtY07}. Our BBGKY hierarchy (\ref{E:hierarchy_1}) is new to the literature and the proof of its uniqueness involves a representation and manipulations of the hierarchy in terms of trees. The technique is related but different from that in \cite{lEbShtY07},  which used the Feynman diagrams.

Note that, by Theorem \ref{T:hierarchy}, $\gamma^{(n,m)}_t$ can be extended  continuously to   $t=0$. Uniqueness of solution for the hierarchy will be established on a subset of the space
$$C([0,T], \mathcal{D}):= \bigoplus_{(n,m)\in \mathbb{N}\times \mathbb{N}}C([0,T],\,\bar{D}_+^n\times \bar{D}_-^m)$$
equipped with the product topologies induced by the uniform norm $\|\,\cdot\,\|_{(T,n,m)}$ on $[0,T]\times \bar{D}_+^n\times \bar{D}_-^m$.

\begin{thm}\label{T:Uniqueness_hierarchy}(Uniqueness of the infinite hierarchy)
Given any $T>0$. Suppose $\beta_t=\{\beta_t^{(n,m)}\}\in C([0,T], \mathcal{D})$ is a solution to the infinite hierarchy (\ref{E:hierarchy_1}) with zero initial condition (i.e. $\beta_0=\Phi=0$) and satisfies $\|\beta_t^{(n,m)}\|_{(T,n,m)}\leq C^{n+m}$ for some $C\geq 0$. Then we have $\|\beta_t^{(n,m)}\|_{(T,n,m)}=0$ for every $n, m\in \mathbb{N}$.
\end{thm}

The remaining of this subsection is devoted to give a proof of this theorem.

\textbf{Convention in this subsection:}
$\beta=\{\beta^{(n,m)}\}$ will always denote the functions stated in Theorem \ref{T:Uniqueness_hierarchy}.   For notational simplicity,  we will also assume $\lambda=2$ and $\rho_{\pm}=1$. The proof for the general case is the same.
we will also drop $T$ from the notation $\|\beta_t^{(n,m)}\|_{(T,n,m)}$.

It is convenient to rewrite the infinite hierarchy (\ref{E:hierarchy_1}) in a more compact form as
\begin{equation}\label{E:hierarchy_2}
\gamma_t^{(n,m)} = P^{(n,m)}_t\Phi^{(n,m)}-\int_0^t
P^{(n,m)}_{t-s}\left(\sum_{i=1}^nV^{+}_{i}\gamma_s^{(n,m+1)}+ \sum_{j=1}^mV^{-}_{j}\gamma_s^{(n+1,m)}\right)\,ds \; ,
\end{equation}
where $V^{+}_{i}\gamma^{(n,m+1)}$ is a measure concentrated on $\partial_+^i$ defined as
\begin{eqnarray*}
V^{+}_{i}\gamma^{(n,m+1)} &:= &\gamma^{(n,m+1)}(\vec{a},(\vec{b},a_i))\,d\sigma_{(n,m)}\Big|_{\partial_+^i}(\vec{a},\vec{b})\\
&=& \gamma^{(n,m+1)}(\vec{a},(\vec{b},a_i))\,d\sigma\Big|_{I}(a_i)\,d(\vec{a}\setminus a_i)\,d\vec{b} .
\end{eqnarray*}
Here $\sigma_{(n,m)}\Big|_{\partial_+^i}$ is the surface measure of $\partial(D_{+}^n\times D_{-}^m)$ restricted to $\partial_+^i$. Similarly, $V^{-}_{j}\gamma^{(n+1,m)}$ is a measure concentrating on $\partial_-^j$ defined as
\begin{eqnarray*}
V^{-}_{j}\gamma^{(n+1,m)} &:= &\gamma^{(n+1,m)}((\vec{a},b_j),\vec{b})\,d\sigma_{(n,m)}\Big|_{\partial_-^j}(\vec{a},\vec{b})\\
&=& \gamma^{(n+1,m)}((\vec{a},b_j),\vec{b})\,d\sigma\Big|_{I}(b_j)\,d\vec{a}\,d(\vec{b}\setminus b_j) .
\end{eqnarray*}

\subsubsection{Duhamel tree expansion}

\begin{figure}[h]
	\begin{center}
    \vspace{-0.5em}
	\includegraphics[scale=.4]{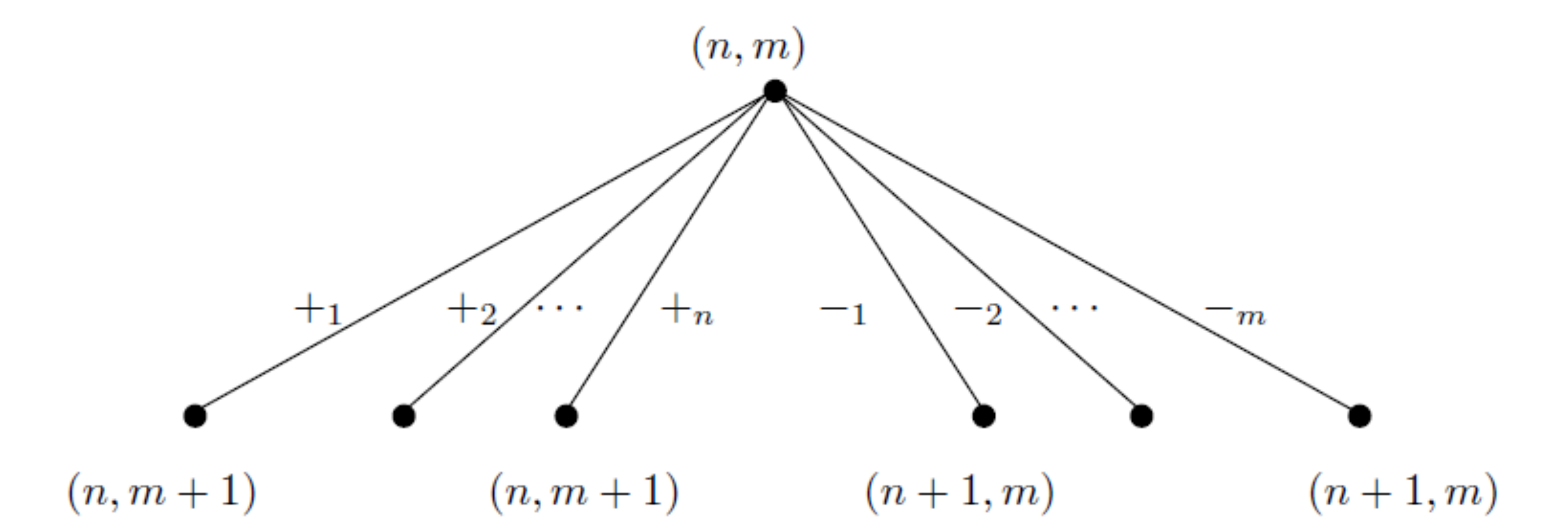}
    \vspace{-0.5em}
    \caption{$\mathbb{T}^{(n,m)}_1$}\label{fig:TreeT_nm1}
    \vspace{-0.5em}
	\end{center}
	\end{figure}

We now  describe the infinite hierarchy in  detail. It is natural and illustrative to represent the infinite hierarchy in terms of a tree structure, with the `root' at the top and the `leaves' at the bottom. Fix two positive integers $n$ and $m$. We construct a sequence of finite trees $\{\mathbb{T}^{(n,m)}_N:\;N=0,1,2,\cdots \}$ recursively as follows.
\begin{enumerate}
\item   $\mathbb{T}^{(n,m)}_0$ is the root, with label $(n,m)$.
\item   $\mathbb{T}^{(n,m)}_1$ is constructed from $\mathbb{T}^{(n,m)}_0$ by attaching $n+m$ new vertices (call them leaves of $\mathbb{T}^{(n,m)}_1$) to it. More precisely, we attach $n+m$ new vertices to the root by drawing $n$ $`+'$ edges and $m$ $`-'$ edges from the root. Those new leaves drawn by the $`+'$ edges are labeled $(n,m+1)$, while those drawn by the $`-'$ edges are labeled $(n+1,m)$. We also label the edges as $\{+_i\}_{i=1}^n$ and $\{-_j\}_{j=1}^m$ (See Figure \ref{fig:TreeT_nm1}).
\item   When $N=2$, we view each of the $n+m$ leaves of $\mathbb{T}^{(n,m)}_1$ as a `root' (with a new label, being either $(n,m+1)$ or $(n+1,m)$), and then attach new leaves (leaves of $\mathbb{T}^{(n,m)}_2$) to it by drawing $`\pm'$ edges. Hence there are $(m+n)(m+n+1)$ new leaves, coming from $n^2+m(n+1)$ new $`+'$ edges and $n(m+1)+m^2$ new $`-'$ edges.
\item   Having drawn $\mathbb{T}^{(n,m)}_{N-1}$, we construct $\mathbb{T}^{(n,m)}_N$ by attaching new edges and new leaves from each leaf of $\mathbb{T}^{(n,m)}_{N-1}$ by the same construction, viewing a leaf of $\mathbb{T}^{(n,m)}_{N-1}$ as a `root'.
\end{enumerate}

In $\mathbb{T}^{(n,m)}_N$, the root is connected to each leaf by a unique path $\vec{\theta}=(\theta_1,\,\theta_2,\,\cdots,\,\theta_N)$ formed by the $`\pm'$ edges. Moreover, such a path passes through a sequence of labels formed by the leaves of $\mathbb{T}^{(n,m)}_k$ ($k=1,2,\cdots ,N$). We denote these labels by $\vec{l}(\vec{\theta})=(l_1(\vec{\theta}),\,l_2(\vec{\theta}),\cdots,\,l_N(\vec{\theta}))$. For example, when $(n,m)=(2,5)$, $N=3$ and the path is $\vec{\theta}=(+_1,\,-_6,\,-_5)$. Then $\vec{l}(\vec{\theta})=((2,6),\,(3,6),\,(4,6))$ and the path connects the root to a leaf of $\mathbb{T}^{(2,5)}_3$ with label $(4,6)$.
Note that the label is not one-to-one. For example,  $\vec{l}  (+_1,\,-_6,\,-_5)
=\vec{l} (+_2,\,-_6,\,-_4)$.

For mnemonic reason, we use the same notation $\mathbb{T}^{(n,m)}_N$ to denote the collection of paths that connect the root to a leaf in $\mathbb{T}^{(n,m)}_N$. By induction, the total number of paths (or the total number of leaves) is
\begin{equation}\label{E:NumberOfLeaves}
(n+m)(n+m+1)\cdots (n+m+N-1)=\frac{(n+m+N-1)!}{(n+m-1)!}.
\end{equation}

Iterating (\ref{E:hierarchy_2}) $N$ times gives
\begin{eqnarray}\label{E:hierarchy_4}
 \beta_t^{(n,m)}
&=& -\int_{t_2=0}^t
P^{(n,m)}_{t-t_2}\left(\sum_{i=1}^nV^{+}_{i}\beta_{t_2}^{(n,m+1)}+ \sum_{j=1}^mV^{-}_{j}\beta_{t_2}^{(n+1,m)}\right)  dt_2\notag\\
&=& \cdots \notag \\
&=& (-1)^N \int_{t_2=0}^t\int_{t_3=0}^{t_2}\cdots \int_{t_{N+1}=0}^{t_N}
dt_2 \cdots dt_{N+1} \notag \\
&&\qquad \qquad \sum_{\vec{\theta}\in \,\mathbb{T}^{(n,m)}_N}\,
P^{(n,m)}_{t-t_2}\,V_{\theta_1}\,P^{l_1(\vec{\theta})}_{t_2-t_3}\,V_{\theta_2}\,P^{l_2(\vec{\theta})}_{t_3-t_4}\,V_{\theta_3}\cdots P^{l_{N-1}(\vec{\theta})}_{t_N-t_{N+1}}\,V_{\theta_N}\,\beta^{l_N(\vec{\theta})}_{t_{N+1}},
\end{eqnarray}
where $V_{\theta_i}$ (for $i=1,\,2,\cdots,\,N$) is defined by $V_{+_i}=V^+_i$ and $V_{-_j}=V^-_j$. For example, when $(n,m)=(2,5)$, $N=3$ and the path is $\vec{\theta}=(+_1,\,-_6,\,-_5)$. Then
\begin{equation*}
P^{(n,m)}V_{\theta_1}P^{l_1(\vec{\theta})}V_{\theta_2}P^{l_2(\vec{\theta})}V_{\theta_3}\beta^{l_3(\vec{\theta})} =
P^{(2,5)}V^+_{1}P^{(2,6)}V^-_{6}P^{(3,6)}V^-_{5}\beta^{(4,6)}.
\end{equation*}

\subsubsection{Telescoping via Chapman-Kolmogorov equation}
By (\ref{E:NumberOfLeaves}), the right hand side of (\ref{E:hierarchy_4}) is a sum of $(n+m)(n+m+1)\cdots (n+m+N-1)$ terms of multiple integrals. We will apply the bound
$\|\beta_t^{(p,q)}\|_{(p, q)}\leq C^{p+q}$
 to each term, and then simplify the integrand using
 Chapman-Kolmogorov  equation.

We demonstrate this for the twelve terms for the case $(n,m,N)=(1,2,2)$. The twelve terms on the right hand side of (\ref{E:hierarchy_4}) are
\begin{eqnarray} \label{E:Tree_122}
P^{(1,2)}_{t-t_2}\Big\{&& V^{+}_{1}P^{(1,3)}_{t_2-t_3}\left(V^{+}_{1}\,\beta^{(1,4)}_{t_3}+(V^{-}_{1}+V^{-}_{2}+V^{-}_{3})\,\beta^{(2,3)}_{t_3}\right) \notag\\
&+& V^{-}_{1}P^{(2,2)}_{t_2-t_3}\left((V^{+}_{1}+V^{+}_{2})\,\beta^{(2,3)}_{t_3}+(V^{-}_{1}+V^{-}_{2})\,\beta^{(3,2)}_{t_3}\right) \notag\\
&+& V^{-}_{2}P^{(2,2)}_{t_2-t_3}\left((V^{+}_{1}+V^{+}_{2})\,\beta^{(2,3)}_{t_3}+(V^{-}_{1}+V^{-}_{2})\,\beta^{(3,2)}_{t_3}\right)\,
\Big\}.
\end{eqnarray}
The first four terms came from the leftmost leaf of the previous level, we group them together to obtain, for $(x,y_1,y_2)\in \bar{D}_+\times \bar{D}_-^2$,
\begin{eqnarray*}
&& \Big|P^{(1,2)}_{t-t_2}V^{+}_{1}P^{(1,3)}_{t_2-t_3}\left(V^{+}_{1}\,\beta^{(1,4)}_{t_3}+(V^{-}_{1}+V^{-}_{2}+V^{-}_{3})\,\beta^{(2,3)}_{t_3}\right)
(x,\,y_1,\,y_2)\Big|\\
&\leq &
C^5 \,\int d\sigma(x')\,dy_1'\,dy_2'\;p^{(1,2)}(t-t_2,\,(x,y_1,y_2),(x',y_1',y_2'))\;\\
&& \qquad \bigg(\int d\sigma(a)\,db_1\,db_2\,db_3 + \int da\,d\sigma(b_1)\,db_2\,db_3 +\int da\,db_1\,d\sigma(b_2)\,db_3 + \int da\,db_1\,db_2\,d\sigma(b_3)
\bigg)\\
&&   \qquad \qquad   p^{(1,3)}(t_2-t_3,\,(x',y_1',y_2',x'),(a,b_1,b_2,b_3))\\
&=&
C^5\int d\sigma(x')\,p^{+}(t-t_2,x,x')\bigg(
\int d\sigma(a)\,p^{+}(t_2-t_3,x',a)+\int d\sigma(b_1)\,p^{-}(t-t_3,y_1,b_1)\\
&&\qquad \qquad \qquad +\int d\sigma(b_2)\,p^{-}(t-t_3,y_2,b_2)\,+\int d\sigma(b_3)\,p^{-}(t_2-t_3,x',b_3)\,
\bigg).
\end{eqnarray*}
Note the telescoping effect upon using the Chapman-Kolmogorov equation for the middle two terms  in the last equality above  gives rise to $t-t_3$ rather than $t_2-t_3$.

\begin{figure}[h]
	\begin{center}
    \vspace{-0.5em}
	\includegraphics[scale=.4]{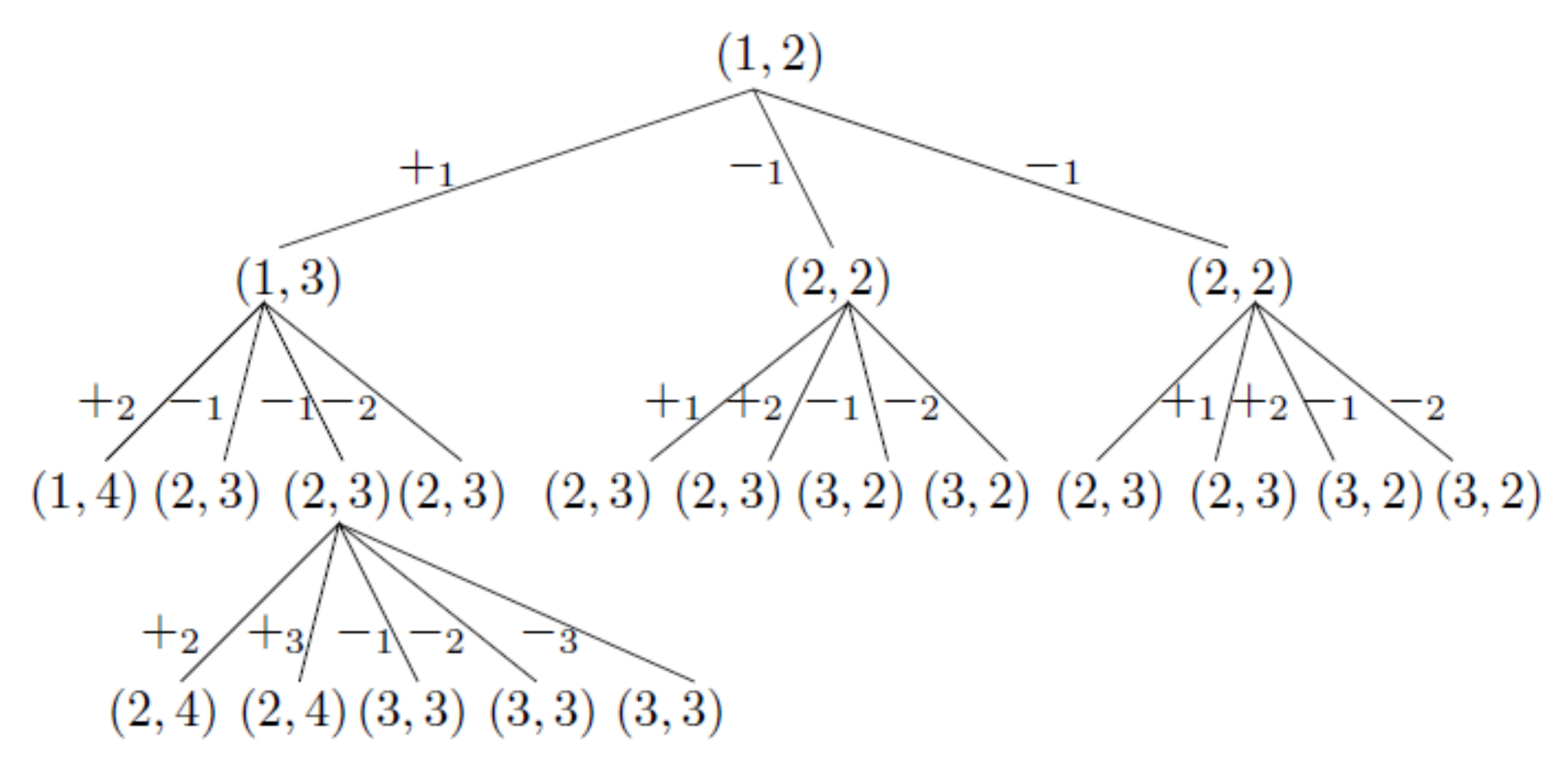}
    \vspace{-0.5em}
    \caption{$\mathbb{S}^{(1,2)}_2$ together with 5 leaves in $\mathbb{S}^{(1,2)}_3$}\label{fig:TreeS_12_2}
    \vspace{-0.5em}
	\end{center}
	\end{figure}

We  apply (\ref{E:Surface_integral_boundedness}) to obtain
\begin{align*}
&\Big\|P^{(1,2)}_{t-t_2}V^{+}_{1}P^{(1,3)}_{t_2-t_3}\left(V^{+}_{1}\,\beta^{(1,4)}_{t_3}+(V^{-}_{1}+V^{-}_{2}+V^{-}_{3})\,\beta^{(2,3)}_{t_3}\right)
\Big\|_{(1,2)} \\
&\leq
C^{5}\,\frac{C_+}{\sqrt{t-t_2}}\,\left(\frac{C_+}{\sqrt{t_2-t_3}}+ \frac{2\,C_-}{\sqrt{t-t_3}}+ \frac{C_-}{\sqrt{t_2-t_3}}\right),
\end{align*}
where $C_{\pm}=C(D_{\pm},T)$ are positive constants. Repeat the above argument for the remaining eight terms of (\ref{E:Tree_122}), we obtain
\begin{align}\label{E:Telescoping_122}
\|\beta_t^{(1,2)}\|_{(1,2)}
&\leq
C^{5}\,\int_{t_2=0}^{t}\int_{t_3=0}^{t_2}\frac{C_+}{\sqrt{t-t_2}}\,\left(\frac{C_+}{\sqrt{t_2-t_3}}+ \frac{2\,C_-}{\sqrt{t-t_3}}+ \frac{C_-}{\sqrt{t_2-t_3}}\right) \notag\\
& \quad +\frac{2\,C_-}{\sqrt{t-t_2}}\,\left(\frac{C_+}{\sqrt{t-t_3}}+ \frac{C_+}{\sqrt{t_2-t_3}}+\frac{C_-}{\sqrt{t-t_3}}+ \frac{C_-}{\sqrt{t_2-t_3}}\right).
\end{align}
The key is to visualize the twelve terms on the right as 12 paths of $\mathbb{T}^{(1,2)}_2$ \textit{with the edges relabeled}. We denote this relabeled tree by $\mathbb{S}^{(1,2)}_2$ (See Figure \ref{fig:TreeS_12_2}, ignoring the five leaves in $\mathbb{S}^{(1,2)}_3$ at the moment). More precisely, since all twelve terms on the right are of the form $\frac{C_{\pm}}{\sqrt{t_{p}-t_2}}\,\frac{C_{\pm}}{\sqrt{t_{q}-t_3}}$, we only need to record the indexes $(p,q)$ and the $\pm$ sign. For example, the first four terms can be represented by $$(+_1\,+_2,\;+_1\,-_1,\;+_1\,-_1,\;+_1\,-_2).$$
Each $+_1\,-_1$ corresponds to $\frac{C_+}{\sqrt{t-t_2}}\,\frac{C_-}{\sqrt{t-t_3}}$ and hence it appears twice. In $\mathbb{S}^{(1,2)}_2$, these four paths are formed by a $+_1$ edge followed by four edges with labels $\{+_2,\,-_1,\,-_1,\,-_2\}$.

\begin{figure}[h]
	\begin{center}
    \vspace{-0.5em}
	\includegraphics[scale=.4]{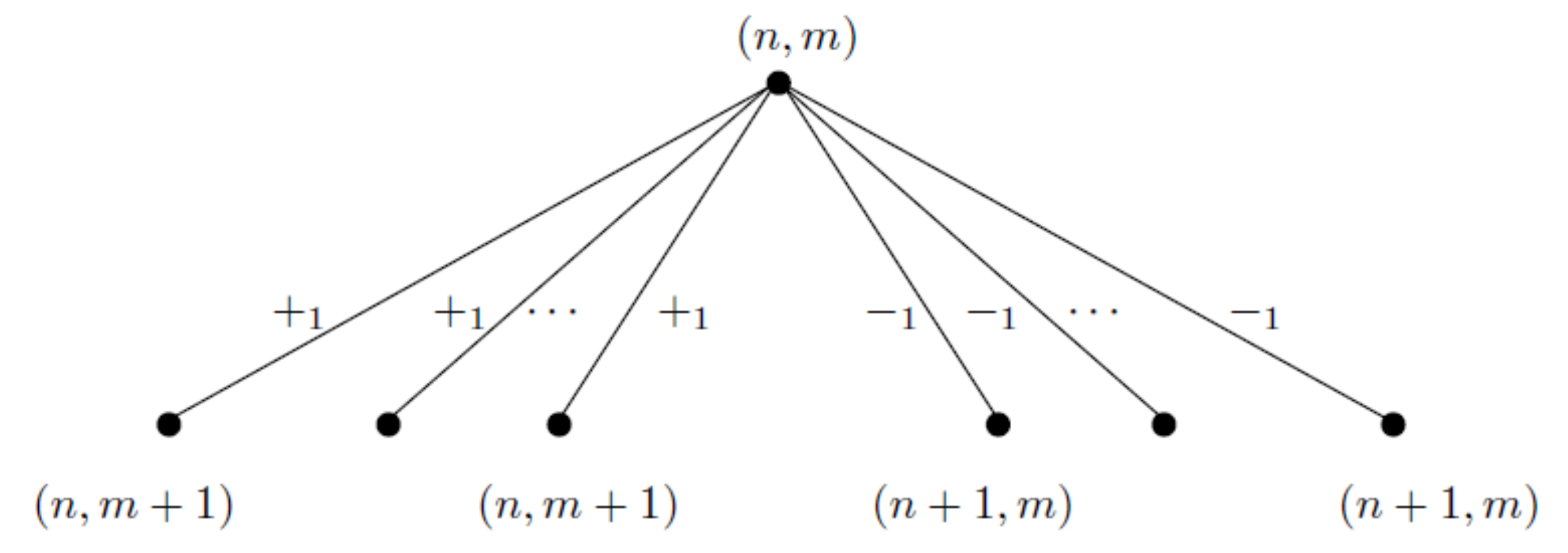}
    \vspace{-0.5em}
    \caption{$\mathbb{S}^{(n,m)}_1$}\label{fig:TreeS_nm1}
    \vspace{-0.5em}
	\end{center}
	\end{figure}

In general, we obtain $\mathbb{S}^{(n,m)}_N$ by relabeling the edges of $\mathbb{T}^{(n,m)}_N$, while keeping the labels for the vertices and the $\pm$ sign for the edges. The relabeling of edges are performed as follows:
\begin{enumerate}
\item   At level $1$, we assign the number `1' to all the edges connected to the root. Hence we have $n$ `$+_1$' edges and $m$ `$-_1$' edges, rather than the labels $\{+_i\}_{i=1}^n$ and $\{-_j\}_{j=1}^m$ (See Figure \ref{fig:TreeS_nm1} for $\mathbb{S}^{(n,m)}_1$).
\item   At level $k\geq 2$, consider the set $\Lambda^+ := \{+_1,\,\cdots,\,+_1,\,+_2,\,+_3,\,\cdots,\,+_k\}$ in which we have $n$ copies of $+_1$ (hence there are $n+k-1$ elements in $\Lambda^+$, in which $n$ of them are $+_1$). Similarly, let $\Lambda^- := \{-_1,\,\cdots,\,-_1,\,-_2,\,-_3,\,\cdots,\,-_k\}$ in which we have $m$ copies of $-_1$. For an arbitrary leaf $\xi$ of $\mathbb{T}^{(n,m)}_{k-1}$, let $R^{\xi}$ be the labels of (the edges of) the path from the root to $\xi$ in $\mathbb{S}^{(n,m)}_{k-1}$, counting with multiplicity. Finally, the collection of new labels of the edges below $\xi$, denoted by $L^{\xi}$, is chosen in such a way that
    $$ \Lambda^+ \cup \Lambda^- = R^{\xi} \cup L^{\xi}\qquad (\text{counting multiplicity}).$$
    Since $|R^{\xi}|=k-1$ and $|L^{\xi}|=n+m+k-1$ (again, counting multiplicity), the cardinalities of the two sides match: $$(n+k-1)+(m+k-1)=(k-1)+(n+m+k-1).$$
    Induction shows that $R^{\xi}\subset \Lambda^+ \cup \Lambda^-$ and the choice for $L^{\xi}$ is unique.
For example, for leaf $\xi =(1, 3)$ of $\mathbb{T}^{(1,2)}_1$,
$R^\xi=\{+_1\}$,  $\Lambda^+ := \{+_1,\,+_2\}$ and
$\Lambda^- := \{-_1, \, ,-_1,\,-_2 \}$.
So $L^\xi=\{ +_2, \, -_1, \, ,-_1,\,-_2 \}$, which gives the new labels to the edges below $\xi$;
see Figure \ref{fig:TreeS_12_2}.

\end{enumerate}

As a further illustration,
we continue to `grow' $\mathbb{S}^{(1,2)}_2$ (see Figure \ref{fig:TreeS_12_2}) by adding suitably labeled edges to leaves of $\mathbb{S}^{(1,2)}_2$. Precisely, let $\xi$ be a leaf of $\mathbb{S}^{(1,2)}_2$.
\begin{itemize}
\item   If $R^{\xi}=\{+_1,\,+_2\}$, then $L^{\xi}=\{+_3,\,-_1,\,-_1,\,-_2,\,-_3\}$ (this is the case for the leftmost leaf, which has  label $(1,4)$)
\item   If $R^{\xi}=\{+_1,\,-_1\}$, then $L^{\xi}=\{+_2,\,+_3,\,-_1,\,-_2,\,-_3\}$ (shown in Figure \ref{fig:TreeS_12_2}).
\item   If $R^{\xi}=\{+_1,\,-_2\}$, then $L^{\xi}=\{+_2,\,+_3,\,-_1,\,-_1,\,-_3\}$.
\item   If $R^{\xi}=\{-_1,\,+_1\}$, then $L^{\xi}=\{+_2,\,+_3,\,-_1,\,-_2,\,-_3\}$.
\item   If $R^{\xi}=\{-_1,\,+_2\}$, then $L^{\xi}=\{+_1,\,+_3,\,-_1,\,-_2,\,-_3\}$.
\item   If $R^{\xi}=\{-_1,\,-_1\}$, then $L^{\xi}=\{+_1,\,+_2,\,+_3,\,-_2,\,-_3\}$.
\item   If $R^{\xi}=\{-_1,\,-_2\}$, then $L^{\xi}=\{+_1,\,+_2,\,+_3,\,-_1,\,-_3\}$.
\end{itemize}

For mnemonic reason, we use the same notation $\mathbb{S}^{(n,m)}_N$ to denote the collection of paths from the root to the leaves of $\mathbb{S}^{(n,m)}_N$. Any such path is represented by the \textit{ordered} (new)  labels of the edges. We now `forget' the sign of the edges and only record the integer labels. For example, the path $(-_1,\,+_1,\,-_2,\,+_3)$ is replaced by $(1,\,1,\,2,\,3)$.

Using the hypothesis
$\|\beta_t^{(n,m)}\|_{(n,m)}\leq C^{n+m}$
(of Theorem  \ref{T:Uniqueness_hierarchy}) and applying Chapman-Kolmogorov equation to (\ref{E:hierarchy_4}), and then applying (\ref{E:Surface_integral_boundedness}), we obtain the following lemma by the same argument that we used to obtain (\ref{E:Telescoping_122}).

\begin{lem}\label{L:beta_I_N}
\begin{equation}
\|\beta_t^{(n,m)}\|_{(n,m)}
\leq C^{n+m+N}\,\,(C_+\vee C_-)^N\,I^{(n,m)}_N(t),
\end{equation}
where
\begin{equation}
I^{(n,m)}_N(t) :=  \int_{t_2=0}^{t}\cdots \int_{t_{N+1}=0}^{t_N}
\sum_{\vec{\upsilon}\,\in \mathbb{S}^{(n,m)}_N}\frac{1}{\sqrt{(t_{\upsilon_1}-t_2)\,(t_{\upsilon_2}-t_3)\cdots(t_{\upsilon_N}-t_{N+1})}} dt_2
\cdots dt_{N+1} .
\end{equation}
\end{lem}

Our goal is to show that $I^{(n,m)}_N(t)\leq (C\,t)^{N/2}$ for some $C=C(n,m)>0$. This will imply $\|\beta_t^{(n,m)}\|_{(n,m)}=0$ for $t>0$ small enough. Clearly we have
\begin{eqnarray*}
I^{(n,m)}_N(t) &\leq& \frac{(n+m+N-1)!}{(n+m-1)!} \int_{t_2=0}^{t}\cdots \int_{t_{N+1}=0}^{t_N}
\frac{1}{\sqrt{(t-t_2)\,(t_2-t_3)\cdots(t_N-t_{N+1})}} dt_2
\cdots dt_{N+1}  \\
&=&   \frac{(n+m+N-1)!}{(n+m-1)!}\,\frac{(\pi\,t)^{N/2}}{\Gamma(\frac{N+2}{2})}.
\end{eqnarray*}
Unfortunately, this crude bound is asymptotically larger than $(Ct)^{N/2}$ for any $C>0$.

\subsubsection{Comparison with a `dominating' tree}

Note that
\begin{eqnarray*}
I^{(1,2)}_3(t) &\leq& \int_{t_2=0}^{t}\int_{t_3=0}^{t_2}\int_{t_4=0}^{t_3}
\left(\frac{3}{\sqrt{t-t_2}}\right)\,\left(\frac{2}{\sqrt{t_2-t_3}}+ \frac{2}{\sqrt{t-t_3}}\right) \notag\\
&& \qquad \left(\frac{2}{\sqrt{t_3-t_4}}+ \frac{2}{\sqrt{t_2-t_4}}+ \frac{1}{\sqrt{t-t_4}}\right)
dt_2 \, dt_3 \, dt_4 .
\end{eqnarray*}
This is obtained by comparing the labels in $\mathbb{S}^{(1,2)}_3$ with a `dominating'
labeling,
in which the labels of the edges below \textit{every} leaf of $\mathbb{S}^{(1,2)}_2$ is $\{3,\,3,\,2,\,2,\,1\}$ (the $\pm$ sign is discarded). This trick enables us to group the terms at each level.

For the general case, let $\xi$ be an arbitrary leaf of $\mathbb{S}^{(n,m)}_{k-1}$. Note that in $L^{\xi}$, each of the integers $2,\,3,\,\cdots,\,k$ appears at most twice and the integer $1$ appears at most $n+m$ times. We compare $L^{\xi}$ with the `dominating' label  $\tilde{L}^\xi$ defined below:

\begin{tabular}{|r|c|l|}\hline
&&\\
Level $k$ & $|\tilde{L}^\xi|=|L^{\xi}|$ &  $\tilde{L}^\xi$ \\
\hline\hline
1 & $n+m$ & $1,\,1,\,1,\cdots,\,1$ \\
2 & $n+m+1$ & $2,\,2,\,1,\,1,\,1,\cdots,\,1$ \\
3 & $n+m+2$ & $3,\,3,\,2,\,2,\,1,\,1,\,1,\cdots,\,1$\\
$\cdots$ & $\cdots$ & $\cdots $\\
$n+m-1$ & $2(n+m)-2$ & $n+m-1,\,n+m-1,\,n+m-2,\,\cdots,\,3,\,3,\,2,\,2,\,1,\,1$\\
$n+m$ & $2(n+m)-1$ & $n+m,\,n+m,\,n+m-1,\,n+m-1,\,\cdots,\,3,\,3,\,2,\,2,\,1$\\
$n+m+1$ & $2(n+m)$ & $n+m+1,\,n+m+1,\,n+m,\,n+m,\,\cdots,\,3,\,3,\,2,\,2$\\
$n+m+2$ & $2(n+m)+1$ & $n+m+2,\,n+m+2,\,n+m+1,\,\cdots,\,3,\,3,\,2$\\
$\cdots$ & $\cdots$ & $\cdots $\\
$N$ & $n+m+N-1$ & $N,\,N,\,N-1,\,N-1,\,\cdots,\,c,\,b,\,a$\\
\hline
\end{tabular}

In the last row, if $n+m+N-1$ is even, then $a=b=(N-n-m+3)/2$ and $c=b+1$; if $n+m+N-1$ is odd, then $a=(N-n-m+2)/2$ and $b=c=a+1$.

We can now group the terms in each level $k$ as a sum of $k$ terms to obtain
\begin{eqnarray*}
I^{(n,m)}_N(t)
&\leq& \int_{t_2=0}^{t}\int_{t_3=0}^{t_2}\cdots \int_{t_{N+1}=0}^{t_N}
\left(\frac{n+m}{\sqrt{t-t_2}}\right)\,
\left(\frac{2}{\sqrt{t_2-t_3}}+\frac{n+m-1}{\sqrt{t-t_3}}\right)\,\\
&& \qquad \left(\frac{2}{\sqrt{t_3-t_4}}+\frac{2}{\sqrt{t_2-t_4}}+\frac{n+m-2}{\sqrt{t-t_4}}\right)\,\cdots\\
&& \qquad \left(\frac{2}{\sqrt{t_{n+m}-t_{n+m+1}}}+\cdots +\frac{2}{\sqrt{t_{2}-t_{n+m+1}}}+\frac{1}{\sqrt{t-t_{n+m+1}}}\right)\\
&& \qquad \prod_{k=n+m+1}^{N}
\left(\frac{2}{\sqrt{t_{k}-t_{k+1}}}+\frac{2}{\sqrt{t_{k-1}-t_{k+1}}}+\cdots +\frac{2}{\sqrt{t_{2}-t_{k+1}}}\right) dt_2 \cdots dt_{N+1}.
\end{eqnarray*}
In the last term, we have used the observation that when $k>n+m$, the smallest element in $\tilde{L}^\xi$ is at least 2 and so the sum stops before reaching $1/\sqrt{t-t_{k+1}}$. From this and the simple estimates like
$$\frac{2}{\sqrt{t_3-t_4}}+\frac{2}{\sqrt{t_2-t_4}}+\frac{n+m-2}{\sqrt{t-t_4}}\leq \frac{2(n+m)}{3}\left(\frac{1}{\sqrt{t_3-t_4}}+\frac{1}{\sqrt{t_2-t_4}}+\frac{1}{\sqrt{t-t_4}}\right),$$
we have derived the following
\begin{lem}\label{L:I_N_J_N}
\begin{equation}
I^{(n,m)}_N(t)  \leq  \frac{(n+m)^{(n+m)}}{(n+m)!}\,2^N\,J_N(t),
\end{equation}
where
\begin{equation}
J_N(t) := \int_0^t \int_0^{t_2}\cdots\int_0^{t_N}\prod_{i=2}^{N+1}\left(\sum_{j=1}^{i-1}\frac{dt_i}{\sqrt{t_j-t_i}}\right) .
\end{equation}
\end{lem}

\subsubsection{Estimating $J_N$}

Our goal in this section is show that $J_N(t)\leq (C\,t)^{N/2}$ for some $C>0$. Our proof relies on the following recursion formula pointed out
to us by David Speyer:
\begin{equation}\label{E:Recursion_J_N(t)}
J_N(t) = \sum_{k=1}^N \frac{1}{k!} \sum_{ \substack{n_2+n_3 + \cdots + n_{k+1} = N \\ n_2, n_3, \cdots, n_{k+1} \geq 1} }\;\prod_{j=2}^{k+1} \int_{0}^t \frac{J_{n_j-1}(t_j)}{\sqrt{t-t_j}}\,dt_j.
\end{equation}
We assume (\ref{E:Recursion_J_N(t)}) for now and use it to establish the following lemma. The proof of (\ref{E:Recursion_J_N(t)}) will be given immediately after it.

\begin{lem}\label{L:Rate_J_N}
$J_N(t)$ is homogeneous in the sense that
\begin{equation}
J_N(t) = J_N \cdot t^{N/2} \quad \text{where  }J_N:= J_N(1).
\end{equation}
Moreover,
\begin{equation}
2^N\leq  J_N \leq  \frac{(N+1)^{N}\,\pi^{N/2}}{(N+1)!}.
\end{equation}
\end{lem}

\begin{pf}
$J_N(t) = J_N \cdot t^{N/2}$ is obvious from \eqref{E:Recursion_J_N(t)}
after a change of variable. Let $\mathcal{M}_N$ be the collection of functions $f: \{ 2,\,3,\cdots,\,N+1 \} \to \{1,\,2,\cdots,\,N \}$ satisfying $f(i)< i$. We can rewrite $J_N(t_1)$ as
\begin{equation}\label{E:Recursion_J_N(t)_1}
J_N(t_1)=\sum_{f\in \mathcal{M}_N}
\int_{t_2=0}^{t_1} \int_{t_3=0}^{t_2}\cdots\int_{t_{N+1}=0}^{t_N} \prod_{i=2}^{N+1}\frac{dt_i}{\sqrt{t_{f(i)}-t_i}}.
\end{equation}
This is a sum of $N!$ terms. When we put $t_1=1$, the smallest term is
$$\int_{t_2=0}^1 \int_{t_3=0}^{t_2}\cdots\int_{t_{N+1}=0}^{t_N}\prod_{j=2}^{N+1}\frac{dt_j}{\sqrt{1-t_{j}}}=\frac{1}{N!}\left(\int_0^1\frac{1}{\sqrt{1-s}}\,ds\right)^N=\frac{2^N}{N!}.$$
Hence we have the lower bound $2^N\leq J_N$. Unfortunately, the largest term is exactly
$$\int_{t_2=0}^1 \int_{t_3=0}^{t_2}\cdots\int_{t_{N+1}=0}^{t_N}\prod_{i=1}^N\frac{dt_j}{\sqrt{t_{j-1}-t_{j}}}=\frac{\pi^{N/2}}{\Gamma(\frac{N+2}{2})}.$$
which grows faster than $C^N$ for any $C\in (0,\infty)$.
Hence for the upper bound, we will employ the recursion formula (\ref{E:Recursion_J_N(t)}). We apply the homogeneity to the right hand side of (\ref{E:Recursion_J_N(t)}) to obtain
$$J_N = \sum_{k= 1}^N \frac{1}{k!} \sum_{\substack{ n_2+n_3 + \cdots + n_{k+1} = N \\  n_2, \cdots, n_{k+1} \geq 1 }} \;\prod_{j=2}^{k+1} J_{n_j-1} \int_{0}^1 \frac{t_j^{(n_j-1)/2}}{\sqrt{1-t_j}}\, d t_j.$$
The integrals are now simple one dimensional and can be evaluated:
\begin{equation}\label{E:Recursion_J_N}
J_N = \sum_{k= 1}^N \frac{1}{k!} \sum_{\substack{ n_2+n_3 + \cdots + n_{k+1} = N \\  n_2, \cdots, n_{k+1} \geq 1 }} \;\prod_{j=2}^{k+1} \left( J_{n_j-1} \cdot \frac{\sqrt{\pi}\,\Gamma((n_j+1)/2)}{\Gamma((n_j+2)/2)} \right).
\end{equation}

Since $\frac{\sqrt{\pi}\, \Gamma((n_j+1)/2)}{\Gamma((n_j+2)/2)} \leq \sqrt{\pi}$, we have $J_N\leq K_N$, where $K_N$ is defined by the recursion
$$
K_N = \sum_{k=1}^N \frac{1}{k!} \sum_{\substack{ n_2+n_3 + \cdots + n_{k+1} = N \\  n_2, \cdots, n_{k+1} \geq 1 }} \;\prod_{j=2}^{k+1} \left( K_{n_j-1} \cdot \sqrt{\pi}\right)
$$
with $K_0=1$. The generating function $\phi(x) := \sum_{N=0}^{\infty} K_N x^N$ of $K_N$ clearly satisfies $\phi(x) = \exp(\sqrt{\pi} x \phi(x))$. We thus see that $\phi(x) = W(-\sqrt{\pi}x)/(-\sqrt{\pi}x)$, where $W$ is the Lambert $W$ function. By Lagrange inversion Theorem (see Theorem 5.4.2 of \cite{rpS99}), $W(z)=\sum_{k=1}^{\infty}(-k)^{k-1}z^k/k!$ (for $|z|<1/e$). Hence by comparing coefficients in the series expansion of $\phi(x)$, we have $K_N = \frac{(N+1)^{N}\,\pi^{N/2}}{(N+1)!}$ as desired.
\end{pf}

\begin{remark} \rm
By Stirling's formula, $ \frac{(N+1)^{N}\,\pi^{N/2}}{(N+1)!}\sim (\sqrt{\pi}e)^N$ (where $a(N)\sim b(N)$ means $\lim_{N\to\infty}\frac{a(N)}{b(N)}=1$). Hence $J_N \leq C^N$ for some $C>0$. Monte Carlo simulations suggests that $J_N \sim \pi^N$.
The recursion (\ref{E:Recursion_J_N}) also makes it clear that $J_N$'s are all in $\mathbb{Q}[\pi]$ (polynomials in $\pi$ with rational coefficients) and makes it easy to compute them recursively. This is because $\frac{\sqrt{\pi}\,\Gamma((n_j+1)/2)}{\Gamma((n_j+2)/2)}$ is rational if $n_j$ is odd and is a rational multiple of $\pi$ if $n_j$ is even. For example, $J_1=2$, $J_2=2+\pi$ and $J_3=4+\frac{10\pi}{3}$.  \qed
\end{remark}

We now turn  to the proof the recursion formula (\ref{E:Recursion_J_N(t)}) which is restated in the following lemma.
\begin{lem}\label{L:Recursion_J_N(t)}
\begin{equation}
J_N(t) = \sum_{k=1}^N  \frac{1}{k!} \sum_{ \substack{n_2+n_3 + \cdots + n_{k+1} = N \\ n_2, n_3, \cdots, n_{k+1} \geq 1} }\;\prod_{j=2}^{k+1} \int_{0}^t \frac{J_{n_j-1}(t_j)}{\sqrt{t-t_j}}\,dt_j
\end{equation}
provided that we set $J_0(t)=1$.
\end{lem}

\begin{pf}
The proof is based on standard combinatorial methods for working on sums over planar rooted trees (see \cite{rpS99}).

\textbf{Step 1: Summing over labeled trees. }Recall (\ref{E:Recursion_J_N(t)_1}). There are $N!$ elements in $\mathcal{M}_N$. We can visualize each of them as a rooted tree with vertex set $\{1,\,2,\,\cdots,\,N+1\}$ and a directed edge  from $i$ to $f(i)$ for each $i$. For example, the 6 elements of $\mathcal{M}_3$ can be represented by

\begin{center}
\begin{tikzpicture}
[edge from parent/.style={draw, <-, thick},node distance=18mm, level distance=14mm]
\tikzstyle{level 1}=[sibling distance=12mm]
\node (a)  {1}
child  {node {2}}
child  {node {3}}
child  {node {4}};
\hskip 0.2truein

\node (b) [right=of a] {1}
child {node {2}
child {node {4}}
}
child {node {3}};

\node (c) [right=of b] {1}
child {node {2}}
child {node {3}
child {node {4}}
};

\node (d) [right=of c] {1}
child {node {2}
child {node {3}}
}
child {node {4}};
\node (e) [right=of d] {1}
child {node {2}
child {node {3}}
child {node {4}}
};

\node (f) [right=of e] {1}
child {node {2}
child {node {3}
child {node {4}}
}};
\end{tikzpicture}
\end{center}

The trees are drawn so that arrows point upwards and the children of a given vertex are listed from left to right. Note that the second and forth tree of the list are the same up to relabeling the vertices.   The idea is to group terms in (\ref{E:Recursion_J_N(t)_1}) like this together. First, we rewrite (\ref{E:Recursion_J_N(t)_1}) in terms of trees. Let $\mathcal{D}_N$ be the set of `decreasing trees', which are trees whose vertices are labeled by $\{1,\,2,\,\cdots,\,N+1\}$ and such that $i<j$ whenever there is an edge $i\leftarrow j$. Then
\begin{equation}\label{E:Recursion_J_N(t)_2}
J_N(t_1)=\sum_{T\in \mathcal{D}_N} \int_{t_1 \geq t_2 \geq \cdots \geq t_{N+1} \geq 0} \prod_{(i\leftarrow j) \in Edge(T)} \frac {dt_j}{\sqrt {t_i - t_j}}.
\end{equation}

\textbf{Step 2: Summing over unlabeled trees. }A planar tree is a rooted unlabeled tree where, for each vertex, the children of that vertex are ordered. We draw a planar tree so that its children are ordered from left to right. Here are the 5 planar rooted trees on 3+1 vertices:

\begin{center}
\begin{tikzpicture}
[edge from parent/.style={draw, <-, thick},node distance=20mm, level distance=14mm]
\tikzstyle{level 1}=[sibling distance=12mm]
\node (a)  {$\bullet$}
child  {node {$\bullet$}}
child  {node {$\bullet$}}
child  {node {$\bullet$}};

\hskip 0.2truein

\node (b) [right=of a] {$\bullet$}
child {node {$\bullet$}
child {node {$\bullet$}}
}
child {node {$\bullet$}};
\node (c) [right=of b] {$\bullet$}
child {node {$\bullet$}}
child {node {$\bullet$}
child {node {$\bullet$}}
};
\node (d) [right=of c] {$\bullet$}
child {node {$\bullet$}
child {node {$\bullet$}}
child {node {$\bullet$}}
};
\node (e) [right=of d] {$\bullet$}
child {node {$\bullet$}
child {node {$\bullet$}
child {node {$\bullet$}}
}};
\end{tikzpicture}
\end{center}
Let $\mathcal{T}_{k}$ be the set of planar rooted trees with $k$ vertices. In general, there are $\frac{(2N)!}{N!(N+1)!}$ (the Catalan number) elements in $\mathcal{T}_{N+1}$, see exercise 6.19 in \cite{rpS99}. We now group all the integrals in (\ref{E:Recursion_J_N(t)_2}) with the same planar tree. For example, two different labeled trees (the second and forth in our list of labeled trees) both give the same unlabeled planar tree (which is the second in the above list). We redraw this unlabeled planar tree $T_0$ below and attach letters $\{a,\,b,\,c,\,d\}$ to $T_0$ for later use.
\begin{center}
\begin{tikzpicture}
[edge from parent/.style={draw, <-, thick}]
\node {a}
child  {node {b}
child  {node {d}}
}
child  {node {c}};
\end{tikzpicture}
\end{center}

The integrands corresponding to the second and forth  labeled trees are
\begin{equation*}
\frac{dt_2 \,dt_3\, dt_4}{\sqrt{(t_1 - t_2)(t_2 - t_4)(t_1 - t_3)}} \quad\text{and}\quad  \frac{dt_2\, dt_3 \,dt_4}{\sqrt{(t_1 - t_2)(t_2 - t_3)(t_1 - t_4)}}.
\end{equation*}
They are the same to
\begin{equation*}
\frac{dt_b \,dt_c \,dt_d}{\sqrt{(t_a - t_b)(t_a - t_c)(t_b - t_d)}},
\end{equation*}
once we relabel the variables by the vertices of $T_0$. That is, $(1, 2, 3, 4)\rightarrow (a, b, c, d)$ for the first term and $(1, 2, 3, 4) \rightarrow (a, b, d, c)$ for the second.

We now go back and keep track of the bounds of integration. In the first integral, they are $t_a \geq t_b \geq t_c \geq t_d$ and, in the second integral, they are $t_a \geq t_b \geq t_d \geq t_c$. We can group these together as
\begin{equation*}
t_a \geq t_b ,  \quad  t_a \geq t_c ,   \quad  t_b \geq t_d , \quad t_b>t_c,
\end{equation*}
which is the same as $t_a\geq t_b\geq t_c$ and $t_b\geq t_d$.

In general, the inequality constraints we have are of two types. First, whenever we have an edge $u \leftarrow v$, we get the inequality $t_u \geq t_v$. Second, if $v$ and $w$ are children of $u$ with $v$ to the left of $w$, then $t_v \geq t_w$. Let $P(T, t_1)$ be the polytope cut out by these inequalities where $t_1$ is the variable at the root. We have proved
\begin{equation}\label{E:Recursion_J_N(t)_3}
J_N(t_1) = \sum_{T\in \mathcal{T}_{N+1}} \int_{P(T, t_1)} \prod_{(u \leftarrow v) \in Edge(T)} \frac{dt_v}{\sqrt{t_u - t_v}},
\end{equation}

\textbf{Step 3: Grouping terms for which the root has degree $k$. }
We abbreviate
$$\omega(T, t_1) := \prod_{(u \leftarrow v) \in Edge(T)} \frac{dt_v}{\sqrt{t_u - t_v}}$$
if $T$ has more than one vertex (otherwise $\omega(T, t_1):= 1$).
Then (\ref{E:Recursion_J_N(t)_3}) translates into
\begin{equation}\label{e:3.47}
J_N(t_1) = \sum_{T\in \mathcal{T}_{N+1}} \int_{P(T, t_1)} \omega(T, t_1).
\end{equation}

Fix an integer $k$ and let $T$ be a tree whose root has degree $k$. Removing the root leaves behind
$k$ children, denoted in chronical order by $t_2, \cdots t_{k+1}$,
and $k$ planar subtrees $T_j$ having $t_j$ as its root for $2\leq j\leq k+1$.
Then
\begin{equation*}
\int_{P(T, t_1)} \omega(T, t_1) = \int_{t_1 \geq t_2 \geq \cdots \geq t_{k+1}\geq 0} \prod_{j=2}^{k+1} \frac {dt_j}{\sqrt {t_1 - t_j}} \int_{P(T_j, t_j)} \omega(T_j, t_j).
\end{equation*}
Hence, group together the terms where the root has degree $k$, we have
\begin{eqnarray}\label{E:Recursion_J_N(t)_4}
J_N(t_1)
&=& \sum_{k=1}^N \sum_{T_2, \cdots, T_{k+1}} \int_{t_1 \geq t_2 \geq \cdots \geq t_{k+1}\geq 0} \prod_{j=2}^{k+1} \frac {dt_j}{\sqrt {t_1 - t_j}} \int_{P(T_j, t_j)} \omega(T_j, t_j).
\end{eqnarray}
Here the summation conditions include that $\sum_{j=2}^{k+1}|T_j|=N$ and each different ordering of $(T_2, \,T_3,\, \cdots ,\, T_{k+1})$ are considered to be different, where $|T_j|$ is the number of vertices in $T_j$. This abbreviation applies whenever $\sum_{T_2, \, \cdots,\, T_{k+1}}$ appears.

On other hand, we have by applying \eqref{e:3.47}
 to each $J_{n_j-1}(t_j)$ below that
\begin{eqnarray}\label{E:Recursion_J_N(t)_5}
&& \sum_{k=1}^N  \frac{1}{k!} \sum_{ \substack{n_2+n_3 + \cdots + n_{k+1} = N \\  n_2, \cdots, n_{k+1} \geq 1} }\;\prod_{j=2}^{k+1} \int_0^{t_1}  \frac{J_{n_j-1}(t_j)}{\sqrt{t_1-t_j}}\,dt_j  \notag \\
&=&  \sum_{k=1}^N \frac{1}{k!} \sum_{\substack{n_2+n_3 + \cdots + n_{k+1} = N \\  n_2, \cdots, n_{k+1} \geq 1}}\;\prod_{j=2}^{k+1} \sum_{| T_j| = n_j} \int_0^{t_1} \frac {dt_j}{\sqrt {t_1 - t_j}} \int_{P(T_j, t_j)} \omega(T_j, t_j) \notag\\
&=&  \sum_{k=1}^N \frac{1}{k!} \sum_{T_2, \cdots, T_{k+1}} \prod_{j=2}^{k+1} \int_0^{t_1}  \frac {dt_j}{\sqrt {t_1 - t_j}} \int_{P(T_j, t_j)} \omega(T_j, t_j) \notag\\
&=&  \sum_{k=1}^N \frac{1}{k!} \sum_{T_2, \cdots, T_{k+1}} \int_{[0, t_1]^{k+1}} \prod_{j=2}^{k+1} \frac{dt_j}{\sqrt {t_1 - t_j}} \int_{P(T_j, t_j)} \omega(T_j, t_j).
\end{eqnarray}

\textbf{Step 4: Identifying the integrals. }
It remains to show that (\ref{E:Recursion_J_N(t)_4}) is equal to (\ref{E:Recursion_J_N(t)_5}).
Let ${\cal S}_k$ denote the space of permutations of $\{2,\,3,\,\cdots,\,k+1\}$.
Then the right hand side of \eqref{E:Recursion_J_N(t)_5} is equal to
\begin{eqnarray*}
&& \sum_{k=1}^N \frac{1}{k!} \sum_{T_2, \cdots, T_{k+1}}
\sum_{\sigma\in {\cal S}_k} \int_{t_1\geq t_{\sigma (2)}\geq \cdots \geq t_{\sigma (k+1)}\geq 0} \prod_{j=2}^{k+1} \frac{dt_{\sigma (j)}}{\sqrt {t_1 - t_{\sigma (j)}}} \int_{P(T_{\sigma(j)}, t_{\sigma (j)})} \omega( T_{\sigma (j)}, t_{\sigma (j)})  \\
&=&\sum_{k=1}^N \frac{1}{k!} \sum_{\sigma \in {\cal S}_k} \sum_{T_2, \cdots, T_{k+1}}
\int_{t_1\geq s_2\geq \cdots \geq s_{k+1}\geq 0} \prod_{j=2}^{k+1} \frac{ds_j}{\sqrt {t_1 - s_j}} \int_{P(T_{\sigma(j)}, s_j)} \omega( T_{\sigma (j)}, s_j)  \\
&=&\sum_{k=1}^N \frac{1}{k!} \sum_{\sigma \in {\cal S}_k} \sum_{T_2, \cdots, T_{k+1}}
\int_{t_1\geq s_2\geq \cdots \geq s_{k+1}\geq 0} \prod_{j=2}^{k+1} \frac{ds_j}{\sqrt {t_1 - s_j}} \int_{P(T_j, s_j)} \omega( T_j, s_j)  \\
&=& \sum_{k=1}^N \sum_{T_2, \cdots, T_{k+1}}
\int_{t_1\geq t_2 \geq \cdots \geq t_{k+1}\geq 0} \prod_{j=2}^{k+1} \frac{dt_j}{\sqrt {t_1 - t_j}} \int_{P(T_j, t_j)} \omega( T_j, t_j),
\end{eqnarray*}
which is $J_N(t)$ by \eqref{E:Recursion_J_N(t)_4}. This completes the  proof of the lemma.
\end{pf}

\subsubsection{Proof of uniqueness}

\noindent{\it Proof of Theorem \ref{T:Uniqueness_hierarchy}}.
By Lemma \ref{L:beta_I_N}, Lemma \ref{L:I_N_J_N} and Lemma \ref{L:Rate_J_N}, we have
\begin{align}
\|\beta_t^{(n,m)}\|_{(n,m)} \leq C_1(n,m)\;C_2(D_+,D_-,T)^N\;t^{N/2}
\end{align}
for all $t\in[0,T]$ and $N\in\mathbb{N}$. This implies that there is a constant $\tau>0$ so that
 $\|\beta_t^{(n,m)}\|_{(n,m)}=0$ for $t\leq \tau$ and  for all $(n,m)\in \mathbb{N}\times \mathbb{N}$.  Note that $\Tilde{\beta}_t := \beta_{\tau+t}$ also satisfies the hierarchy (\ref{E:hierarchy_2}), and that $\Tilde{\beta}_0=0$. Using the hypothesis $\|\beta_t^{(n,m)}\|_{(T,n,m)}\leq C^{n+m}$, we can extend to obtain $\|\beta_t^{(n,m)}\|_{(n,m)}=0$ for $t\in[0,T]$.
\qed

\section{Hydrodynamic Limits}

This section is devoted to the proof of Theorem \ref{T:conjecture}. Throughout this section, we assume the conditions of Theorem \ref{T:conjecture} hold.

\subsection{Constructing martingales}

Since $\eta_t=\eta^{\eps}_t$ has a finite state space, we know that for all bounded function $F:\R_+\times E^\eps\rightarrow \R$ that is smooth in the first coordinate with $\sup_{(s,x)}\big|\frac{\partial F}{\partial s}(s,x)\big|<C<\infty$, we have two $\mathfrak{F^{\eta}_t}$-martingales below:
\begin{equation} \label{E:mtg_MF}
M(t) := F(t,\eta_t)-F(0,\eta_0)-\int_0^t{\dfrac{\partial F}{\partial s}(s,\eta_s)+\mathfrak{L}F(s,\cdot)(\eta_s)ds}
\end{equation}
and
\begin{equation} \label{E:mtg_NF}
N(t) := M(t)^2-\int_0^t{\mathfrak{L}(F^2(s,\cdot))(\eta_s)-2F(s,\eta_s)\mathfrak{L}F(s,\cdot)(\eta_s)}ds ,
\end{equation}
where $\mathfrak{L}=\mathfrak{L}^\eps$ is the generator defined in \eqref{E:generator0}.
See Lemma 5.1 (p.330) of \cite{cKcL98} or Proposition 4.1.7 of \cite{EK86} for a proof. We will use this fact to construct some  important martingales in Lemma \ref{L:OderOfM} below.

In general, suppose $X=(X_t)_{t\geq 0}$ is a CTRW in a finite state space $E$, whose one step transition probability is $p_{xy}$ and mean holding time at $x$ is $h(x)$. Its infinitesimal generator of $X$ is the discrete operator
$$\mathcal{A}f(x) := \frac{1}{h(x)}\sum_{y\in E}p_{xy}(f(y)-f(x)).$$
The formal adjoint $\mathcal{A}^{*}$ of $\mathcal{A}$ is defined by
$$
\mathcal{A}^{*}f(x) := \sum_{y\in E}\Big(\frac{1}{h(y)}p_{yx}f(y)-\frac{1}{h(x)}p_{xy}f(x)\Big).
$$
It can be easily checked that
\begin{equation}\label{E:Dual_Discrete_Generator}
\<f, \,\mathcal{A}^{*}g\>=\<\mathcal{A}f,\, g\>,\quad \text{where  } \<f,g\>:= \frac{1}{N}\sum_{x \in E}{f(x)g(x)} .
\end{equation}

We denote by $\mathcal{A^{\pm}_{\eps}}$  the generator of the CTRW $X^{\pm,\eps}$ on $D^\eps_\pm$, respectively, and by $\mathcal{A^{*,\pm}_{\eps}}$  the corresponding formal adjoint. In our case, $h(x)=h_\eps (x)=\eps^2/d$ for all $x$. We can check that if $f\in C^2(D_{\pm})$, then
\begin{equation}
\lim_{\eps\to 0}\mathcal{A^{\pm}_{\eps}}f(x^{\eps}) = \mathcal{A^{\pm}} f(x) \text{  whenever  }x^{\eps}\in D^{\eps}_{\pm}\text{ converges to } x\in D_{\pm}.
\end{equation}

\begin{lem} \label{L:OderOfM}
For any $\phi \in\mathcal{B}_b(D_+)$,
\begin{equation}\label{E:discrete}
M(t):= M^{+,N}_{\phi}(t) := \<\phi,\,\X^{N,+}_t\> - \<\phi,\,\X^{N,+}_0\>
-\int_0^t{\<\mathcal{A}^{+}_{\eps} \phi,\, \X^{N,+}_s\>}ds + \lambda\,\int_0^t\<J^{N,+}_s,\,\phi\>\,ds
\end{equation}
is an $\mathfrak{F^{\eta}_t}$-martingale for $t\geq 0$, where
$J^{N,+}$ is the measure-valued process defined by \eqref{e:2.34}.
 Moreover, if $\phi\in C^{1}(\bar{D}_{+})$, then there is a constant $C>0$ independent of $N$ so that
 for every  $T>0$,
\begin{equation}\label{e:4.6}
\E\Big[\sup_{t\in[0,T]}M^2(t)\Big]\leq \frac{C\,T}{N}.
\end{equation}
Similar statements hold for $\X^{N,-}$.
\end{lem}

\begin{pf}
The lemma follows   by applying (\ref{E:mtg_MF}) and  (\ref{E:mtg_NF}) to the function
$$F(s,\eta):=f(\eta):=\frac{1}{N} \sum_{x \in D_+}{\phi^+(x)\eta^+(x)}\quad,\,(s,\eta)\in [0,\infty)\times E^{\eps}.$$
We spell out the details here for completeness.
Observe that $f(\eta_t)=\<\phi,\X^{N,+}_t\>$.
Fix $x_0 \in D^{\eps}_+$, and define $\eta^+_{x_0}$ to be the function from $E^{\eps}$ to $\R$ which maps $\eta$ to $\eta^+(x_0)$. Then by the definition of $\mathfrak{L}=\mathfrak{L}^{\eps}$ in (\ref{E:generator0}),
\begin{equation}\label{E:generator1}
  \mathfrak{L}\eta^+_{x_0}(\eta)=\left(\mathcal{A^{*,+}_{\eps}}\eta^+\right)(x_0)-
  \sum_{\{z\in I^{\eps}:\,z_+=x_0\}}\,\dfrac{\lambda}{\eps}\,\Psi(z)\, \eta^+_t(z_+)\eta^-_t(z_-) \, .
\end{equation}
Similarly, for $y_0 \in D^{\eps}_-$, we have
\begin{equation}\label{E:generator2}
  \mathfrak{L}\eta^-_{y_0}(\eta)=\left(\mathcal{A^{*,-}_{\eps}}\eta^-\right)(y_0)-\sum_{\{z\in I^{\eps}:\,z_-=y_0\}}\,\, \dfrac{\lambda}{\eps}\, \Psi(z)\,\eta^+_t(z_+)\eta^-_t(z_-) .
\end{equation}
Hence, by linearity of $\mathfrak{L}$, (\ref{E:generator1}) and then (\ref{E:Dual_Discrete_Generator}), we have
\begin{eqnarray*}
 \mathfrak{L}f(\eta)
&=&\dfrac{1}{N}\sum_{x \in D_+}{\phi(x)(\mathfrak{L}\eta^+_{x})(\eta)}\\
&=&\dfrac{1}{N}\sum_{x \in D_+}{\phi(x)(\mathcal{A^{*,+}_{\eps}}\eta^+(x))}-\dfrac{\lambda}{N\eps}\sum_{z\in I^{\eps}} \,\Psi(z)\, \eta^+(z_+)\eta^-(z_-)\,\phi(z_+)\\
&=&\dfrac{1}{N}\sum_{x \in D_+}{\eta^+(x)(\mathcal{A^{+}_{\eps}}\phi(x))}-\dfrac{\lambda}{N\eps}\,\sum_{z\in I^{\eps}} \,\Psi(z)\,\eta^+(z_+)\eta^-(z_-)\,\phi(z_+) .
\end{eqnarray*}
Hence
\begin{equation}\label{E:formula_for_generator}
\mathfrak{L}f(\eta_s)=\<\mathcal{A^{+}_{\eps}}\phi, \X^{N,+}_s\>- \dfrac{\lambda}{N\,\eps^d}\,\<J^{N,+}_s,\,\phi\>
\end{equation}
and $M(t)$ is an $\mathfrak{F^{\eta}_t}$-martingale by (\ref{E:mtg_MF}).
  Next, we compute $\E[\<M\>_t]$. Note that
  $$\L (f^2)(\eta)= \dfrac{1}{N^2}\sum_{a\in D^{\eps}_+}\sum_{b\in D^{\eps}_+}\phi^+(a)\phi^+(b)\,\L(\eta_a\eta_b)(\eta),$$
  where $\L(\eta_a\eta_b)$ can be computed explicitly using (\ref{E:generator0}). Hence from (\ref{E:mtg_NF}), we can check that $$\E[M^2(t)]=\E[\<M\>_t]=\E\left[\int_0^t{\mathfrak{L}(f^2)(\eta_s)-2f(\eta_s)\mathfrak{L}f(\eta_s)}\,ds\right]=\int_0^t{\E[g(\eta_r)]}dr,$$ where
  \begin{eqnarray} \label{E:OderOfM}
     g(\eta)
    &=& \dfrac{1}{N^2}\,\left(\sum_{y,z \in D_+}\eta^+(z)h^{-1}(z)p_{zy}(\phi(y)-\phi(z))^2 + \dfrac{\lambda}{\eps} \sum_{z\in I^{\eps}} \,\Psi(z)\, \eta(z_+)\eta^-(z_-)\,(-\phi(z_+))^2 \right) \notag \\
    &\leq& \dfrac{1}{N^2}\,\left(\eps^2\|\nabla\phi\|_{\infty}^2 \sum_{x\in D_+}\eta(x)h^{-1}(x) + \dfrac{\lambda}{\eps}\, \|\phi\|_{\infty}^2 \sum_{z\in I^{\eps}} \,\Psi(z)\, \eta^+(z_+)\eta^-(z_-) \right) \notag\\
    &\leq& \frac{d\,\|\nabla\phi\|_{\infty}^2}{N}+ \frac{\lambda\,\|\phi\|_{\infty}^2}{N^2}\left(\frac{1}{\eps}\sum_{z\in I^{\eps}} \,\Psi(z)\, \eta^+(z_+)\eta^-(z_-)\right) . \notag
  \end{eqnarray}
  After taking expectation for $g(\eta_r)$, the first term is the last display is of order at most $1/N$ since $\phi\in C^{1}(\bar{D}_{+})$, while the second term inside the bracket is of order at most $1/N$, uniformly in $r\in[0,t]$, by (\ref{E:Jump_SecondMoment}). Hence $\E[M^2(t)]\leq \frac{C}{N}$ for some $C=C(\phi,d,D_{\pm},\lambda)$.   Doob's maximal inequality
then gives \eqref{e:4.6}.
\end{pf}

\begin{remark} \rm
From the second term of (\ref{E:formula_for_generator}), we see that if the parameter of the killing time is of order $\lambda/\eps$, then we need $N\eps^d$ to be comparable to 1.   \qed
\end{remark}

\subsection{Tightness}

The following simple observation is useful for proving tightness when the transition kernel of the process has a singularity at $t=0$. It says that we can break down the analysis of the fluctuation of functionals of a process on $[0,T]$ into two parts. One part is near $t=0$, and the other is away from $t=0$ where we have a bound for a higher moment. Its proof, which is based on the Prohorov's theorem, is simple and is omitted (detail can be found in \cite{wtF14}).
\begin{lem}\label{L:tightness_criteria}
    Let $\{Y_N\}$ be a sequence of real valued processes such that $t\mapsto \int_0^tY_N(r)\,dr$ is continuous on $[0,T]$ a.s., where $T\in[0,\infty)$. Suppose the following  holds:
    \begin{itemize}
        \item[(i)]  There exists $q>1$ such that $\varlimsup_{N\to\infty}\E[\int_h^T\,|Y_N(r)|^q\,dr]<\infty$ for any $h>0$,\\
        \item[(ii)] $\lim_{\alpha\searrow 0}\varlimsup_{N\to\infty}\P(\int_0^{\alpha}\,|Y_N(r)|\,dr>\eps_0)=0$ for any $\eps_0>0$.
    \end{itemize}
    Then $\big\{\int_0^tY_N(r)\,dr; t\in [0, t]\big\}_{N\in\mathbb{N}}$ is tight in $C([0,T],\R)$.
\end{lem}

Here is our tightness result for $\{(\X^{N,+},\,\X^{N,-})\}$. We need Lemma \ref{L:tightness_criteria} in the proof mainly because we do not know if $\varlimsup_{N\to\infty}\E \int_0^T \< \mathcal{A}^{+}_{\eps} \,\varphi_{+},\,\X^{N,+}_s\>^2\,ds$ is finite or not.
\begin{thm}\label{T:tight}
The sequence $\{(\X^{N,+},\,\X^{N,-})\}$ is relatively compact in
$D([0,T],\mathfrak{E})$ and any subsequential limit of the laws of $\{(\X^{N,+},\,\X^{N,-})\}$
carries
on $C([0,T],\mathfrak{E})$. Moreover, for all $\varphi_{\pm}\in C^2(\bar{D}_{\pm})$,
$$\left\{\int_0^t\< J^{N,+}_s,\,\varphi_{\pm}\>\,ds\right\},\quad \left\{\int_0^t\< \mathcal{A}^{+}_{\eps} \,\varphi_{+},\,\X^{N,+}_s\>\,ds\right\} \text{  and}\quad \left\{\int_0^t\<\mathcal{A}^{-}_{\eps} \,\varphi_{-},\,\X^{N,-}_s\>\,ds\right\}$$
are all tight in $C([0,T],\R)$.
\end{thm}

\begin{pf}
   We write $\X^{\pm}$ in place of $\X^{N,\pm}$ for convenience. By Stone-Weierstrass Theorem, $C^2(\bar{D}_{\pm})$ is dense in $C(\bar{D}_{\pm})$ in uniform topology. It suffices to check that $\{(\<\X^{+},\phi^+\>,\,\<\X^{-},\phi^-\>)\}$ is relatively compact in $D([0,T],\R^2)$ for all $\phi^{\pm} \in C^2(\bar{D}_{\pm})$ (cf. Proposition 1.7 (p.54) of \cite{cKcL98}) for this weak tightness criterion).
   By Prohorov's theorem (see Theorem 1.3 and remark 1.4 of \cite{cKcL98}), $\{(\<\X^{+},\phi^+\>,\,\<\X^{-},\phi^-\>)\}$ is relatively compact in $D([0,T],\R^2)$ if (1) and (2) below holds:
  \begin{itemize}
  \item[(1)]    For all $t\in[0,T]$ and $\eps_0>0$, there exists a compact set  $K(t,\eps_0)\subset \R^2 $ such that  $$\sup_{N}\P\left(\,(\<\X^{+}_t,\phi^+\>,\,\<\X^{-}_t,\phi^-\>) \notin K(t,\eps_0)\,\right)<\eps_0;$$
  \item[(2)]    For all $\eps_0>0$,  $$\lim_{\delta\to 0}\varlimsup_{N\to\infty}\P\left(\sup_{\substack{|t-s|<\delta\\0\leq s,t\leq T}}\Big| \Big(\<\X^{+}_t,\phi^+\>,\,\<\X^{-}_t,\phi^-\>\Big)-\Big(\<\X^{+}_s,\phi^+\>,\,\<\X^{-}_s,\phi^-\>\Big) \Big|_{\R^2} >\eps_0\right)=0.
$$
  \end{itemize}

  We first check (1). Since $\phi^{\pm}$ is bounded on $\bar{D}_{\pm}$ and $|\<\X^{+}_t,1\>|\leq 1$ for all $t\in[0,\infty)$, we can always take $K=[-\|\phi^{+}\|_{\infty},\,\|\phi^{+}\|_{\infty}]\times[-\|\phi^{-}\|_{\infty},\,\|\phi^{-}\|_{\infty}]$.

  To verify (2), since $|(x_1,y_1)-(x_2,y_2)|_{\R^2}\leq|x_1-x_2|+|y_1-y_2|$, we only need to focus on $\X^{+}$. By Lemma \ref{L:OderOfM},
    \begin{equation}\label{E:tight}
     \Big|\<\phi,\X^{+}_t\> - \<\phi,\X^{+}_s\>\Big|= \Big|
     \int_s^t{\<\mathcal{A}^{+}_{\eps} \,\phi, \X^{+}_r\>}dr
     -\int_s^t{\dfrac{\lambda}{N\,\eps^d}\,\<J^{N,+}_r,\,\phi\>}\,dr
     + \left(M_{\phi}(t)-M_{\phi}(s)\right)\Big| .
    \end{equation}
  So we only need to verify (2) with $\<\phi,\X^{+}_t\> - \<\phi,\X^{+}_s\>$ replaced by each of the 3 terms on RHS of the above equation (\ref{E:tight}).

  For the first term of (\ref{E:tight}), we apply Lemma \ref{L:tightness_criteria} for the case $q=2$ and $Y_N(r)=\<\mathcal{A}^{+}_{\eps}\,\phi,\,\X^{+}_r\>$. Since $\phi \in C^2(\bar{D}_{+})$, we have $$\sup_{x\in D^{\eps}\setminus \partial D^{\eps}}|\mathcal{A}^{+}_{\eps} \phi(x)| \leq C(\phi) \text{  and  } \sup_{x\in \partial D^{\eps}}|\eps\,\mathcal{A}^{+}_{\eps} \phi(x)| \leq C(\phi)$$
  for some constant $C(\phi)$ which only depends only on $\phi$. Using Lemma \ref{L:q_near_I}, we have
  $$
\E \left[ \<|\A^+_{\eps}\phi|,\,\X^{N,+}_r\> \right] \leq \frac{1}{N}\sum_{i=1}^N\,\sum_{D^{\eps}}|\A^+_{\eps}\phi(\cdot)|p^{\eps,+}(r,x_i,\cdot)m_{\eps}(\cdot)\leq C_1(d,D_+,\phi)+\frac{C_2(d,D_+,\phi)}{\eps\vee r^{\frac{1}{2}}},
$$
  which is in $L^1[0,T]$ as a function in $r$. This implies hypothesis (ii) of Lemma \ref{L:tightness_criteria}, via the Chebyshev's inequality. Hypothesis (i) of Lemma \ref{L:tightness_criteria} can be verified easily using the upper bound (\ref{E:uniformbound}) for the correlation function, or by direct comparison to the process without annihilation:
  \begin{eqnarray*}
  &&\E \left[ \<|\A_{\eps}\phi|,\,\X^{+}_r\>^2 \right] \\
  &\leq& \left(\frac{1}{N}\sum_{i=1}^N\E[\A_{\eps}\phi(X^i_r)]\right)^2
  +\frac{1}{N^2}\sum_{i=1}^N\E \left[ (\A_{\eps}\phi)^2(X^i_r) \right]-\frac{1}{N^2}\sum_{i=1}^N\left(\,\E[\A_{\eps}\phi(X^i_r)]\,\right)^2 \\
  &\leq& C(d,D,\phi)\left(1+\frac{1}{\sqrt{r}}+\frac{1}{r}\right).
  \end{eqnarray*}

  For the second term of (\ref{E:tight}), by (\ref{E:Jump_SecondMoment}) we have $\varlimsup_{N\to\infty}\E\left[\int_0^T\,\<1,\,J^N_r\>^2\,dr\right]<\infty$. Hence (2) holds for this term by Lemma \ref{L:tightness_criteria}.

    For the third term of (\ref{E:tight}), by Chebyshev's inequality, Doob's maximal inequality and Lemma \ref{L:OderOfM}, we have
    \begin{eqnarray*}
    \P\left(\sup_{|t-s|<\delta}|M_{\phi}(t)-M_{\phi}(s)|>\eps_0\right) &\leq& \frac{1}{\eps_0^2}\,\E\left[\left(\sup_{|t-s|<\delta}|M_{\phi}(t)-M_{\phi}(s)|\right)^2\right]\\
    &\leq& \frac{1}{\eps_0^2}\,\E\left[\left(2\sup_{t\in[0,T]}|M_{\phi}(t)|\right)^2\right]\\
    &\leq& 16\E[M_{\phi}(T)^2]\leq \frac{C}{N} .
    \end{eqnarray*}
    We have proved that (2) is satisfied. Hence $\{(\X^{N,+},\,\X^{N,-})\}$ is relatively compact. Using (2) and the metric of $\mathfrak{E}$, we can check that any subsequential limit $\mathrm{L}^{*}$ of the laws of $\{(\X^{N,+},\,\X^{N,-})\}$ concentrates on $C([0,\infty),\mathfrak{E})$.
    \end{pf}

\begin{remark} \rm
In general, to prove tightness for $(X_n,\,Y_n)$ in $D([0,T], A\times B)$, it is NOT enough to prove tightness separately for $(X_n)$ and $(Y_n)$ in $D([0,T], A)$ and $D([0,T], B)$ respectively.  (However, the latter condition implies tightness in $D([0,T], A)\times D([0,T], B)$ trivially). See Exercise 22(a) in Chapter 3 of \cite{EK86}.
For example, $\left(\1_{[1+\frac{1}{n},\infty)},\;\1_{[1,\infty)}\right)$ converges in $D_{\R}[0,\infty)\times D_{\R}[0,\infty)$ but not in $D_{\R^2}[0,\infty)$. The reason is that the two processes can jump at different times ($t=1$ and $t=1+\frac{1}{n}$) that become identified in the limit (only one jump at $t=1$); this can be avoided if one of the two processes is $C$-tight (i.e. has only continuous limiting values), which is satisfied in our case since $\X^{N,+}$ and $\X^{N,-}$ turns out to be both $C$-tight.  \qed
\end{remark}

\begin{remark} \rm
Even without condition (ii) of Theorem \ref{T:conjecture} for $\eta_0$, we can still verify hypothesis (i) of Lemma \ref{L:tightness_criteria}. Actually, applying (\ref{E:discrete}) to suitable test functions, we have
\begin{equation}
\lim_{\alpha\to 0}\varlimsup_{N\to\infty}\E\left[\int_0^{\alpha}\<J^{N,+}_s,\,1\>\,ds\right] =0 .
\end{equation}
\qed
\end{remark}

\subsection{Identifying the limit}

Suppose $(\X^{\infty,+},\,\X^{\infty,-}) $ is a subsequential limit of $(\X^{N,+},\,\X^{N,-})$, say the convergence is along the subsequence $\{N'\}$. By the Skorokhod representation Theorem, the continuity of the limit in $t$ and \cite[Theorem 3.10.2]{EK86}, there exists a  probability space $(\Omega, \mathfrak{F}, \P)$ such that
\begin{equation}
    \lim_{N' \to \infty} \sup_{t \in [0,T]} \Big\| (\X^{+,N'}_t,\,\X^{-,N'}_t)-(\X^{\infty,+}_t,\,\X^{\infty,-}_t) \Big\|_{\mathfrak{E}} =0  \quad \P
\hbox{-a.s},
\end{equation}
Hence we have for any $t>0$ and $\phi\in C(\bar{D}_+)$,
\begin{equation*}
\lim_{N'\to\infty}\,\E[\<\mathfrak{X}^{+}_t,\phi\>]= \E[\<\mathfrak{X}^{+,\infty}_t,\phi\>]  \quad \text{  and  } \quad
\lim_{N'\to\infty}\,\E[(\<\mathfrak{X}^{+}_t,\phi\>)^2]= \E[(\<\mathfrak{X}^{+,\infty}_t,\phi\>)^2].
\end{equation*}
Combining with Corollary \ref{cor:correlation}, we have
$$\<\mathfrak{X}^{+,\infty}_t,\phi\>= \<u_{+}(t),\,\phi\>_{\rho_+} \quad \P \hbox{-a.s. for every } t\geq 0 \hbox{ and for   } \phi.
$$
Here we have used the simple fact that if $\E[X]= (\E[X^2])^{1/2}=a$, then $X=a$ a.s.

Suppose $\{\phi_k\}$ is a countable dense subset of $C(\bar{D}_{+})$. Then for every $t\geq 0$,
$$
\<\X^{\infty,+}_t,\phi_k\>= \<u_{+}(t),\,\phi_k\>_{\rho_+}
\quad \hbox{for every  } k\geq  1 \   \P \hbox{-a.s.}
$$
Since $\X^{\infty,+} \in C((0,\infty),M_+(\bar{D}_{+}))$, we can pass to rational numbers to obtain
$$
\<\X^{\infty,+}_t,\phi_k\>= \<u_{+}(t),\,\phi_k\>_{\rho_+}   \quad  \hbox{for every  $t\geq 0$ and  } k\geq 1   \   \P \hbox{-a.s.}
$$
Hence, $$\X^{\infty,+}_t(dx)= u_{+}(t,x)\,\rho_+(x) dx  \quad \hbox{for every } t\geq 0  \ \quad \P \hbox{-a.s.}
$$
Similarly, $$\X^{\infty,-}_t(dy)= u_{-}(t,y)\,\rho_-(y)dy  \quad \hbox{for every } t\geq 0 \ \quad \P \hbox{-a.s.}
$$
In conclusion, any subsequential limit is the dirac delta measure
$$
\delta_{u_{+}(t,x) dx, u_{-}(t,y)dy}\,\in M_1(D([0,\infty),\mathfrak{E})).
$$
 This together with Theorem \ref{T:tight} completes the proof of  Theorem \ref{T:conjecture}.

\section{Local Central Limit Theorem} \label{S:5}

Suppose $D\subset \R^d$ is a bounded Lipschitz domain and $\rho\in W^{1,2}(D)\cap C^1(\bar{D})$ is strictly positive. Suppose $X$ is a $(I_{d\times d},\rho)$-reflected diffusion and $X^{\eps}$ be an $\eps$-approximation of $X$, described
 in the subsection \ref{UnderlyingMotion}.  In this section, we  prove the local central limit theorem (Theorem \ref{T:LCLT_CTRW}), the Gaussian upper bound (Theorem \ref{T:UpperHKE}) and the H\"older continuity (Theorem \ref{T:HolderCts}) for $p^{\eps}$. The proofs
are  standard once we establish a discrete analogue of a relative isoperimetric inequality (Theorem \ref{T:Isoperimetric_Discrete}) for bounded Lipschitz domains.

\subsection{Discrete relative isoperimetric inequality}

Note that any Lipschitz domain enjoys the uniform cone property and any bounded
$(d-1)$-dimensional manifold of class $C^{0,1}$ has finite perimeter. Hence, by Corollary 3.2.3 (p.165) and Theorem 6.1.3 (p.300) of \cite{vgM85}, we have the following relative isoperimetric inequality.
\begin{prop}[Relative isoperimetric inequality]\label{prop:relative_iso_ineq}
Let $D\subset \R^d$ be a bounded Lipschitz domain and $r\in(0,1]$. Then
\begin{equation}
S(r,D) := \sup_{U\in \mathfrak{G}}\dfrac{|U|^{\frac{d-1}{d}}}{\sigma(\partial U \cap D)} <\infty ,
\end{equation}
where $\mathfrak{G}$ is the collection of open subsets $U\subset D$ such that $|U|\leq r|D|$ and $\partial U \cap D$ is a manifold of class $C^{0,1}$ (i.e. each point $x\in \partial U \cap D$ has a neighborhood which can be represented by the graph of a Lipschitz function). Moreover, $S(r,D)=S(r,aD)$ for all $a>0$.
\end{prop}

In this subsection, we establish a discrete analogue for the relative isoperimetric inequality (Theorem \ref{T:Isoperimetric_Discrete}).

First, we study the scaled graph $aD^{\eps}$ and gather some basic properties of a continuous time random walk on a finite set.

\subsubsection{CTRW on scaled graph $aD^{\eps}$}

Recall the $m_{\epsilon}$-symmetric CTRW $X^{\eps}$ on $D^{\eps}$ defined in the Subsection \ref{UnderlyingMotion}. The Dirichlet form $(\D^{(\epsilon)},\, l^2(m_{\epsilon}))$ of $X^{\eps}$ in $l^2(m_{\epsilon})$ is given by
\begin{eqnarray}\label{E:DirichletForm_CTRW}
  \D^{(\epsilon)} (f,g) :=   \dfrac{1}{2} \sum_{x,y\in D^{\epsilon}}(f(y)-f(x))(g(y)-g(x))\,\mu_{xy},
\end{eqnarray}
where $\mu_{xy}=\mu^{D^{\eps}}_{x,y}$ are the conductance on the graph $D^{\eps}$ defined in the Subsection \ref{UnderlyingMotion}. The stationary measure $\pi=\pi^{D^{\eps}}$ of $X^{\eps}$ is given by $\pi(x)=m_{\eps}(x)/m(D^{\eps})$, where $m(D^{\eps}):=\sum_{x\in D^{\eps}}m_{\eps}(x)$.

We now consider the scaled graph $aD^{\eps}= (aD)^{a\eps}$, which is an approximation to the bounded Lipschitz domain $aD$ by square lattice $a\eps\Z^d$. Clearly the degrees of vertices are given by $v^{aD^{\eps}}{(ax)}=v^{D^{\eps}}(x)$. Define the function $\rho_{(aD)}$ on $aD$ by $\rho_{(aD)}(ax):= \rho(x)$. Then define the CTRW $X^{aD^{\eps}}$ using $\rho_{(aD)}$ as we have done for $X^{\eps}$ using $\rho$.

The mean holding time of $X^{aD^{\eps}}$ is
$(a \eps)^2 /d$.  Clearly, the symmetrizing measure $m^{aD^{\eps}}$ and the stationary probability measure $\pi^{aD^{\eps}}$
have the scaling property $m^{aD^{\eps}}(ax)=a^d\,m^{D^{\eps}}(x)$ and $\pi^{aD^{\eps}}(ax)=\pi^{D^{\eps}}(x)$. Let $p_{aD}^{a\eps}$ be the transition density of $X^{aD^{\eps}}$ with respect to the symmetrizing measure $m^{aD^{\eps}}$. Then
\begin{equation}
a^d\,p_{aD}^{a\eps}(a^2t,ax,ay)=p_{D}^{\eps}(t,x,y)
\end{equation}
for every $t>0$, $\eps>0$, $a>0$ and $x,y\in D^{\eps}$ .

We will simply write $m$ and $\pi$ for the symmetrizing measure and the stationary probability measure when there is no ambiguity for the underlying graph.

\subsubsection{An extension lemma}

Following the notation of \cite{scL97}, we let $G$ be a finite set, $K(x,y)$ be a Markov kernel on $G$ and $\pi$ the stationary measure of $K$. Note that a Markov chain on a finite set induces a natural graph structure as follows.
Let $Q(e) := \frac{1}{2}(K(x,y)\pi(x)+K(y,x)\pi(y))$ for any $e=(x,y)\in G\times G$. Define the set of directed edges $E := \{e=(x,y)\in G\times G:\, Q(e)>0 \}$.

We use the following 2 different notions for the ``boundary'' of $A\subset G$:
\begin{eqnarray*}
\partial_e A &:=& \{ e=(x,y)\in E:\,x\in A,y\in G\setminus A \text{ or } y\in A,x\in G\setminus A \} , \\
\partial A &:=& \{x\in A:\, \exists y\in G\setminus A \text{ such that } (x,y)\in E \} .
\end{eqnarray*}
Observe that each edge in $\partial_e A$ is counted twice. Set
$$ Q(\partial_eA) := \frac{1}{2}\sum_{e\in \partial_eA} q(e)=\frac{1}{2}\sum_{x\in A,\,y\in G\setminus A}(K(x,y)\pi(x)+K(y,x)\pi(y)) .
$$

\begin{definition}\label{Def:IsoperimetricConstant}
For any $r\in(0,1)$, define
\begin{equation}
S_{\pi}(r,G):= \sup_{\{A\subset G: \pi(A)\leq r\}}\dfrac{2 |A|^{(d-1)/d}}{|\partial_eA|} \quad \hbox{ and }
\quad
\tilde{S}_{\pi}(r,G):= \sup_{\{A\subset G: \pi(A)\leq r\}}\dfrac{\pi(A)^{(d-1)/d}}{Q(\partial_e A)}.
\end{equation}
We call $1/\tilde{S}_{\pi}(r,G)$ an \textbf{isoperimetric constant} of the chain $(K,\pi)$. It provides rich information about the geometric properties of $G$ and the behavior of the chain (cf. \cite{scL97}).
\end{definition}

In our case, $G=aD^{\eps}$, $\pi(x)=\frac{m(x)}{m(aD^{\eps})}$ and $K(x,y)=p_{xy}$ in $aD^{\eps}$, where $p_{x,y}$ is the one-step transition probabilities of $X^{aD^{\eps}}$ defined in subsection \ref{UnderlyingMotion}. For $a=1$ and $A\subset D^{\eps}$, we have
\begin{eqnarray*}
\partial_e A &=& \{(x,y)\in A\times D^{\eps}\setminus A \cup D^{\eps}\setminus A \times A:\text{ the line segment } (x,y]\subset D \}, \\
\partial A &=& \{x\in A:\, \exists y\in D^{\eps}\setminus A \text{ such that } |x-y|=\eps \text{ and the line segment }[x,y]\subset D\, \}, \\
\tilde{\partial} A &:=& \{x\in A:\,\exists y\in \eps\Z^d \text{ such that } |x-y|=\eps \text{ and } (x,y]\cap \partial D \neq \emptyset\} , \\
\Delta A &:=& \tilde{\partial}A \setminus \partial A.
\end{eqnarray*}
In this notation, we have $\partial D^{\eps}=\emptyset$, $\tilde{\partial} D^{\eps}=\{x\in D^{\eps}:\,v(x)<2d\}$, $A\cap\tilde{\partial} D^{\eps}=\tilde{\partial}A$ and $\tilde{\partial} D^{\eps}= \Delta A \cup (\partial A\cap \tilde{\partial}D^{\eps}) \cup (\tilde{\partial}D^{\eps} \setminus A)$. See Figure \ref{fig:ExtendDomain2_1} for an illustration.

\begin{definition}
We say that $A\subset D^{\eps}$ is {\it grid-connected} if $\partial_e A_1 \cap \partial_e A_2 \neq \phi$ whenever $A=A_1 \cup A_2$.
\end{definition}
It is easy to check that $A$ is grid-connected if and only if for every $ x,\,y\in A$, there exists  $\{x_1=x,\,x_2,\,\cdots\,,\,x_{m-1},\,x_m=y\}\subset A$ such that each line segments $[x_j,\,x_{j+1}]\subset D$ and $|x_j-x_{j+1}|=\eps$.

The following is a key lemma which allows us to derive the relative isoperimetric inequality for the discrete setting from that in the continuous setting, and hence leads us to Theorem \ref{T:Isoperimetric_Discrete}.

\begin{lem}(Extension of sub-domains)\label{L:ExtendDomain}
Let $\pi_{srw}$ be the stationary measure of the simple random walk (SRW) on $D^{\eps}$.
For any $r\in(0,1)$, there exist positive constants $\eps_1(d,D,r)$, $M_1(d,D,r)$ and $M_2(d,D,r)$ such that if $\eps \in (0,\eps_1)$, then for all grid-connected $A\subset D^{\eps}$ with  $\pi_{srw}(A)\leq r$, we can find a connected open subset $U\subset D$ which contains $A$ and satisfies:
\begin{itemize}
\item [(a)]  $\partial U \cap D$ is $\mathcal{H}^{d-1}$-rectifiable,
\item [(b)]   $|U|\leq \frac{49r+1}{50}|D|$,
\item [(c)]   $\eps^{d}\,|A| \leq M_1\,|U|$,
\item [(d)]   $M_2\,\eps^{d-1}|\partial A| \geq \sigma(\partial U \cap D)$ .
\end{itemize}
\end{lem}

\begin{pf}
Since the proof for each $r\in(0,1)$ is the same, we just give a proof for the case $r=1/2$.

For   $x\in \eps\Z^d$, let  $U_x := \prod_{i=1}^d\,(x_i-\frac{\eps}{2},\,x_i+\frac{\eps}{2})$   be the cube in the dual lattice which contains $x$. Since $A$ is grid-connected, we have $(W_1)^{o}$ is connected in $\R^d$, where $(W_1)^{o}$ is the interior of $W_1:= \cup_{x\in A} (\bar{U_x}\cap D)$. (See Figure \ref{fig:ExtendDomain2_1} for an illustration.)

\begin{figure}[h]
\begin{minipage}{.5\textwidth}
    \centering
     \scalebox{0.4}{\includegraphics{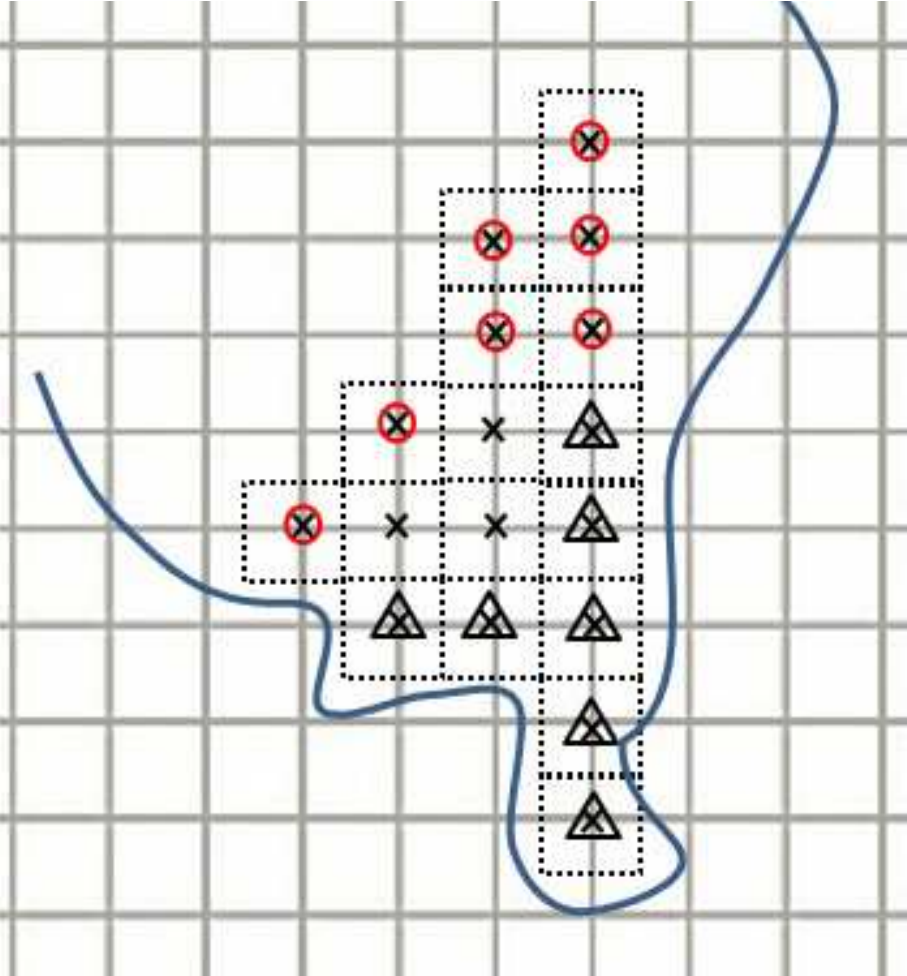}}
     \caption{$W_1:= \cup_{x\in A} (\bar{U_x}\cap D)$ }\label{fig:ExtendDomain2_1}
\end{minipage}
\begin{minipage}{.5\textwidth}
    \centering
     \scalebox{0.4}{\includegraphics{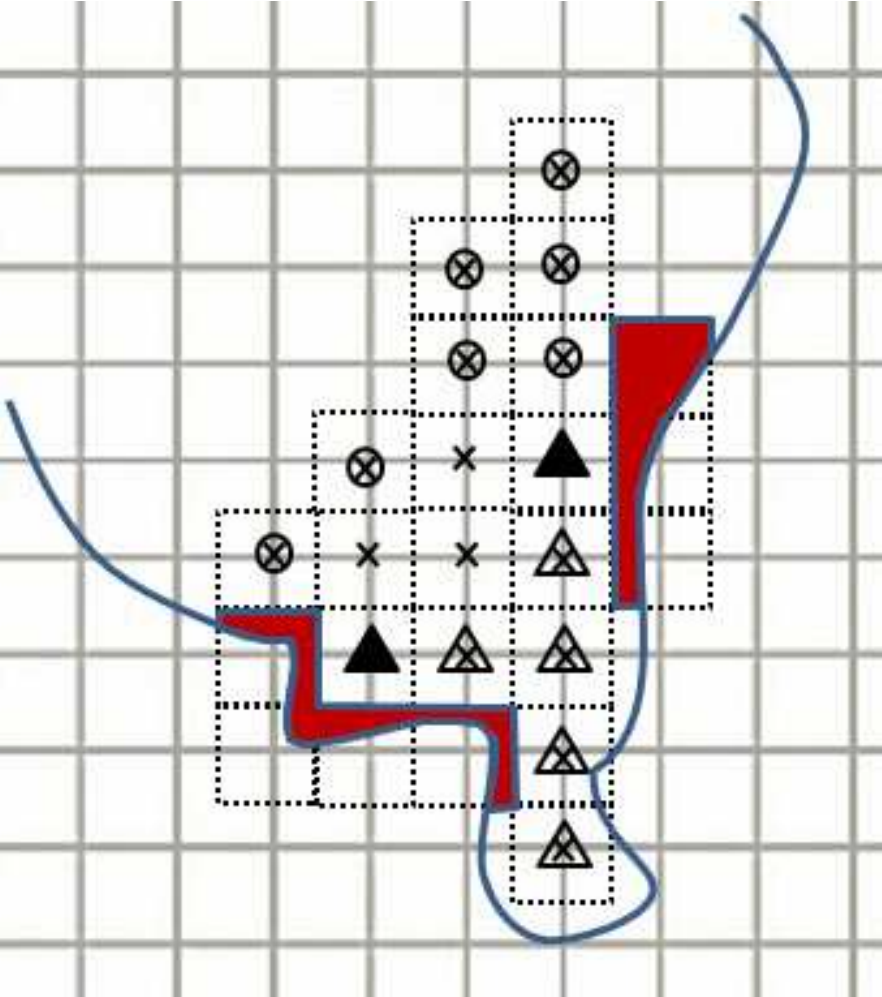}}
     \caption{$W_3$ is the shaded part}\label{fig:ExtendDomain2_2}
\end{minipage}
\end{figure}

Note that we cannot simply take $U=(W_1)^o$  because (d) may fail, for example when $\Delta A$ contributes too much to $\partial U \cap D$,  i.e.,  when $(W_1)_{\Delta}:=\partial W_1\cap \left(\bigcup_{x\in\Delta A}\partial{U_x}\right)$ is large. However, $\Delta A\subset \tilde{\partial}D^{\eps}$ is close to $\partial D$ and so we can fill in the gaps between $\Delta A$ and $\partial D$ to eliminate those contributions. In this process, we may create some extra pieces for $\partial U \cap D$, but we will show that those pieces are small enough. Following this observation, we will eventually take $U=(W_1\cup W_2)^o$ where $W_2\subset D_h$ for some small enough $h>0$.

Since $D$ is a bounded Lipschitz domains, we can choose $h>0$ small enough so that $|D_h|<|D|/200$.
Moreover, $\pi_{srw}(A)<1/2$ implies $\eps^d|A|\leq \eps^d|\partial (D^{\eps})|+m_{srw}(D^{\eps})/2$. So we can choose $\eps$ small enough so that $|W_1|\leq \eps^d|A|\leq \frac{101}{200}|D|$. Hence $U$ satisfies (b). By Lipschitz property again, there exists $M_1>0$ such that $|U_x\cap D|\geq |U_x|/M_1=\eps^d/M_1$ for any $x\in D^{\eps}$. Hence (c) is satisfied.

It remains to construct $W_2$ in such a way that $W_2\subset D_h$ for some small enough $h>0$ (more precisely, for $h$ small enough so that $|D_h|<|D|/200$) and that (a) and (d) are satisfied. We will construct $W_2$ in 3 steps:

Step 1: (Construct $W_3$ to seal the opening between $\partial D$ and the subset of $(W_1)_{\partial}$ which are close to $\partial D$. See Figure \ref{fig:ExtendDomain2_2}.)
Write $\Delta A = \Delta_1A \cup \Delta_2 A$ where $\Delta_2 A:= \Delta A \setminus \Delta_1A$ and
$$\Delta_1A:= \{x\in \Delta A:\,\exists y\in \partial A \text{ such that } \max\{|x_i-y_i|:\,1\leq i \leq d\}=1\}.$$
Points in $\Delta_1A$ are marked in solid black in Figure \ref{fig:ExtendDomain2_2}. For $x\in \Delta_1A$, consider the following closed cube centered at $x$:
$$T_x := \bigcup_{y\in  \wh{B}(x,\,10R\eps)}\bar{U_y}\quad, \text{ where } R=\sqrt{d}(M+1)$$
Let $\Theta_{x}$ be the union of all connected components of $T_{x}\cap D$ whose closure intersects $\bar{U_x}$ and define
$$W_3:=  \bigcup_{x\in\Delta_1A}\,\Theta_{x}.
$$

Step 2: (Fill in the gaps between $\partial D$ and $(W_1)_{\Delta}$ near $\Delta_2$. See Figure \ref{fig:ExtendDomain2_3})
Note that $\cup_{x\in\Delta_1A}\partial U_{x}$ does not contribute to $\partial (W_1\cup W_3)\cap D$.
Let $W_4$ be the union of all connected components of $D\setminus (W_1\cup W_3)$ whose closure intersects $\bar{U_x}$ for some
$x\in \Delta_2 A$.

\begin{figure}[h]
\begin{minipage}{.5\textwidth}
    \centering
     \scalebox{0.4}{\includegraphics{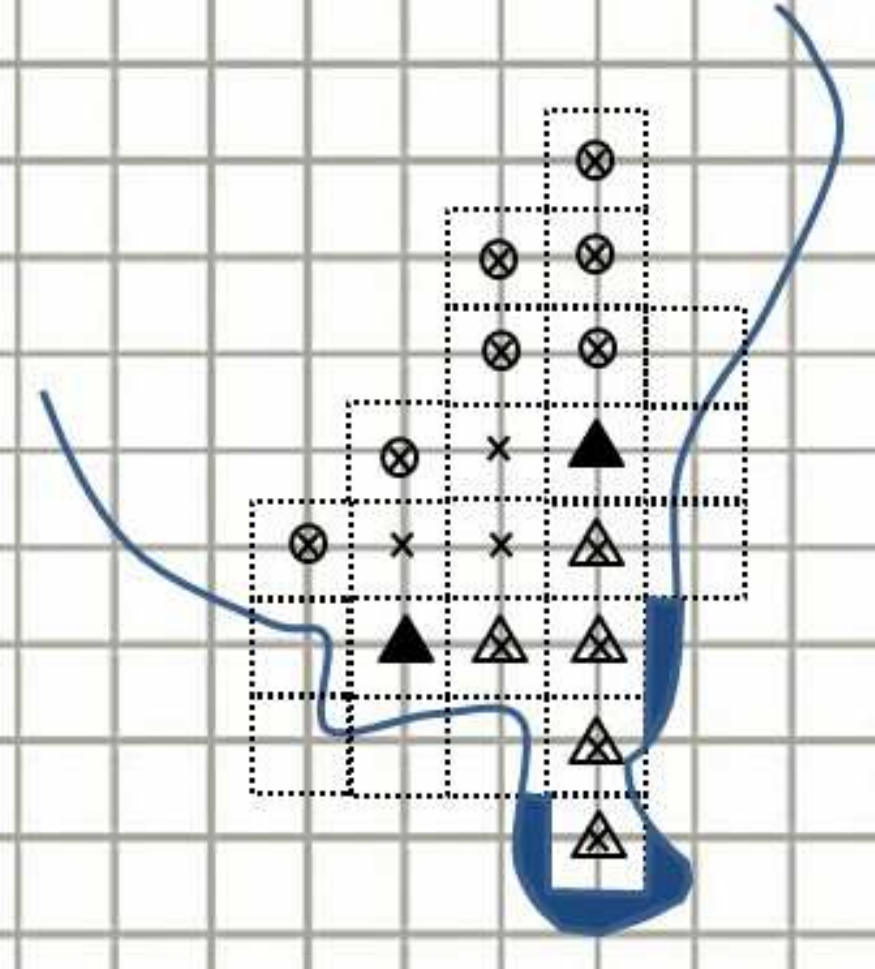}}
     \caption{$W_4$ is the shaded part}\label{fig:ExtendDomain2_3}
\end{minipage}
\begin{minipage}{.5\textwidth}
    \centering
     \scalebox{0.4}{\includegraphics{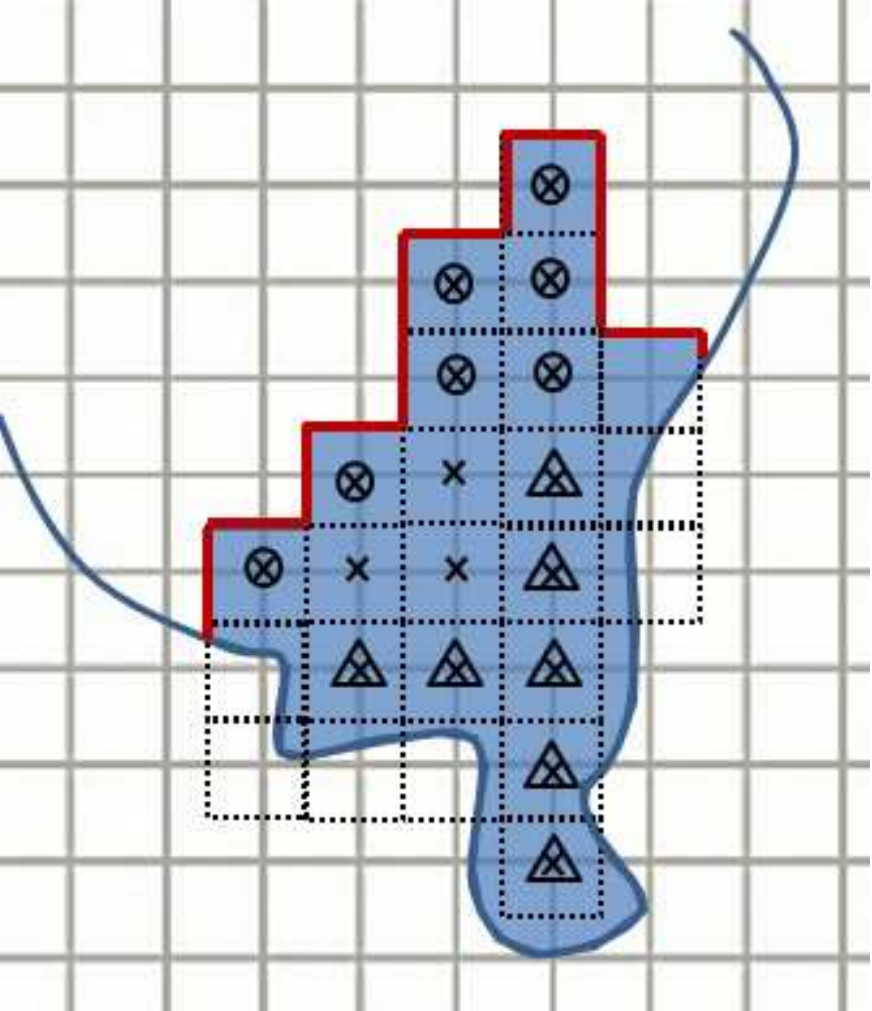}}
     \caption{$U$ is the shaded part}\label{fig:ExtendDomain2_4}
\end{minipage}
\end{figure}

Step 3:
Finally, take $W_2:= W_3\cup W_4$, and set $U:= (W_1\cup W_3\cup W_4)^o$. (See Figure \ref{fig:ExtendDomain2_4}.)

It is clear that $U$ is connected and $\partial U\cap D \subset \bigcup_{x\in \eps \Z^d}\partial U_x$ is piecewise linear, so (a) is satisfied.
For any $W\subset D$, we have $\partial W \cap D=W_{\partial}\cup W_{\Delta} \cup W_{\nabla}$,  where
$$W_{\partial}:=\partial W\cap \left(\bigcup_{x\in\partial A}\partial{U_x}\right), \quad
W_{\Delta}:=\partial W\cap \left(\bigcup_{x\in\Delta A}\partial{U_x}\right) \text{   and    }
W_{\nabla}:=\partial W\setminus \left(\bigcup_{x\in\partial A\cup\Delta A}\partial{U_x}\right).
$$
Therefore, $\sigma(\partial W \cap D)\leq\sigma(W_{\partial})+ \sigma(W_{\Delta})+ \sigma(W_{\nabla})$ whenever the corresponding surface measures are defined.
It is clear that by construction we have
\begin{itemize}
\item   $(W_1)_{\nabla}=\emptyset$,
\item   $(W_1\cup W_3)_{\partial}\subset (W_1)_{\partial},\, (W_1\cup W_3)_{\Delta}\subset (W_1)_{\Delta_2},\,(W_1\cup W_3)_{\nabla}\subset \bigcup_{x\in \Delta_1}\bigcup_{y\in \wh{B}(x,10R\eps)} \partial U_y$ \\ where $(W_1)_{\Delta_2}$ is defined analogously as $(W_1)_{\Delta}$, with $\Delta$ replaced by $\Delta_2$,
\item   $U_{\partial}\subset (W_1\cup W_3)_{\partial},\, U_{\Delta}=\emptyset,\, U_{\nabla}\subset (W_1\cup W_3)_{\nabla}$.
\end{itemize}
Now $\sigma(U_{\partial})\leq\sigma((W_1)_{\partial})\leq |\partial A|\,2d\eps^{d-1}$. Moreover,
each $x\in \partial A$ is adjacent to at most $3^d-1$ points in $\Delta_1A \cup \Delta A$, and for each $x\in \Delta_1A$, there are at most $|\wh{B}(10R\eps)|\leq (20R+1)^d$ cubes in $T_x$. So we have
$$
\sigma(U_{\nabla})\leq\sigma((W_1\cup W_3)_{\nabla})\leq (3^d-1)|\partial A|\,(20R+1)^d\,2d\eps^{d-1} .
$$
Hence (d) is satisfied.

Since $diam\,(T_x)<20R\sqrt{d}\eps$, we have $W_3\subset D_{(20R\sqrt{d}+1)\eps}$. To complete the proof, it suffices to show that $W_4\subset D_{(10R)\eps}$. This is equivalent to show that any curve in $D\setminus \bar{W_1\cup W_3}$ starting from any point in $(W_1)_{\Delta_2}$ must  lie in $D_{(10R)\eps}$.

Let $\gamma[0,1]$ be an arbitrary continuous curve starting at an arbitrary point $p\in (W_1)_{\Delta_2}$ such that $\gamma(0,1) \subset D\setminus \bar{W_1}$ and $dist(\gamma(t),\,\partial D)>(10R)\eps$ for some $t\in(0,1]$. Define $\Omega_{D^{\eps}}:= \left(\bigcup_{x\in D^{\eps}}\bar{U_x}\right)^o\cap D$. Since $(W_1)_{\Delta}\subset\partial(\Omega_{D^{\eps}})\cap D\subset \bigcup_{z\in\tilde{\partial}D^{\eps}}\partial U_z$ and $\sup_{z\in \tilde{\partial} D^{\eps}}dist(z,\,\partial D)<\eps$, the time when $\gamma$ first exits $D\setminus \bar{\Omega_{D^{\eps}}}$ must be less than $t$ by continuity of $\gamma$. That is,
$$
  \tau:= \inf\left\{s>0:\,\gamma(s)\in \left(\bigcup_{z\in\partial D^{\eps}\setminus A}\partial U_z\right)\cap D\right\} <t.
$$
It suffices to show that $\gamma(0,\tau] \cap \Theta_{x}\neq\emptyset$ for some $x\in \Delta_1A$. We  do so by constructing a continuous curve $\tilde{\gamma}$ which is close to $\gamma$ and passes through $\partial U_x$ for some $x\in \Delta_1A$.

                \begin{figure}[h]
            	\begin{center}
            	\includegraphics[scale=.4]{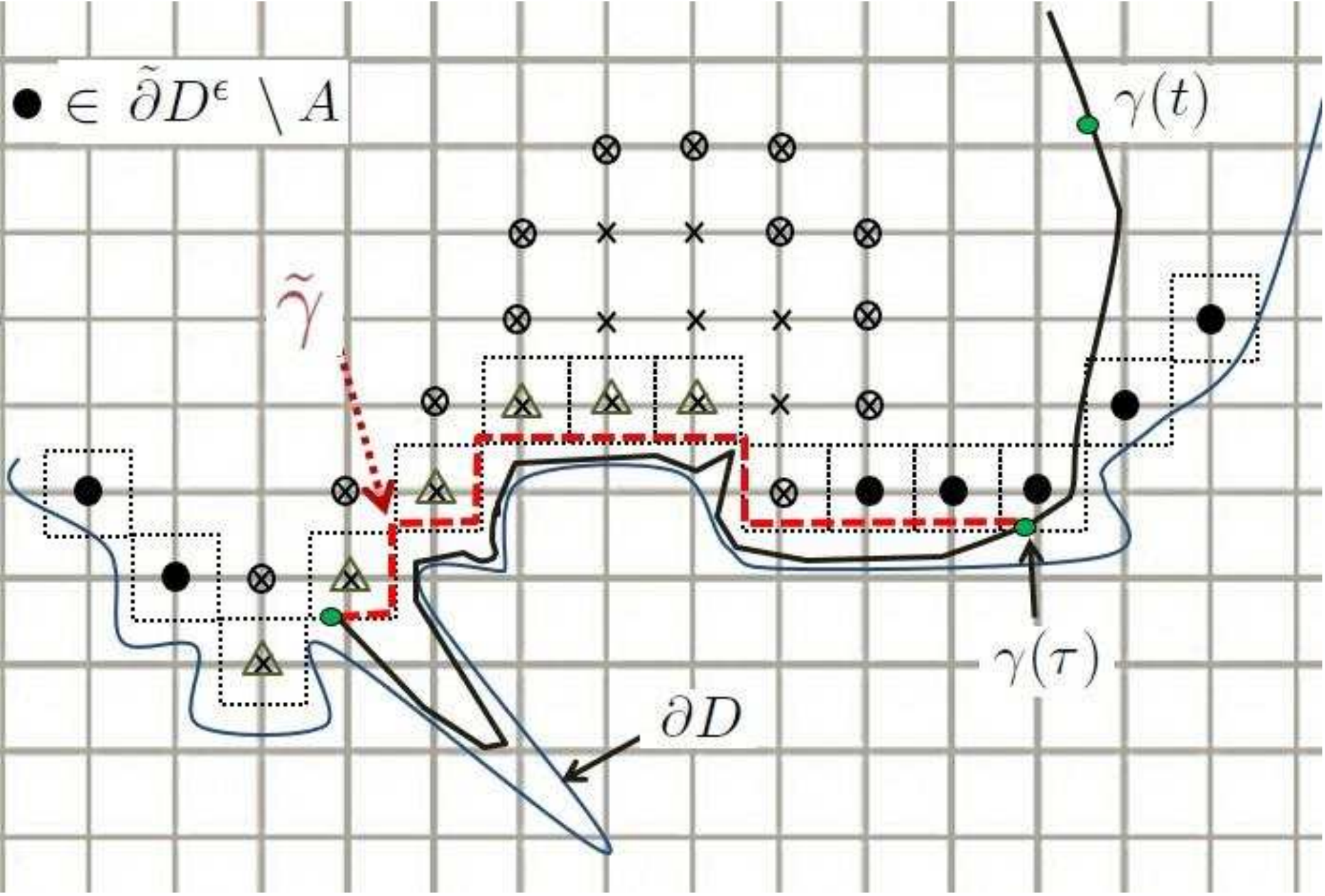}
                \caption{$\gamma$ and a corresponding continuous $\tilde{\gamma}\subset \partial (\Omega_{D^{\eps}})\cap D$}\label{fig:ExtendDomain_curve}
                \vspace{-1em}
            	\end{center}
            	\end{figure}

Since $\sup_{s\in[0,\tau]}dist(\gamma(s),\,\partial D)< 2R\eps$, we can choose $\eps$ small enough (depending only on $D$) and split $[0,\,\tau]$ into finitely many disjoint intervals $I$'s so that the $4R\eps$-tube of each $\gamma(I)$ lies in a coordinate ball $B_{(I)}$ of $D$. For $s\in I$, project $\gamma(s)$ vertically upward (along the $d$-th coordinate of $B_{(I)}$) onto $\partial (\Omega_{D^{\eps}})\cap D$ to obtain $\widehat \gamma(s)$. Note that $\widehat \gamma$ maybe discontinuous even in the interior of $I$. However, it is continuous on $[0,\tau]$ except possibly for finitely many points. Let $\{0\leq s_1<s_2<\cdots<s_m\leq \tau\}$ be the collection of discontinuities for $\widehat \gamma([0,\tau])$. Then $0<|\widehat \gamma(s_j-)-\widehat \gamma(s_j+)|\leq 2R\eps$ and we can connect $\widehat \gamma(s_j-)$ to $\widehat \gamma(s_j+)$ by a continuous curve $\beta_j:\,[0,1]\longrightarrow \partial (\Omega_{D^{\eps}})\cap D\, \cap B(\widehat \gamma(s_j-),\,8R\eps) \cap B(\widehat \gamma(s_j+),\,8R\eps)$.

Define $\tilde{\gamma}:\,[0,\,\tau+m]\longrightarrow \partial (\Omega_{D^{\eps}})\cap D$ to be the continuous curve obtained by concatenating $\widehat \gamma$ and $\{\beta_j:\,j=1,2,\cdots,m\}$ (See Figure \ref{fig:ExtendDomain_curve}). Then $\tilde{\gamma}(0)=\widehat \gamma(0)=p\in (W_1)_{\Delta_2}$ and $\tilde{\gamma}(\tau+m)=\widehat \gamma(\tau)=\gamma(\tau)\in \partial (\Omega_{D^{\eps}})\cap D\setminus (\partial W_1 \cap D)$. By the  continuity of $\tilde{\gamma}$, there is some $t_{*}\in (0,\,\tau+m)$ such that $\tilde{\gamma}(t_{*})\in (W_1)_{\Delta_1}$. (Roughly speaking, on $\tilde{\partial} D^{\eps}$, $\Delta_2A$ is separated from $\tilde{\partial} D^{\eps}\setminus A$ by $\Delta_1A$.)

Now for some $1\leq j\leq m$, we have $\tilde{\gamma}(t_{*})$ and $\gamma(s_j)$ are connected in $D\setminus \Omega_{D^{\eps}} \subset D\setminus \bar{W_1}$, and
$$|\tilde{\gamma}(t_{*})-\gamma(s_j)|\leq |\tilde{\gamma}(t_{*})-\widehat \gamma(s_j)|+|\widehat \gamma(s_j)-\gamma(s_j)|\leq 8R\eps+R\eps.$$
Hence $\tilde{\gamma}(t_{*})\in \partial U_x$ for some $x\in \Delta_1 A$. We therefore have $\gamma(s_j)\in \Theta_{x}$. The proof is now complete.
\end{pf}

\subsubsection{Discrete relative isoperimetric inequality}

Let $\pi_{srw}$ be the stationary measure of the simple random walk (SRW) on the graph under consideration and recall Definition \ref{Def:IsoperimetricConstant}.

\begin{thm}[Discrete relative isoperimetric inequality]\label{T:Isoperimetric_Discrete}
For every $r\in(0,1)$, there exists $\wh{S}_{srw}=\wh{S}_{srw}(d,D,r)\in (0,\infty)$ and $\eps_1=\eps_1(d,D,r)\in (0,\infty)$ such that
\begin{eqnarray*}
&& \sup_{\eps\in(0, \eps_1)}S_{srw}(r,D^{\eps})  \leq  \wh{S}_{srw}, \text{  and}\\
&& \tilde{S}_{srw}(r,D^{\eps}) \leq  \dfrac{2d\,(m_{srw}(D^{\eps}))^{1/d}}{\eps}\wh{S}_{srw}
\quad \hbox{ for every  }\eps\in(0,\eps_1).
\end{eqnarray*}
\end{thm}

\begin{pf}
We can also assume that $A$ is grid-connected. This is because
$$\dfrac{|A|^{(d-1)/d}}{|\partial_e A|} \leq \dfrac{|A_1|^{(d-1)/d}}{|\partial_e A_1|} \,\vee\,\dfrac{|A_2|^{(d-1)/d}}{|\partial_e A_2|}$$
whenever $A=A_1\cup A_2$ with $\partial_e A_1 \cap \partial_e A_2=\phi$.
From Lemma \ref{L:ExtendDomain} and Proposition \ref{prop:relative_iso_ineq}, we have
$$
\sup_{\eps\in (0, \eps_1)}\sup_{\{A\subset D^{\eps}: \pi(A)\leq r\}}\dfrac{|A|^{(d-1)/d}}{|\partial A|} \leq\, M_2\,M_1^{(d-1/d}\,S\left(\frac{49r+1}{50},\,D\right) .
$$
We thus have the first inequality since $4d|\partial A|\geq |\partial_eA|\geq 2|\partial A|$. The second inequality follows from the first since $q(e)=\frac{(\eps)^d}{2d\,m(D^{\eps})}$.
\end{pf}

\medskip

For the CTRW $X^{aD^{\eps}}$ on $aD^{\eps}$, we let $\pi$ be the stationary measure. Observe that, because $\pi^{aD^{\eps}}(aA) = \pi^{D^{\eps}}(A)$ and $m(aD^{\eps})=a^d\,m(D^{\eps})$, we have
\begin{equation}\label{E:Homo_Isoperimetric}
S_{\pi}(r,\,aD^{\eps})= S_{\pi}(r,\,D^{\eps}) \quad\text{and}\quad  \tilde{S}_{\pi}(r,\,aD^{\eps})= \tilde{S}_{\pi}(r,\,D^{\eps})
\end{equation}
for all $a>0$ and $r>0$. Hence we only need to consider the case $a=1$. In view of Theorem \ref{T:Isoperimetric_Discrete} and (\ref{E:Approximate_m(x)}), we have (taking $r=1/2$)

\begin{cor}\label{cor:Isoperimetric_Discrete}
There exist positive constants $\wh S= \wh{S}(d,D,\rho)$ , $\eps_1=\eps_1(d,D,\rho)$ and $\wh{C}=\wh{C}(d,D,\rho)$ such that
\begin{eqnarray}
&& \sup_{\eps\in(0, \eps_1)}S_{\pi}(1/2,\,D^{\eps})  \leq  \wh{S}, \text{  and} \\
&& \tilde{S}_{\pi}(1/2,\,D^{\eps}) \leq \frac{\wh{C}}{\eps} \wh{S}  \quad
\ \hbox{ for every }\eps\in(0,\eps_1). \label{E:Isoperimetric_Discrete2}
  \label{hat{S}}
\end{eqnarray}
\end{cor}

As an immediate consequence of Corollary \ref{cor:Isoperimetric_Discrete} and (\ref{E:Homo_Isoperimetric}), we have the following Poincar\'e inequality.
\begin{cor}\label{cor:Poincare}(Poincar\'e inequality)
There exist $\eps_1=\eps_1(d,D,\rho)>0$  such that
$$ \frac{|D|\,a^{d-2}}{16\,\wh{C}^2\,\wh{S}^2}\,\|f-\<f\>_{\pi}\|_{l^2(\pi)}^2 \leq \D^{a\eps}_{aD}(f)$$
for all $f\in l^2(aD^{\eps},\pi)$, $\eps\in (0, \eps_1)$, $a>0$. Here $\<f\>_{\pi}:= \sum f\,\pi$, $\wh{C}$ and $\wh{S}$ are the same constants in Corollary \ref{cor:Isoperimetric_Discrete}, and $\D^{a\eps}_{aD}$ is the Dirichlet form in $l^2(m^{aD^{\eps}})$ of the CTRW $X^{aD^{\eps}}$ (see (\ref{E:DirichletForm_CTRW})).
\end{cor}

\begin{pf}
By Corollary \ref{cor:Isoperimetric_Discrete}, the isoperimetric constant
$$\mathcal{I}:= \inf_{\pi(A)\leq 1/2}\dfrac{Q(\partial A)}{\pi(A)} \geq 2^{1/d}\,\frac{1}{\tilde{S}}\geq 2^{1/d}\,\frac{\eps}{\wh{C}\,\wh{S}}.$$
Hence, by the Cheeger's inequality (see  \cite[Lemma 3.3.7]{scL97}),
$$ \inf_{f} \dfrac{\D^{a\eps}_{aD}(f)}{\|f-\<f\>_{\pi}\|_{l^2(\pi)}^2} \geq \dfrac{d\,m(aD^{\eps})}{(a\eps)^2}\,\dfrac{\mathcal{I}^2}{8}\geq \frac{|D|}{16}\,\dfrac{a^{d-2}}{\wh{C}^2\,\wh{S}^2}
$$
when $\eps>0$ is small enough.
\end{pf}

The above Poincar\'e inequality already tells us a positive lower bound for the spectral gap of $X^{aD^{\eps}}$ and hence gives us an estimate for the mixing time. However, we will state a stronger result in Proposition \ref{prop:mixing} in the next subsection.

\subsection{Nash's inequality and Poincar\'e inequality}

The discrete relative isoperimetric inequality leads to the following two functional inequalities;  namely, a Poincar\'e inequality and a Nash inequality that are uniform in $\eps$ and in scaling $D\mapsto aD$. The uniformity in scaling helps proving the near diagonal lower bound for $p^{\eps}$.

\begin{thm}(Nash's inequality and Poincar\'e inequality uniform in $\eps$ and in scaling)
There exist $\eps_1=\eps_1(d,D,\rho)>0$ and $C=C(d,D,\rho)>0$ such that
\begin{equation}\label{E:Nash_a1}
\|f-\<f\>_{\pi}\|_{l^2(\pi)}^{2(1+2/d)}\leq 8\,\tilde{S}_{\pi}(1/2,\,D^{\eps})^2 \,\left(  \dfrac{(a\,\eps)^2}{d\,m(aD^{\eps})}\,\D^{a\eps}_{aD}(f) \right)\,\|f\|_{l^1(\pi)}^{4/d}
\end{equation}
\begin{equation}\label{E:Nash_a2}
\|f\|_{l^2(m)}^{2(1+2/d)}\leq C\,\left(\D^{a\eps}_{aD}(f) + \,(\wh{C}\,\wh{S}\,a)^{-2}\,\|f\|_{l^2(m)}^2 \right)\,\|f\|_{l^1(m)}^{4/d}
\end{equation}
for every $ f\in l^2(aD^{\eps})$, $\eps\in (0,\eps_1)$ and $a\in(0,\infty)$, where $\wh{C}$ and $\wh{S}$ are the same constants in Corollary \ref{cor:Isoperimetric_Discrete}; $\<f\>_{\pi}:= \sum f\,\pi$ and $\D^{a\eps}_{aD}$ is the Dirichlet form in $l^2(m^{aD^{\eps}})$ of the CTRW $X^{aD^{\eps}}$ (see (\ref{E:DirichletForm_CTRW})).
\end{thm}

\begin{pf}
Note that $ \dfrac{(a\,\epsilon)^2}{d\,m(aD^{\epsilon})}\,\D^{a\epsilon}_{aD}(f)$ is the Dirichlet form of the unit speed CTRW with the same one-step transition probabilities as that of $X^{aD^{\epsilon}}$. Hence (\ref{E:Nash_a1}) follows directly from \cite[Theorem 3.3.11]{scL97} and (\ref{hat{S}}). For (\ref{E:Nash_a2}), let $R=(2^{\frac{1}{d}}\delta\,a\,\eps)^{-1}$ with $\delta \geq  (\wh{C}\,\wh{S}\,a)^{-1}$.
For any nonempty subset $A\subset aD^{\eps}$,
\begin{eqnarray*}
\dfrac{Q(\partial A)+\frac{1}{R}\pi(A)}{\pi(A)^{\frac{d-1}{d}}} &\geq& \dfrac{1}{\wh{S}} \wedge \dfrac{1}{R}\left(\dfrac{1}{2}\right)^{\frac{1}{d}}
\geq  \left((\wh{C}\,\wh{S})^{-1}\wedge a\,\delta \right)\,\eps
= (\wh{C}\,\wh{S})^{-1}\,\eps.
\end{eqnarray*}
Hence,
\begin{equation}\label{E:NashAssumption}
\sup_{A\subset aD^{\eps}}\dfrac{\pi(A)^{\frac{d-1}{d}}}{Q(\partial A)+\frac{1}{R}\pi(A)} \leq
\dfrac{\wh{C}\,\wh{S}}{\eps}.
\end{equation}
By \cite[Theorem 3.3.10]{scL97},
$$\|f\|_{l^2(\pi)}^{2(1+2/d)}\leq 16\,(\dfrac{C}{\eps})^2\left( \dfrac{(a\,\eps)^2}{d\,m(aD^{\eps})}\,\D^{a\eps}_{aD}(f) + \dfrac{1}{8R^2}\|f\|_{l^2(\pi)}^2 \right)\,\|f\|_{l^1(\pi)}^{4/d}.$$

Using the relations $\|f\|_{l^2(\pi)}^2=(m(aD^{\eps}))^{-1}\|f\|_{l^2(m)}^2$, $\|g\|_{l^1(\pi)}=(m(aD^{\eps}))^{-1}\|f\|_{l^1(m)}$ and (\ref{E:Approximate_m(x)}), we get the desired inequality (\ref{E:Nash_a2}).
\end{pf}

\subsection{Mixing time}

By the Poincar\'e inequalities in (\ref{E:Nash_a1}) and \cite[Corollary 2.3.2]{scL97}, we obtain an estimate on the time needed to reach stationarity.
\begin{prop}\label{prop:mixing}(Mixing time estimate)
There exists $C>0$ which depends only on $d$ such that
$$\Big|p^{a\eps}_{aD}\left(t, x, y \right) - \dfrac{1}{ m(aD^{\eps})\,}\Big| \leq C\,\min\bigg\{\,(a\,\wh{C}\,\wh{S})^d\,{t}^{-d/2},\; \frac{1}{(a\,\eps)^d}\,\exp{\Big(\frac{-d\,t}{8\,(a\wh{C}\wh{S})^2}\Big)} \,\bigg\}  $$
for every $t>0$, $ x,\,y\in aD^{\eps}$, $\eps\in (0, \eps_1)$ and  $a>0$.
Here $\wh{C}$ and $\wh{S}$ are the  constants in Corollary \ref{cor:Isoperimetric_Discrete}.
\end{prop}

\begin{pf}
By (\ref{E:Nash_a1}) and Theorem 2.3.1 of \cite{scL97}, we have
$$\Big| m(aD^{\eps})\,p^{a\eps}_{aD}\left(\frac{(a\eps)^2}{d}t, x, y \right) - 1\Big| \leq \left(\dfrac{d\,(8\,\tilde{S}^2)}{2t} \right)^{d/2}$$
After simplification and using (\ref{E:Isoperimetric_Discrete2}), we obtain the upper bound which is of order ${t}^{-d/2}$. On other hand, by Corollary \ref{cor:Poincare} and \cite[Lemma 2.1.4]{scL97}, we obtain the exponential term on the right hand side.
\end{pf}

\subsection{Gaussian bound and uniform H\"older continuity of $p^{\eps}$}

Equipped with the Nash inequality (\ref{E:Nash_a2}) and the Poincar\'e inequality (\ref{E:Nash_a1}), one can follow a now standard procedure (see, for example, \cite{CKS86} or \cite{tD99}) to obtain two sided Gaussian estimates for $p^{\eps}$. In the following, $C_1$, $C_2$ and $\eps_0$ are positive constants which depends only on $d,\,D,\,\rho$ and $T$.

More precisely, we only need the Nash inequality (\ref{E:Nash_a2}) and Davies' method to obtain the following Gaussian upper bound.
\begin{thm}\label{T:UpperHKE_a}
There exist constants $C_i=C_i(d,D,\rho,T)>0$, $i=1,2$,  and $\eps_0=\eps_0(d,D,\rho,T) \in (0, 1]$ such that
\begin{equation*}
p^{a\eps}_{aD}(t,x,y) \leq \dfrac{C_1}{(a\eps\vee t^{1/2})^d}\,\exp\left(\frac{C_2}{a^2}t-\frac{|y-x|^2}{(a\eps)^2\vee t}\right)
\end{equation*}
for every  $t\geq a\eps$, $\eps\in(0,\eps_0)$, $a>0$ and $x,y \in aD^{\eps}$.
Moreover, the following weaker bound holds for all $t>0$:
\begin{equation*}
p^{a\eps}_{aD}(t,x,y) \leq \dfrac{C_1}{(a\eps\vee t^{1/2})^d}\,\exp\left(\frac{C_2}{a^2}t-\frac{|y-x|}{a\eps\vee t^{1/2}}\right).
\end{equation*}
In particular, this implies the upper bound in Theorem \ref{T:UpperHKE} which is the case when $a=1$.
\end{thm}
We can then apply the Poincar\'e inequality (\ref{E:Nash_a1}) and argue as in section 3 of \cite{tD99} to obtain the near diagonal lower bound. A more comprehensive proof is given in \cite{wtF14}.
\begin{lem}\label{L:NearDiagonalLowerBound}
\begin{equation*}
p^{\eps}(t,x,y) \geq \dfrac{C_2}{(\eps\vee t^{1/2})^d}
\end{equation*}
for every $(t,x,y)\in(0,\infty)\times D^{\eps}\times D^{\eps}$ with $|x-y|\leq C_1\,t^{1/2}$ and $\eps\in (0,\eps_0)$.
\end{lem}
The Gaussian lower bound for $p^{\eps}$ in Theorem \ref{T:LowerHKE} then follows from the Lipschitz property of $D$ and a well-known chaining argument (see, for example, page 329 of \cite{dS88}). Therefore, we  have  the two-sided Gaussian bound for $p^{\eps}$ as stated in Theorem \ref{T:UpperHKE} and Theorem \ref{T:LowerHKE}. It then  follows from a standard `oscillation' argument (cf. Theorem 1.31 in \cite{SZ97} or Theorem II.1.8 in \cite{dS88}) that $p^{\eps}$ is H\"older continuous in $(t,x,y)$, \emph{uniformly in $\eps$}. More precisely,
\begin{thm}\label{T:HolderCts}
There exist positive constants $\gamma=\gamma (d,D,\rho)$, $\eps_0(d,D,\rho)$ and $C(d,D,\rho)$ such that for all $\eps\in(0,\eps_0)$, we have
\begin{equation}\label{E:HolderCts2}
|p^{\eps}(t,x,y)-p^{\eps}(t',x',y')| \leq C\, \dfrac{ (|t-t'|^{1/2}+\|x-x'\|+\|y-y'\|)^\gamma }
{ (t\wedge t')^{\sigma/2}\,[1 \wedge (t\wedge t')^{d/2}] } .
\end{equation}
\end{thm}

\subsection{Proof of local CLT}

The following weak convergence result for RBM with drift is a natural generalization of \cite[Theorem 3.3]{BC08}.
\begin{thm}\label{T:WeakConvergence_RBMDrift}
Let $D\subset \R^d$ be a bounded domain whose boundary $\partial D$ has zero Lebesque measure. Suppose $D$ also satisfies:
$$C^{1}(\bar{D}) \text{  is dense in  } W^{1,2}(D).$$
Suppose $\rho\in W^{1,2}(D)\cap C^1(\bar{D})$ is strictly positive. Then for every $T>0$, as $k\rightarrow \infty$,
\begin{enumerate}
\item[(i)]  $(X^{2^{-k}},\,\P_{m})$ converges weakly to the stationary process $(X,\,\P_{\rho})$ in the Skorokhod space $D([0,T],\bar{D})$.
\item[(ii)]  $(X^{2^{-k}},\,\P_{x_k})$ converges weakly to $(X,\,\P_{x})$ in the Skorokhod space $D([0,T],\bar{D})$ whenever $x_k$ converges to $x\in D$.
\end{enumerate}
\end{thm}

\begin{pf}
For (i), the proof follows from a direct modification of the proof of \cite[Theorem 3.3]{BC08}.
Recall the definition of the one-step transition probabilities $p_{xy}$, defined in the paragraph that contains (\ref{E:Conductance_BaisedRW+}) and (\ref{E:Conductance_BaisedRW-}).
Observe that, since $\rho\in C^1(\bar{D})$, approximations using Taylor's expansions  in the proofs of \cite[Lemma 3.1 and Lemma 2.2]{BC08} continue to work  with the current definition of $p_{xy}$. Thus we have
\begin{eqnarray*}
\lim_{k\to\infty}\,\D^{2^{-k}}(f,f) &=& \frac{1}{2}\,\int_{D}|\nabla f(x)|^2\,\rho(x)dx, \quad \forall\,f\in C^{1}(\bar{D}), \quad \text{and}\\
\lim_{k\to\infty}\,L^{(2^{-k})}f &=& \frac{1}{2}\,\Delta f + \frac{1}{2}\nabla (\log \rho)\cdot\nabla f \quad \text{uniformly in } D,\quad  \forall\,f\in C^{\infty}_{c}(D) .
\end{eqnarray*}
The process $X^{\eps}$ has a L\'evy system $(N^{\eps}(x,dy),\,t)$, where for $x\in D^{\eps}$,
$$N^{\eps}(x,dy)= \frac{d}{\eps^2}\,\sum_{z:\,z\leftrightarrow x}\,p_{xz}\delta_{\{z\}}(dy).$$
Following the same calculations as in the proof of \cite[Theorem 3.3]{BC08}, while noting that \cite[Theorem 6.6.9]{CF12} (in place of \cite[Theorem 1.1]{BC08}) can be applied to handle general symmetric reflected diffusions as in our present case, we get part (i). Part (ii) follows from part (i) by a localization argument (cf. \cite[Remark 3.7]{BC11}).
\end{pf}

\medskip

We can now present the  proof of  the local CLT.

\medskip

\noindent{\it Proof of Theorem \ref{T:LCLT_CTRW}}.
For each $\eps>0$ and $t>0$, we extend $p^{\eps}(t,\cdot,\cdot)$ to $\bar{D}\times \bar{D}$ in such a way that $p^{\eps}$ is nonnegative and continuous on $(0,\infty)\times\bar{D}\times\bar{D}$, and that both the maximum and the minimum values are preserved on each cell in the grid $\eps \Z^d$.
This can be done in many ways, say by the interpolation described in \cite{BK08}, or a sequence of harmonic extensions along the simplexes (described in  \cite{wtF14}).

Consider the family $\{t^{d/2}p^\eps\}_{\eps}$ of continuous functions on $(0, \infty)\times \bar{D} \times \bar{D}$. Theorem \ref{T:UpperHKE} and Theorem \ref{T:HolderCts} give us uniform pointwise bound and equi-continuity respectively. By Arzela-Ascoli Theorem, it is relatively compact. i.e. for any sequence $\{\eps_n\}\subset (0,1]$ which decreases to $0$, there is a subsequence $\{\eps_{n'}\}$ and a continuous $q:\,(0,\infty)\times \bar{D}\times \bar{D}\longrightarrow [0,\infty)$ such that $p^{\eps_{n'}}$ converges to $q$ locally uniformly.

On other hand, by part (ii) of Theorem \ref{T:WeakConvergence_RBMDrift}, if the original sequence $\{\eps_k\}$ is a subsequence of $\{2^{-k}\}$, then $q=p$. More precisely, the weak convergence implies that for all $t>0$,
$$
\int_{D} \phi(y)p(t,x,y)dy=\int_{D} \phi(y)q(t,x,y)dy \quad \text{for all }\phi\in C_c(D) \text{ and } x\in D.
$$
Then by the continuity of both $p$ and $q$ in the second coordinate, we have $q=p$ on $(0,\infty)\times D\times \bar{D}$. Since $p(t,\cdot,\cdot)$ and $q(t,\cdot,\cdot)$ are continuous on $\bar{D}\times \bar{D}$ (cf. \cite{BH91}), we obtain $p=q$ on $(0,\infty)\times \bar{D}\times \bar{D}$. In conclusion, we have $p^{\eps}$ converges to $p$ locally uniformly through the sequence $\{\eps_n=2^{-n}; n\geq 1\}$.
\qed

\bigskip

{\bf Acknowledgements}:  We thank Krzysztof Burdzy, Rekhe Thomas and Tatiana Toro for helpful discussions.
In particular, we are grateful to David Speyer for pointing out the recurrence relation  \eqref{E:Recursion_J_N(t)} to us.
We also thank Guozhong Cao, Samson A. Jenekhe, Christine Luscombe, Oleg Prezhdo and Rudy Schlaf
for  discussions on solar cells. Financial support from NSF Solar Energy Initiative grant DMR-1035196 as well as NSF grant DMS-1206276
is gratefully acknowledged.

\vspace{5mm}
    \textbf{Zhen-Qing Chen}

    Department of Mathematics, University of Washington, Seattle, WA 98195, USA

    Email: zqchen@uw.edu
\vspace{2mm}

    \textbf{Wai-Tong (Louis) Fan}

    Department of Mathematics, University of Washington, Seattle, WA 98195, USA

    Email: louisfan@math.washington.edu


\begin{thebibliography}{99}

  \bibitem{BH91}
    R.F. Bass and P. Hsu.
    Some potential theory for reflecting Brownian motion in H\"older and lipschitz domains.
    \emph{Ann. Probab.} \textbf{19} (1991), 486-508.
  \bibitem{BK08}
    R.F. Bass and T. Kumagai.
    Symmetric Markov chains on $Z^d$ with unbounded range.
    \emph{Trans. Amer. Math. Soc.} \textbf{360} (2008), 2041-2075.
  \bibitem{djB91}
    D.J. Blount.
    Comparison of stochastic and deterministic models of a linear chemical reaction with diffusion.
    \emph{Ann. Probab.} \textbf{19} (1991), 1440-1462.
  \bibitem{BDP92}
    C. Boldrighini, A. De Masi and A. Pellegrinotti.
    Nonequilibrium fluctuations in particle systems modelling reaction-diffusion equations.
    \emph{Stoch. Proc. Appl.} \textbf{42} (1992), 1-30.
  \bibitem{BDPP87}
    C. Boldrighini, A. De Masi, A. Pellegrinotti and E. Presutti.
    Collective phenomena in interacting particle systems.
    \emph{Stoch. Proc. Appl.} \textbf{25} (1987), 137-152.
  \bibitem{BC08}
    K. Burdzy and Z.-Q. Chen.
    Discrete approximations to reflected Brownian motion.
    \emph{Ann. Probab.} \textbf{36} (2008), 698-727.
  \bibitem{BC11}
    K. Burdzy and Z.-Q. Chen.
    Reflected random walk in fractal domains.
    \emph{Ann. Probab.} \textbf{41} (2011), 2791-2819.
   \bibitem{BHM00}
    K. Burdzy, R. Holyst and P. March.
    A Fleming-Viot particle representation of Dirichlet Laplacian.
    \emph{Comm. Math. Phys.} \textbf{214} (2000), 679-703.
  \bibitem{BQ06}
    K. Burdzy and J. Quastel.
    An annihilating-branching particle model for the heat equation with average temperature zero.
    \emph{Ann. Probab.} \textbf{34} (2006), 2382-2405.
  \bibitem{CKS86}
    E.A. Carlen, S. Kusuoka and D.W. Stroock.
    Upper bounds for symmetric markov transition functions.
    \emph{Ann. Inst. Henri Poincar\'{e}-Probab. Statist. 1987.} \textbf{23} (1986), 245-287. MR 898496.
  \bibitem{mfChen03}
    M.-F. Chen.
    \emph{From Markov Chains to Non-Equilibrium Particle Systems. 2nd ed.}
    World Scientific, 2003.
  \bibitem{zqChen93}
    Z.-Q. Chen.
    On reflecting diffusion processes and Skorokhod decompositions.
    \emph{Probab. Theory Relat. Fields} \textbf{94} (1993), 281-316.

  \bibitem{CF12}
    Z.-Q. Chen and M. Fukushima.
    \emph{Symmetric Markov Processes, Time Change and Boundary Theory.}
    Princeton. University Press, 2012.

  \bibitem{CDP11}
    J.T. Cox, R. Durrett and E. Perkins.
    Voter model perturbations and reaction diffusion equations.
    \emph{Ast\'erisque} \textbf{349},  2013.

  \bibitem{gD88}
    G. David.
    Morceaux de graphes lipschitziens et intégrales singulières sur une surface.
    \emph{Rev. Mat. Iberoamericana.} \textbf{4} (1988), 73–114 (French). MR 1009120.
  \bibitem{DS91}
    G. David and S. Semmes.
    Singular integrals and rectifiable sets in $\R^n$: Beyond Lipschitz graphs.
    \emph{Ast\'erisque.} \textbf{193} (1991), 152. MR 1113517 (92j:42016)9120.

  \bibitem{tD99}
    T. Delmotte.
    Parabolic Harnack inequality and estimates of Markov chains on graphs.
    \emph{Revista Matem\'atica Iberoamericana.} \textbf{15}(1) (1999), 181-232.


  \bibitem{pD88a}
    P. Dittrich.
    A stochastic model of a chemical reaction with diffusion.
    \emph{Probab. Theory Relat. Fields.} \textbf{79} (1988), 115-128.
  \bibitem{pD88b}
    P. Dittrich.
    A stochastic partical system: Fluctuations around a nonlinear reaction-diffusion equation.
    \emph{Stochastic Processes. Appl.} \textbf{30} (1988), 149-164.

  \bibitem{rDsL94}
    R. Durrett and S. Levin.
    The importance of being discrete (and spatial).
    \emph{Theoretical population biology} \textbf{46.3} (1994), 363-394.

  \bibitem{lEbShtY07}
    L. Erd\"os, B. Schlein and H. T. Yau.
    Derivation of the cubic non-linear Schrödinger equation from quantum dynamics of many-body systems.
    \emph{Inventiones Mathematicae} \textbf{167}(3) (2007), 515-614.
  \bibitem{EK86}
    S.N. Ethier and T.G. Kurtz.
    \emph{Markov processes. Characterization and Convergence.}
    Wiley, New York, 1986. MR0838085.
  \bibitem{wtF14}
    W.-T. Fan.
    Systems of reflected diffusions with interactions through membranes. PhD thesis. In preparation, 2014.
  \bibitem{fG05}
    F. Golse.
    Hydrodynamic limits.
    \emph{Euro. Math. Soc.} \textbf{1} (2005), 699-717.
  \bibitem{GPV88}
    M.Z. Guo, G.C. Papanicolaou and S.R.S. Varadhan.
    Nonlinear diffusion limit for a system with nearest neighbor interactions.
    \emph{Comm. Math. Phys.}. \textbf{118} (1988), 31-59.
  \bibitem{GSC11}
    P. Gyrya and L. Saloff-Coste.
    \emph{Neumann and Dirichlet Heat Kernels in Inner Uniform Domains.}
    Paris: Société mathématique de France, 2011.
  \bibitem{cKcL98}
    C. Kipnis and C. Landim.
    \emph{Scaling Limits of Interacting Particle Systems.}
    Springer, 1998.
  \bibitem{KOV89}
    C. Kipnis, S. Olla and S.R.S. Varahan.
    Hydrodynamics and large deviations for simple exclusion process.
   \emph{Comm. Pure Appl. Math.} \textbf{42} (1989), 115-137.
  \bibitem{pK86}
    P. Kotelenez.
    Law of large numbers and central limit theorem for linear chemical reactions with diffusion.
    \emph{Ann. Probab.} \textbf{14} (1986), 173-193
  \bibitem{pK88}
    P. Kotelenez.
    High density limit theorems for nonlinear chemical reactions with diffusion.
    \emph{Probab. Theory Relat. Fields} \textbf{78} (1988), 11-37.
  \bibitem{tK71}
    T. Kurtz.
    Limit theorems for sequences of jump Markov processes approximating ordinary differential processes.
    \emph{J. Appl. Probab.} \textbf{8} (1971), 344-356.
   \bibitem{tK81}
    T. Kurtz.
    \emph{Approximation of Population Processes.}
    SIAM, Philadelphia, 1981.
  \bibitem{LX80}
    R. Lang and N.X. Xanh.
    Smoluchowski's theory of coagulation in colloids holds rigorously in the Boltzmann-Grad-limit.
   \emph{Zeitschrift f\"{u}r Wahrscheinlichkeitstheorie und Verwandte Gebiete.} \textbf{54} (1980), 227-280.
  \bibitem{rMmN92}
    R. M. May and M. A. Nowak.
    Evolutionary games and spatial chaos.
    \emph{Nature} \textbf{359.6398} (1992), 826-829.
  \bibitem{vgM85}
    V.G. Mazja.
    \emph{Sobolev Spaces.}
    Springer Verlag, 1985
  \bibitem{scL97}
    L. Saloff-Coste.
    \emph{Lectures on Finite Markov Chains.}
    Lecture Notes in Math., Springer, 1997.
  \bibitem{rpS99}
    R.P. Stanley.
    \emph{Enumerative Combinatorics II.}
    Cambridge University Press, Cambridge, 1999.
  \bibitem{dS88}
    D.W. Stroock.
    Diffusion semigroups corresponding to uniformly elliptic divergence form operators.
    \emph{S\'eminaire de Probabilit\'es XXII.}
    Lecture Notes in Math., vol. 1321, pp. 316-347, Springer, Berlin, 1988.
  \bibitem{SZ97}
    D.W. Stroock and W. Zheng.
    Markov chain approximations to symmetric diffusions.
    \emph{Ann. Inst. Henri. Poincaré-Probab. Statist.} \textbf{33} (1997), 619-649. MR 1473568.

  \bibitem{htY91}
    H.T. Yau.
    Relative entropy and hydrodynamics of Ginzburg-Landau models.
    \emph{Lett. Math. Phys}. \textbf{22} (1991), 63-80.
\end{thebibliography}
\end{document}